\documentclass[11pt, a4paper]{amsart}

\usepackage{amssymb}
\usepackage{mathrsfs}
\usepackage{graphicx}
\usepackage{hyperref}

\hypersetup{
unicode=true,
colorlinks=true,
linkcolor=blue,
citecolor=blue,
urlcolor=blue,
filecolor=blue,
bookmarksnumbered=true,
pdfstartview=FitH,
pdfhighlight=/N
}

\theoremstyle{definition}

\newtheorem{remark}{Remark}

\def\({\left(}
\def\){\right)}
\def\supp{\operatorname{supp}}
\def\NN{\mathbb N}
\def\RR{\mathbb R}
\def\CC{\mathbb C}
\def\PP{\mathbb P}
\def\HH{\mathscr H}
\def\RS{\mathfrak R}
\def\ZZ{\mathbb Z}
\def\zz{\mathbf z}
\def\mb{\beta}
\let\leq\leqslant
\let\geq\geqslant

\let\ge\geqslant
\let\myh\widehat
\let\myt\widetilde
\let\myo\overline

\begin{document}


\title{Some numerical results on the behavior of
zeros of the Hermite--Pad\'e polynomials}

\author{N.~R. Ikonomov}
\address{Institute of Mathematics and Informatics, Bulgarian Academy of Sciences}
\email{nikonomov@math.bas.bg}
\author{R.~K. Kovacheva}
\address{Institute of Mathematics and Informatics, Bulgarian Academy of Sciences}
\email{rkovach@math.bas.bg}
\author{S.~P. Suetin}
\address{Steklov Mathematical Institute, Russian Academy of Sciences}
\email{suetin@mi.ras.ru}
\thanks{The work of S.~P.~Suetin is supported by the
\href{http://www.rscf.ru/contests}{Russian Science Foundation (RScF)}
under a grant 14-50-00005.}

\date{Januray 28, 2015}

\maketitle

\markright{ZEROS OF THE HERMITE--PAD\'E POLYNOMIALS}


\begin{abstract}
We introduce and analyze some numerical results obtained by
the authors experimentally.
These experiments are related to the well known
problem about the distribution of the zeros of Hermite--Pad\'{e} polynomials
for a collection of three functions $[f_0\equiv1,f_1,f_2]$.
The numerical results refer to two cases:
a pair of functions $f_1,f_2$ forms an Angelesco system  (see \eqref{3})
and a pair of functions $f_1=f,f_2=f^2$
forms a (generalized) Nikishin system (see \eqref{4}).
The authors hope that the obtained numerical results
will set up a new conjectures about the limiting
distribution of the zeros of Hermite--Pad\'{e} polynomials.

Bibliography: \cite{Sue13b}~titles; \ref{Fig_bus205c_4000_120_full} pictures.
\end{abstract}


\setcounter{tocdepth}{1}
\tableofcontents

\section{Introduction}\label{s1}

\subsection{}\label{s1s1}
In the present work we introduce and analyze some numerical results
obtained by experimenting. These experiments are connected with the well known
problem about the distribution of the zeros of Hermite--Pad\'{e} polynomials
for a collection of three functions $[f_0\equiv1,f_1,f_2]$, defined and holomorphic
at the infinity point $z=\infty$, $f_1,f_2\in\HH(\infty)$.
Our numerical experiments are restricted only to the Hermite--Pad\'{e} polynomials
of first kind. It is well known that Hermite--Pad\'{e} polynomials\footnote{Sometimes
these polynomials are called ``multiple orthogonal polynomials'',
see \cite{AsFi13}, \cite{Ass11}.}
of first and second kind
are closely related to each other, see \cite[\S 2, equation (2.1.9)]{Nut84}, \cite{Ass06},
\cite{AsFi13}, \cite{Ass11}. The numerical results obtained here
may permit interpretations about the Hermite--Pad\'{e} polynomials of the second kind.

Let
$\PP_n:=\CC_n[z]$, $n\in\NN$,
be the class of all polynomials with complex coefficients of degree $\leq{n}$.
For an arbitrary $n\in\NN$, define the Hermite--Pad\'{e} polynomials of first kind
$Q_{n,j}\in\PP_n$, $j=0,1,2$, $Q_{n,1},Q_{n,2}\not\equiv0$,
by the relation (see \cite{Nut84},
also \cite{GoRa81}, \cite{Sta88}, \cite{Gon03}, \cite{Apt08}, \cite{ApKu11},
\cite{Apt12}, \cite{AsFi13})
\begin{equation}
(Q_{n,0}\cdot1+Q_{n,1}f_1+Q_{n,2}f_2)(z)=O\(\frac1{z^{2n+2}}\),
\qquad z\to\infty.
\label{1}
\end{equation}
In the present work, we will use the terminology and notation
of the work by J. Nuttall \cite{Nut84} (see also H. Stahl \cite{Sta88}).
According to their paper, the classical Pad\'{e} polynomials $P_{n,0},P_{n,1}$ of the function
$f\in\HH(\infty)$ are, in fact, the Hermite--Pad\'{e} polynomials for the collection of two functions
$[f_0\equiv1,f_1=f]$, such that $P_{n,0},P_{n,1}\in\PP_n$, $P_{n,1}\not\equiv0$ and:
\begin{equation}
(P_{n,0}\cdot1+P_{n,1}f)(z)=O\(\frac1{z^{n+1}}\),
\qquad z\to\infty.
\label{2}
\end{equation}

We present the results from our numerical experiments about two cases.

In the first case the functions $f_1$ and $f_2$ have the following form:
\begin{equation}
f_1(z)=\frac{z}{\sqrt{(z-a_1)(z-b_1)}},\qquad
f_2(z)=\frac{z}{\sqrt{(z-a_2)(z-b_2)}},
\label{3}
\end{equation}
where $a_1,b_1,a_2,b_2\in\CC$, $a_1\neq b_1$, $a_2\neq b_2$, and
$\{a_1,b_1\}\cap\{a_2,b_2\}=\varnothing$. Therefore, the pair of functions
$f_1,f_2$ forms an {\it Angelesco system}
(see \cite{Kal79}, \cite{GoRa81}, \cite{Sta88}, \cite{GoRaSo97}, \cite{Gon03},
\cite{Apt08}).

In the second case, $f_1=f$, $f_2=f^2$, where
\begin{equation}
f(z)=(z^2-1)^{1/4}(z-a)^{-1/2},\quad f(\infty)=1,\quad a\not\in\RR,
\label{4}
\end{equation}
and the pair of functions $f_1,f_2$ forms a (generalized) {\it Nikishin system}
(see \cite{Nik80}, \cite{GoRaSo97},
\cite{Gon03}, \cite{Apt08}, \cite{ApLy10}, \cite{Rak11},
\cite{FiLo11}, \cite{DeLoLo12}).

\begin{figure}[!ht]
\centerline{
\includegraphics[width=15cm,height=15cm]{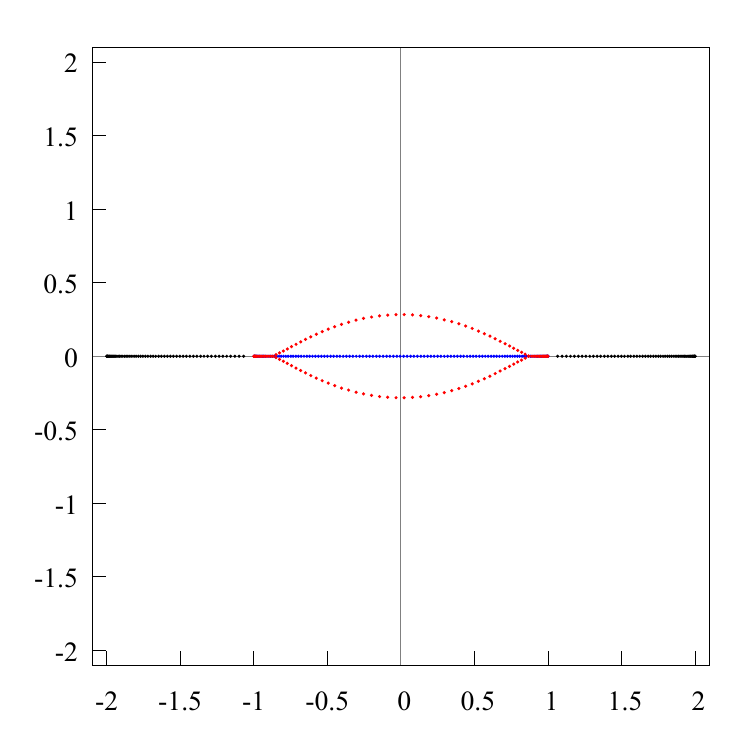}}
\vskip-6mm
\caption{The distribution of the zeros of the Hermite-Pad\'{e} polynomials $Q_{120,0}$ (blue points),
$Q_{120,1}$ (red points) and $Q_{120,2}$ (black points) for a set of three functions
$[1,f_1,f_2]$, where
$f_1(z)=\sqrt{(z-1)/(z+1)}$, $f_2(z)=\sqrt{(z-2)/(z+2)}$.
The sets of branch points of the functions $f_1$ and $f_2$ do not intersect each other,
and thus the pair $f_1,f_2$ create an Angelesco system.
However, the distribution of the zeros of the Hermite-Pad\'{e} polynomials
for this system is the same as the distribution for the Nikishin system
$f_1(z)=\((z-1)/(z+1)\)^{1/3}\((z-2)/(z+2)\)^{1/3}$,
$f_2(z)=\((z-1)/(z+1)\)^{2/3}\((z-2)/(z+2)\)^{1/3}$;
see fig. \ref{Fig_nik(2_1)_4000_120_full}--\ref{Fig_nik(2_1)_4000_120_blrd}.
}
\label{Fig_ang(2_2)_4000_120_full}
\end{figure}

\begin{figure}[!ht]
\centerline{
\includegraphics[width=15cm,height=15cm]{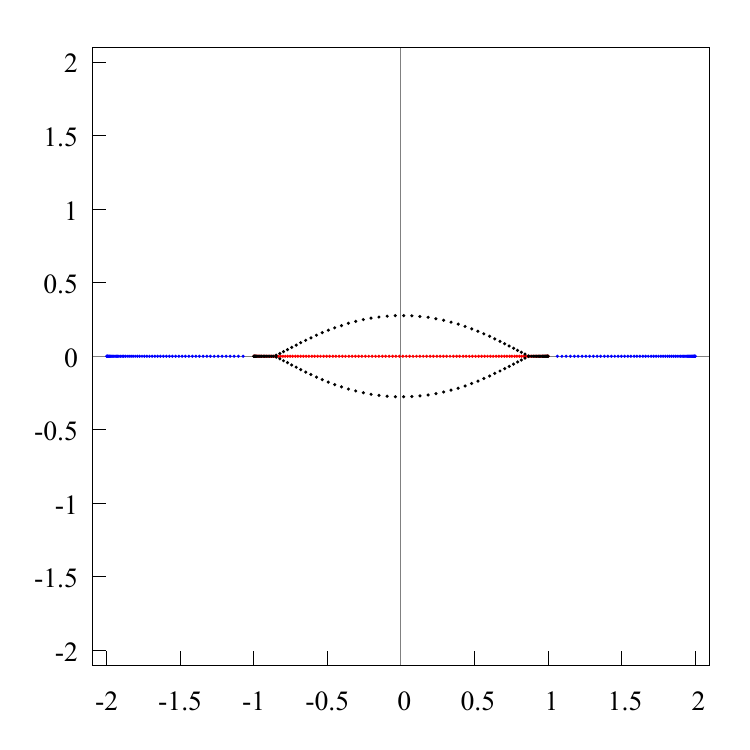}}
\vskip-6mm
\caption{The distribution of the zeros of the Hermite-Pad\'{e} polynomials $Q_{120,0}$ (blue points),
$Q_{120,1}$ (red points) and $Q_{120,2}$ (black points) for a set of three functions
$[1,f_1,f_2]$, where
$f_1(z)=\((z-1)/(z+1)\)^{1/3}\((z-2)/(z+2)\)^{1/3}$,
$f_2(z)=\((z-1)/(z+1)\)^{2/3}\((z-2)/(z+2)\)^{1/3}$.
The pair $f_1,f_2$ create a Nikishin system,
because the branch points of the functions are equivalent.
However, the distribution of the zeros of the Hermite-Pad\'{e} polynomials
for this system is the same as the distribution of the zeros for the Angelesco system
$f_1(z)=\sqrt{(z-1)/(z+1)}$, $f_2(z)=\sqrt{(z-2)/(z+2)}$;
see fig. \ref{Fig_ang(2_2)_4000_120_full},
\ref{Fig_ang(2_2)_4000_120_blbk},
\ref{Fig_nik(2_1)_4000_120_blrd}.
}
\label{Fig_nik(2_1)_4000_120_full}
\end{figure}

\begin{figure}[!ht]
\centerline{
\includegraphics[width=15cm,height=15cm]{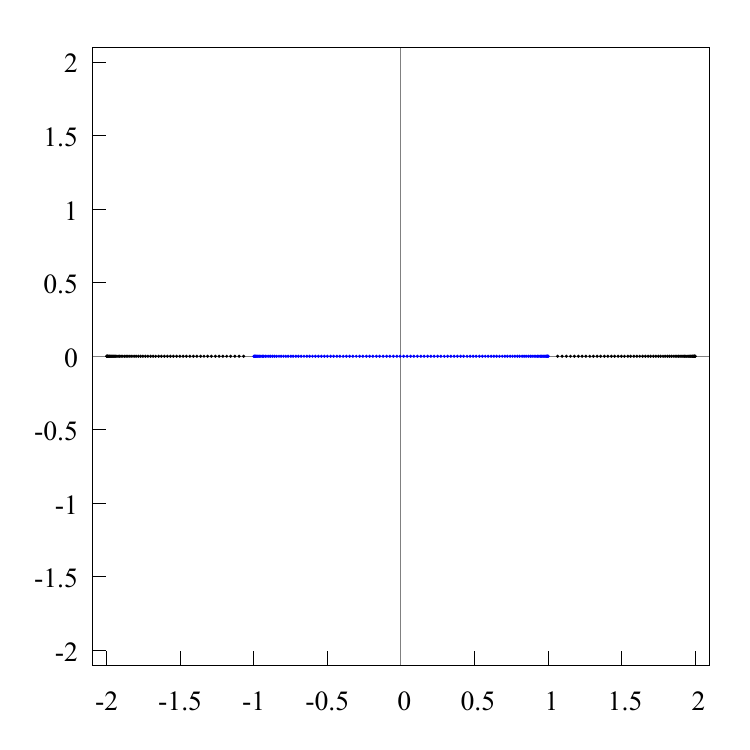}}
\vskip-6mm
\caption{The distribution of the zeros of the Hermite-Pad\'{e} polynomials $Q_{120,0}$ (blue points)
and $Q_{120,2}$ (black points) for a set of three functions
$[1,f_1,f_2]$, where
$f_1(z)=\sqrt{(z-1)/(z+1)}$, $f_2(z)=\sqrt{(z-2)/(z+2)}$.
The pair $f_1,f_2$ create an Angelesco system.
However, the distribution of the zeros of the Hermite-Pad\'{e} polynomials
for this system is the same as the distribution for the Nikishin system
$f_1(z)=\((z-1)/(z+1)\)^{1/3}\((z-2)/(z+2)\)^{1/3}$,
$f_2(z)=\((z-1)/(z+1)\)^{2/3}\((z-2)/(z+2)\)^{1/3}$;
comp. fig. \ref{Fig_nik(2_1)_4000_120_blrd}.
}
\label{Fig_ang(2_2)_4000_120_blbk}
\end{figure}

\begin{figure}[!ht]
\centerline{
\includegraphics[width=15cm,height=15cm]{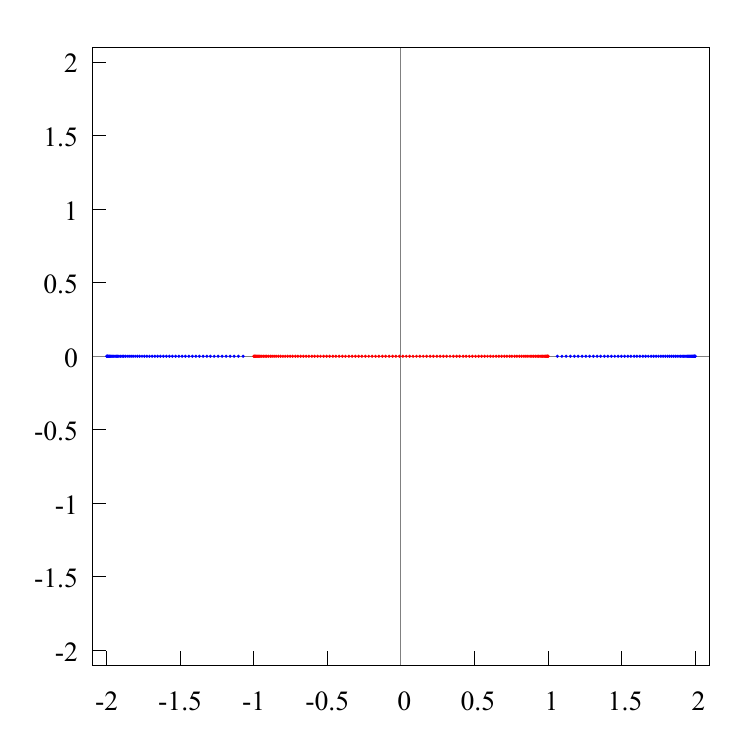}}
\vskip-6mm
\caption{The distribution of the zeros of the Hermite-Pad\'{e} polynomials $Q_{120,0}$ (blue points)
and $Q_{120,1}$ (red points) for a set of three functions
$[1,f_1,f_2]$, where
$f_1(z)=\((z-1)/(z+1)\)^{1/3}\((z-2)/(z+2)\)^{1/3}$,
$f_2(z)=\((z-1)/(z+1)\)^{2/3}\((z-2)/(z+2)\)^{1/3}$.
The pair $f_1,f_2$ create a Nikishin system.
However, the distribution of the zeros of the Hermite-Pad\'{e} polynomials
for this system is the same as the distribution for the Angelesco system
$f_1(z)=\sqrt{(z-1)/(z+1)}$, $f_2(z)=\sqrt{(z-2)/(z+2)}$;
comp. fig. \ref{Fig_ang(2_2)_4000_120_blbk}.
}
\label{Fig_nik(2_1)_4000_120_blrd}
\end{figure}

\begin{remark}\label{rem1}
In the theory of Hermite--Pad\'{e} polynomials for a collection of three
functions $[1,f_1,f_2]$
two opposite situations are usually selected, which are connected
with the distribution of the branch points of the functions $f_1$ and $f_2$.
In the first case the sets of branch points $\mathcal{A}_1=\mathcal{A}(f_1)$ and
$\mathcal{A}_2=\mathcal{A}(f_2)$ {\it do not intersect each other}.
We say that the pair of functions forms an {\it Angelesco system},
see \cite{Gon03}, \cite{Apt08}, \cite{ApLy10}.
In the second case the sets of branch points $\mathcal{A}_1$ and $\mathcal{A}_2$
of the functions $f_1$ and $f_2$ {\it are equivalent}. We say that two functions $f_1,f_2$
forms a {\it Nikishin system}, see \cite{Nik80}, \cite{Gon03}, \cite{Apt08}, \cite{ApLy10}.
In the first case the interaction matrix $M$ for the theory-potential-equilibrium vector problem
is $M=\begin{pmatrix} 2&1\\1&2\end{pmatrix}$. In the second case it is
$M=\begin{pmatrix} 2&-1\\-1&2\end{pmatrix}$,
see \cite{GoRa85}, \cite{Gon03}, \cite{Apt08}, \cite{Rak11}.
The natural expectation is that
such differences between the distribution of the branch points
and the structure of the interaction matrices leads to essentially different
limiting distributions of the zeros of the corresponding Hermite--Pad\'{e} polynomials.
However, this is not always the case. Figures
\ref{Fig_ang(2_2)_4000_120_full}, \ref{Fig_nik(2_1)_4000_120_full},
\ref{Fig_ang(2_2)_4000_120_blbk} and \ref{Fig_nik(2_1)_4000_120_blrd}
show the distribution of the zeros of the polynomials $Q_{120,0}$ (blue points),
$Q_{120,1}$ (red points), $Q_{120,2}$ (black points) for two pair of functions:
one with different branch points
\begin{equation}
\label{ang1}
f_1(z)=\sqrt{(z-1)/(z+1)},\qquad f_2(z)=\sqrt{(z-2)/(z+2)}
\end{equation}
and the other with coincident branch points
\begin{equation}
\label{nik1}
f_1(z)=\(\frac{z-1}{z+1}\)^{1/3}\(\frac{z-2}{z+2}\)^{1/3},\quad
f_2(z)=\(\frac{z-1}{z+1}\)^{2/3}\(\frac{z-2}{z+2}\)^{1/3}.
\end{equation}
From the results on figures
\ref{Fig_ang(2_2)_4000_120_full}--\ref{Fig_nik(2_1)_4000_120_blrd}
it follows that the distribution of the zeros of Hermite--Pad\'{e} polynomials
for these two different systems of functions
will be equivalent (with regard to the fact that supports of the limit measures
for the polynomials $Q_{n,j}$ change over to each other). Note that
the system \eqref{nik1} is a special case of a two Markov functions system,
which was considered by E.~A.~Rakhmanov in \cite{Rak11}. The numerical results
for the system \eqref{nik1} obtained here are in a good agreement with
the results of E.~A.~Rakhmanov~\cite{Rak11}.
\end{remark}

\begin{figure}[!ht]
\centerline{
\includegraphics[width=15cm,height=15cm]{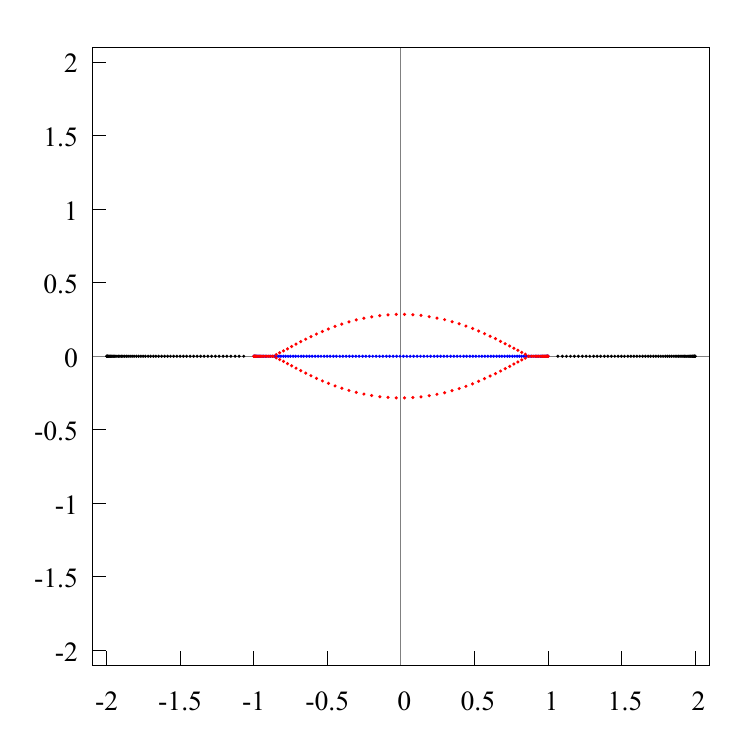}}
\vskip-6mm
\caption{The distribution of the zeros of the Hermite-Pad\'{e} polynomials $Q_{120,0}$ (blue points),
$Q_{120,1}$ (red points) and $Q_{120,2}$ (black points) for a set of three functions
$[1,f_1,f_2]$, where
$f_1(z)=\sqrt[3]{(z-1)/(z+1)}$, $f_2(z)=\sqrt[3]{(z-2)/(z+2)}$.
The sets of branch points of the functions $f_1$ and $f_2$ do not intersect each other,
and thus the pair $f_1,f_2$ create an Angelesco system.
However, the distribution of the zeros of the Hermite-Pad\'{e} polynomials
for this system is the same as the distribution for the Nikishin system
$f_1(z)=\((z-1)/(z+1)\)^{1/3}\((z-2)/(z+2)\)^{1/3}$,
$f_2(z)=\((z-1)/(z+1)\)^{2/3}\((z-2)/(z+2)\)^{1/3}$;
see fig. \ref{Fig_nik(2_1)_4000_120_full}--\ref{Fig_nik(2_1)_4000_120_blrd}.
}
\label{Fig_ang(3_3)_4000_120_full}
\end{figure}

\begin{figure}[!ht]
\centerline{
\includegraphics[width=15cm,height=15cm]{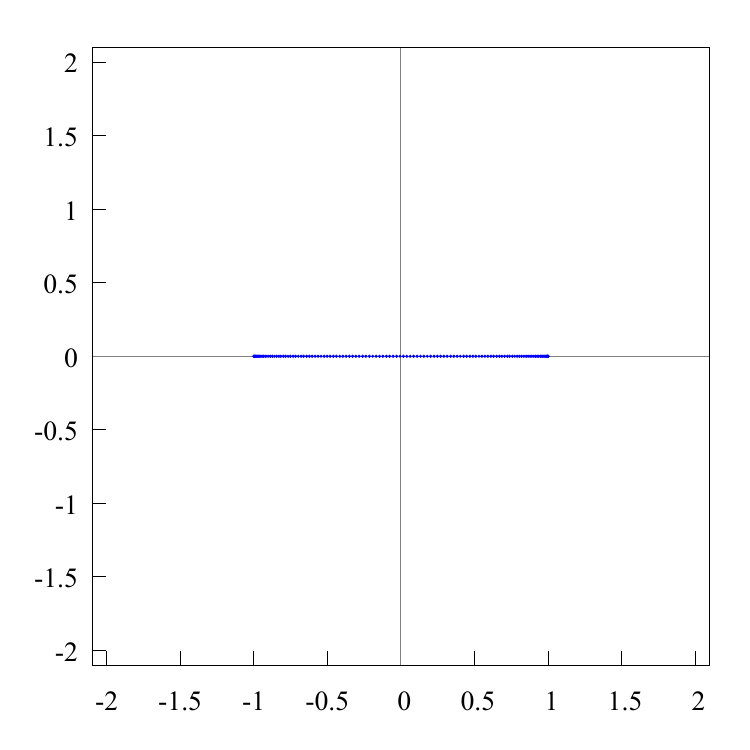}}
\vskip-6mm
\caption{The distribution of the zeros of the Hermite-Pad\'{e} polynomials $Q_{120,0}$ (blue points)
for a set of three functions $[1,f_1,f_2]$, where
$f_1(z)=\sqrt[3]{(z-1)/(z+1)}$, $f_2(z)=\sqrt[3]{(z-2)/(z+2)}$.
The pair $f_1,f_2$ create an Angelesco system.
However, the distribution of the zeros of the Hermite-Pad\'{e} polynomials
for this system is the same as the distribution for the Nikishin system
$f_1(z)=\((z-1)/(z+1)\)^{1/3}\((z-2)/(z+2)\)^{1/3}$,
$f_2(z)=\((z-1)/(z+1)\)^{2/3}\((z-2)/(z+2)\)^{1/3}$;
comp. fig. \ref{Fig_nik(2_1)_4000_120_blrd}.
}
\label{Fig_ang(3_3)_4000_120_bl}
\end{figure}

\begin{figure}[!ht]
\centerline{
\includegraphics[width=15cm,height=15cm]{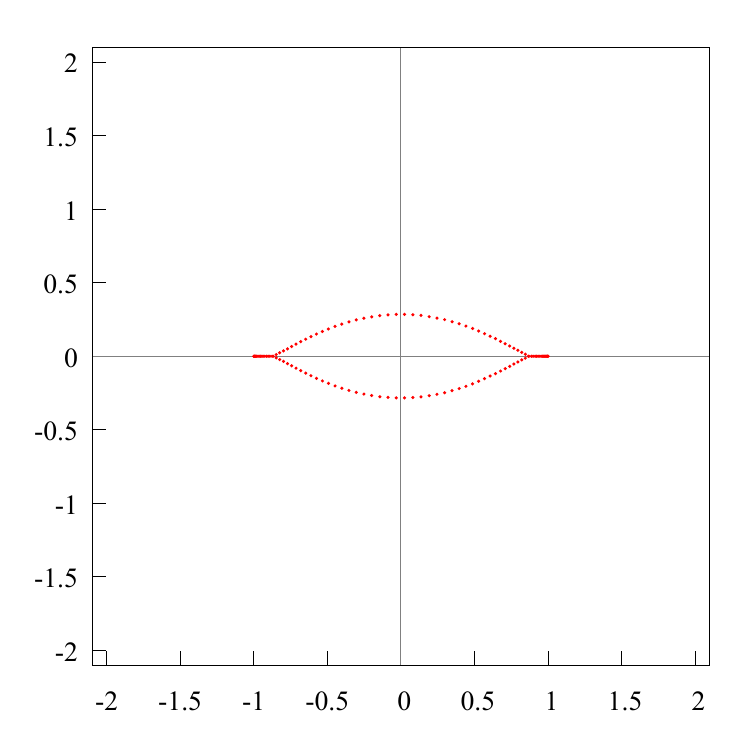}}
\vskip-6mm
\caption{The distribution of the zeros of the Hermite-Pad\'{e} polynomials $Q_{120,1}$ (red points)
for a set of three functions $[1,f_1,f_2]$, where
$f_1(z)=\sqrt[3]{(z-1)/(z+1)}$, $f_2(z)=\sqrt[3]{(z-2)/(z+2)}$.
The pair $f_1,f_2$ create an Angelesco system.
However, the distribution of the zeros of the Hermite-Pad\'{e} polynomials
for this system is the same as the distribution for the Nikishin system
$f_1(z)=\((z-1)/(z+1)\)^{1/3}\((z-2)/(z+2)\)^{1/3}$,
$f_2(z)=\((z-1)/(z+1)\)^{2/3}\((z-2)/(z+2)\)^{1/3}$;
comp. fig. \ref{Fig_nik(2_1)_4000_120_blrd}.
}
\label{Fig_ang(3_3)_4000_120_rd}
\end{figure}

\begin{figure}[!ht]
\centerline{
\includegraphics[width=15cm,height=15cm]{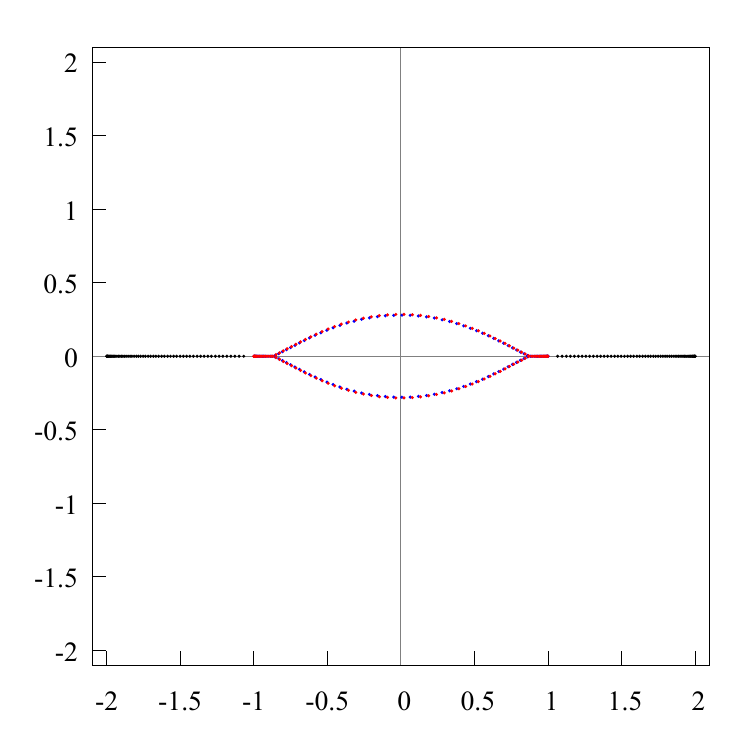}}
\vskip-6mm
\caption{The distribution of the zeros of the Hermite-Pad\'{e} polynomials $Q_{120,0}$ (blue points),
$Q_{120,1}$ (red points) and $Q_{120,2}$ (black points) for a set of three functions
$[1,f_1,f_2]$, where
$f_1(z)=\sqrt{(z-1)/(z+1)}$, $f_2(z)=\sqrt[3]{(z-2)/(z+2)}$.
The sets of branch points of the functions $f_1$ and $f_2$ do not intersect each other,
and thus the pair $f_1,f_2$ create an Angelesco system.
However, the distribution of the zeros of the Hermite-Pad\'{e} polynomials
for this system is different than the distribution of the pair
$f_1(z)=\sqrt[3]{(z-1)/(z+1)}$, $f_2(z)=\sqrt[3]{(z-2)/(z+2)}$.
comp. fig. \ref{Fig_ang(3_3)_4000_120_full}--\ref{Fig_ang(3_3)_4000_120_rd}.
This confirms that the distribution of the zeros of the Hermite-Pad\'{e} polynomials 
depends not only on the geometrical position of the branch points,
but also of the type of the branch points.
}
\label{Fig_ang(2_3)_4000_120_full}
\end{figure}

\begin{figure}[!ht]
\centerline{
\includegraphics[width=15cm,height=15cm]{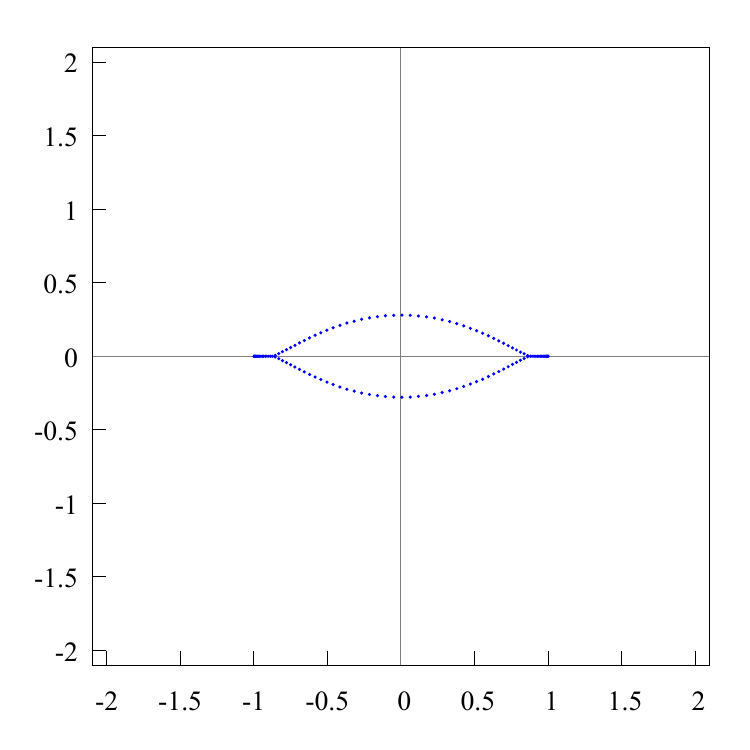}}
\vskip-6mm
\caption{The distribution of the zeros of the Hermite-Pad\'{e} polynomials $Q_{120,0}$ (blue points)
for a set of three functions $[1,f_1,f_2]$, where
$f_1(z)=\sqrt{(z-1)/(z+1)}$, $f_2(z)=\sqrt[3]{(z-2)/(z+2)}$.
The pair $f_1,f_2$ create an Angelesco system.
However, the distribution of the zeros of the Hermite-Pad\'{e} polynomials
for this system is different than the distribution of the pair
$f_1(z)=\sqrt[3]{(z-1)/(z+1)}$, $f_2(z)=\sqrt[3]{(z-2)/(z+2)}$.
comp. fig. \ref{Fig_ang(3_3)_4000_120_full}--\ref{Fig_ang(3_3)_4000_120_rd}.
This confirms that the distribution of the zeros of the Hermite-Pad\'{e} polynomials 
depends not only on the geometrical position of the branch points,
but also of the type of the branch points.
}
\label{Fig_ang(2_3)_4000_120_bl}
\end{figure}

\begin{figure}[!ht]
\centerline{
\includegraphics[width=15cm,height=15cm]{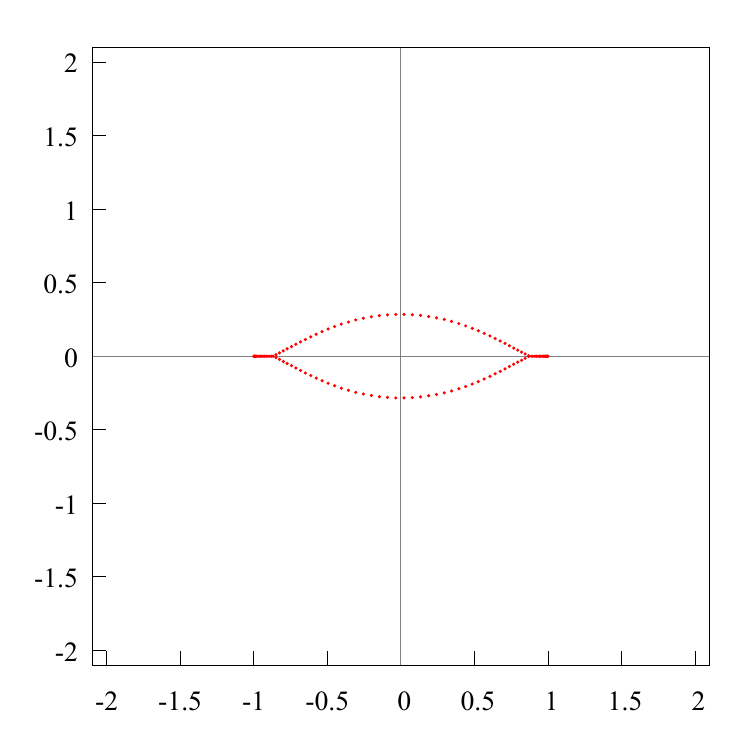}}
\vskip-6mm
\caption{The distribution of the zeros of the Hermite-Pad\'{e} polynomials $Q_{120,1}$ (red points)
for a set of three functions $[1,f_1,f_2]$, where
$f_1(z)=\sqrt{(z-1)/(z+1)}$, $f_2(z)=\sqrt[3]{(z-2)/(z+2)}$.
The pair $f_1,f_2$ create an Angelesco system.
However, the distribution of the zeros of the Hermite-Pad\'{e} polynomials
for this system is different than the distribution of the pair
$f_1(z)=\sqrt[3]{(z-1)/(z+1)}$, $f_2(z)=\sqrt[3]{(z-2)/(z+2)}$.
comp. fig. \ref{Fig_ang(3_3)_4000_120_full}--\ref{Fig_ang(3_3)_4000_120_rd}.
This confirms that the distribution of the zeros of the Hermite-Pad\'{e} polynomials 
depends not only on the geometrical position of the branch points,
but also of the type of the branch points.
}
\label{Fig_ang(2_3)_4000_120_rd}
\end{figure}

\begin{figure}[!ht]
\centerline{
\includegraphics[width=15cm,height=15cm]{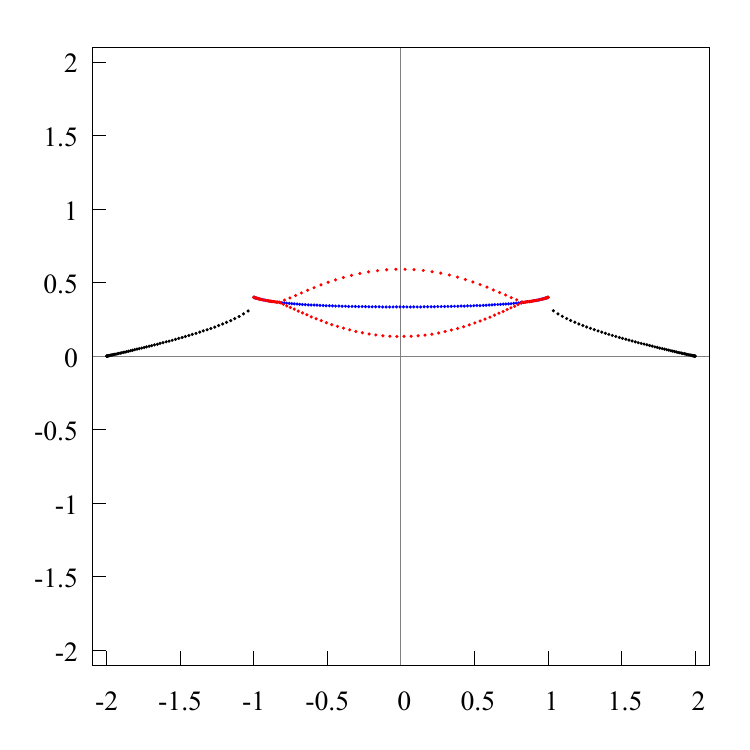}}
\vskip-6mm
\caption{The distribution of the zeros of the Hermite-Pad\'{e} polynomials $Q_{120,0}$ (blue points),
$Q_{120,1}$ (red points) and $Q_{120,2}$ (black points) for a set of three functions
$[1,f_1,f_2]$, where
$f_1(z)=\sqrt[3]{(z-(1+i\cdot 0.4))/(z-(-1+i\cdot 0.4))}$,
$f_2(z)=\sqrt[3]{(z-2)/(z+2)}$.
The set of the branch points of the functions $f_1$ and $f_2$ do not intersect each other,
thus the pair $f_1,f_2$ create an Angelesco system;
comp. \ref{Fig_ang(3_3)_4000_120_full}--\ref{Fig_ang(3_3)_4000_120_rd}.
}
\label{Fig_ang_im(_4)_3_3_4000_120_full}
\end{figure}

\begin{figure}[!ht]
\centerline{
\includegraphics[width=15cm,height=15cm]{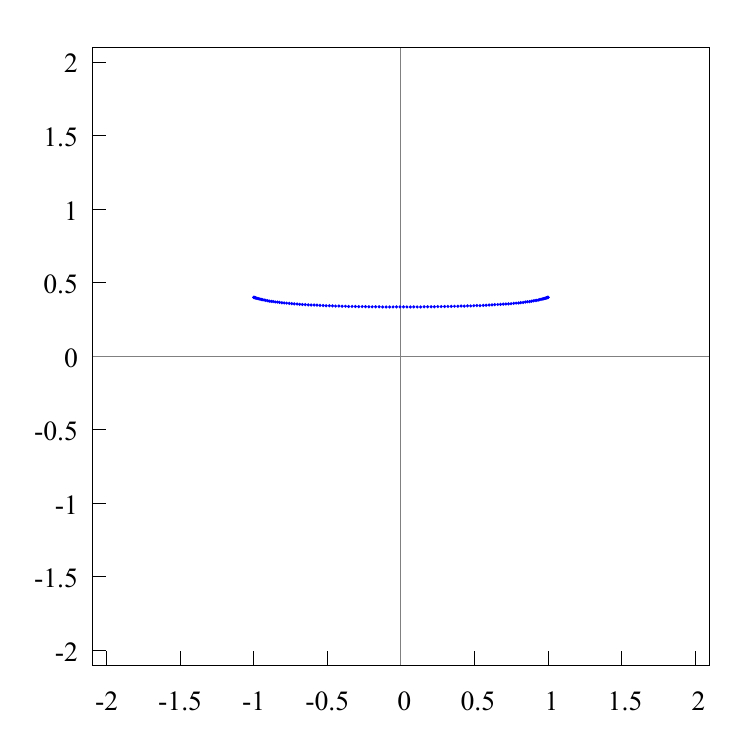}}
\vskip-6mm
\caption{The distribution of the zeros of the Hermite-Pad\'{e} polynomials $Q_{120,0}$ (blue points)
for a set of three functions $[1,f_1,f_2]$, where
$f_1(z)=\sqrt[3]{(z-(1+i\cdot 0.4))/(z-(-1+i\cdot 0.4))}$,
$f_2(z)=\sqrt[3]{(z-2)/(z+2)}$.
The set of the branch points of the functions $f_1$ and $f_2$ do not intersect each other,
thus the pair $f_1,f_2$ create an Angelesco system;
comp. \ref{Fig_ang(3_3)_4000_120_full}--\ref{Fig_ang(3_3)_4000_120_rd}.
}
\label{Fig_ang_im(_4)_3_3_4000_120_bl}
\end{figure}

\begin{figure}[!ht]
\centerline{
\includegraphics[width=15cm,height=15cm]{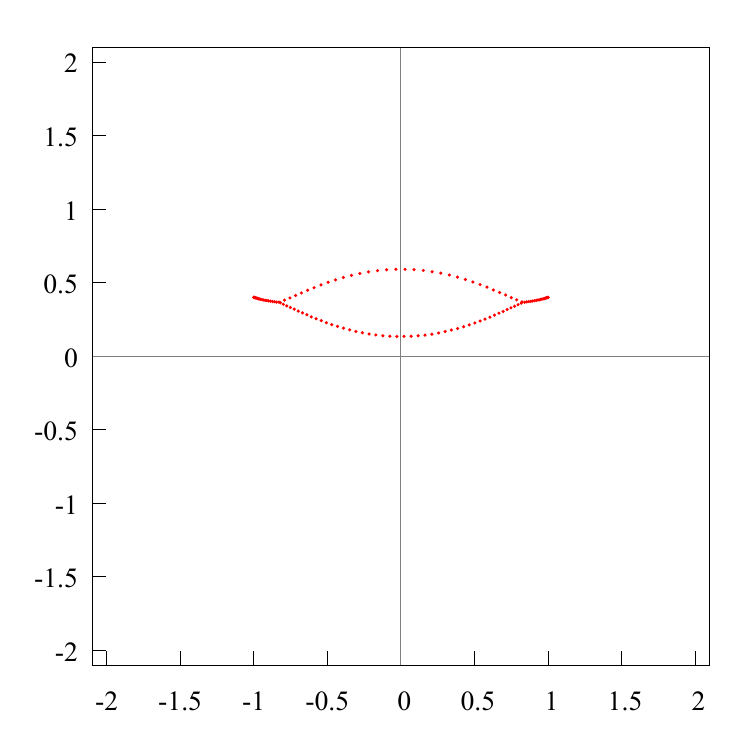}}
\vskip-6mm
\caption{The distribution of the zeros of the Hermite-Pad\'{e} polynomials $Q_{120,1}$ (red points)
for a set of three functions $[1,f_1,f_2]$, where
$f_1(z)=\sqrt[3]{(z-(1+i\cdot 0.4))/(z-(-1+i\cdot 0.4))}$,
$f_2(z)=\sqrt[3]{(z-2)/(z+2)}$.
The set of the branch points of the functions $f_1$ and $f_2$ do not intersect each other,
thus the pair $f_1,f_2$ create an Angelesco system;
comp. \ref{Fig_ang(3_3)_4000_120_full}--\ref{Fig_ang(3_3)_4000_120_rd}.
}
\label{Fig_ang_im(_4)_3_3_4000_120_rd}
\end{figure}

\begin{figure}[!ht]
\centerline{
\includegraphics[width=15cm,height=15cm]{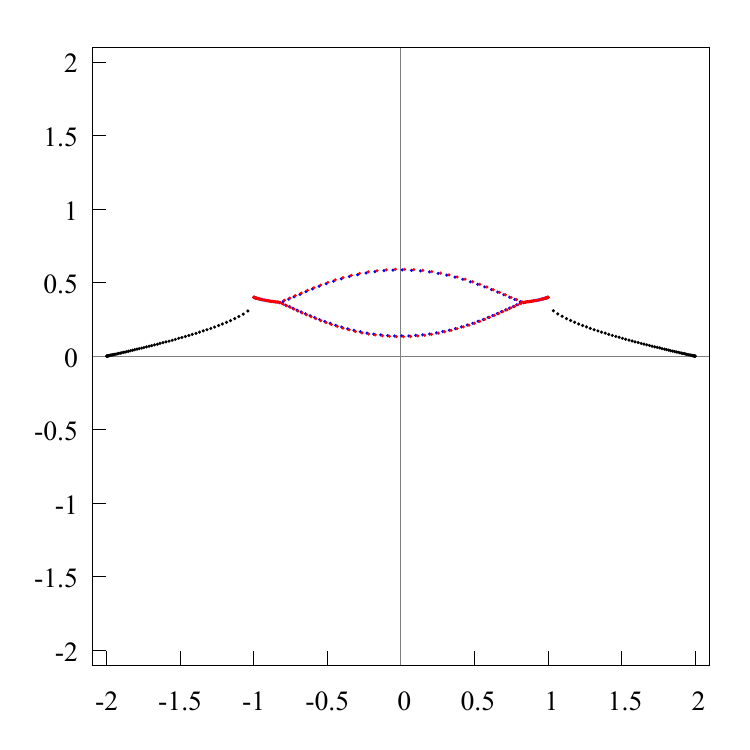}}
\vskip-6mm
\caption{The distribution of the zeros of the Hermite-Pad\'{e} polynomials $Q_{120,0}$ (blue points),
$Q_{120,1}$ (red points) and $Q_{120,2}$ (black points) for a set of three functions
$[1,f_1,f_2]$, where
$f_1(z)=\sqrt{(z-(1+i\cdot 0.4))/(z-(-1+i\cdot 0.4))}$,
$f_2(z)=\sqrt[3]{(z-2)/(z+2)}$.
The set of the branch points of the functions $f_1$ and $f_2$ do not intersect each other,
thus the pair $f_1,f_2$ create an Angelesco system;
comp. \ref{Fig_ang(2_3)_4000_120_full}--\ref{Fig_ang(2_3)_4000_120_rd}.
}
\label{Fig_ang_im(_4)_2_3_4000_120_full}
\end{figure}

\begin{figure}[!ht]
\centerline{
\includegraphics[width=15cm,height=15cm]{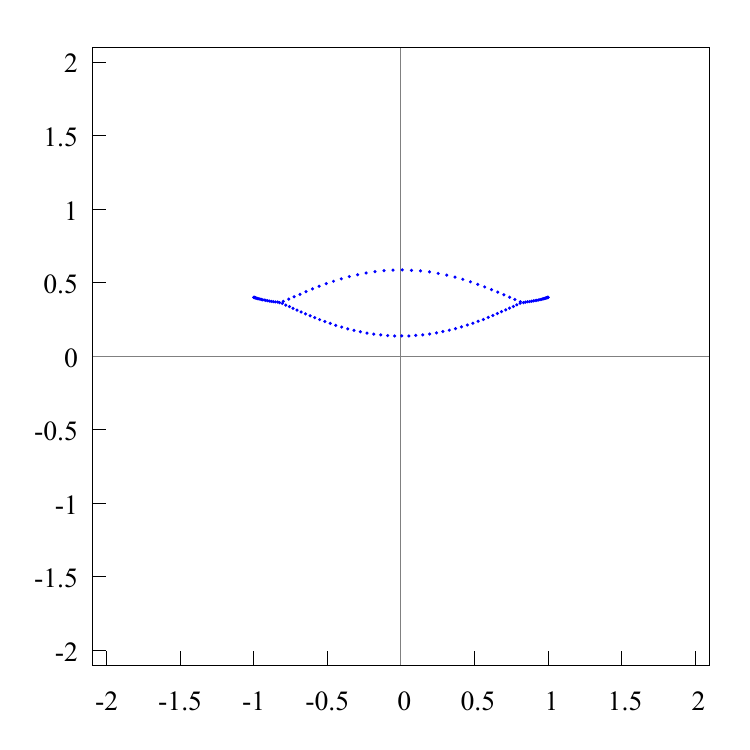}}
\vskip-6mm
\caption{The distribution of the zeros of the Hermite-Pad\'{e} polynomials $Q_{120,0}$ (blue points)
for a set of three functions $[1,f_1,f_2]$, where
$f_1(z)=\sqrt{(z-(1+i\cdot 0.4))/(z-(-1+i\cdot 0.4))}$,
$f_2(z)=\sqrt[3]{(z-2)/(z+2)}$.
The set of the branch points of the functions $f_1$ and $f_2$ do not intersect each other,
thus the pair $f_1,f_2$ create an Angelesco system;
comp. \ref{Fig_ang(2_3)_4000_120_full}--\ref{Fig_ang(2_3)_4000_120_rd}.
}
\label{Fig_ang_im(_4)_2_3_4000_120_bl}
\end{figure}

\begin{figure}[!ht]
\centerline{
\includegraphics[width=15cm,height=15cm]{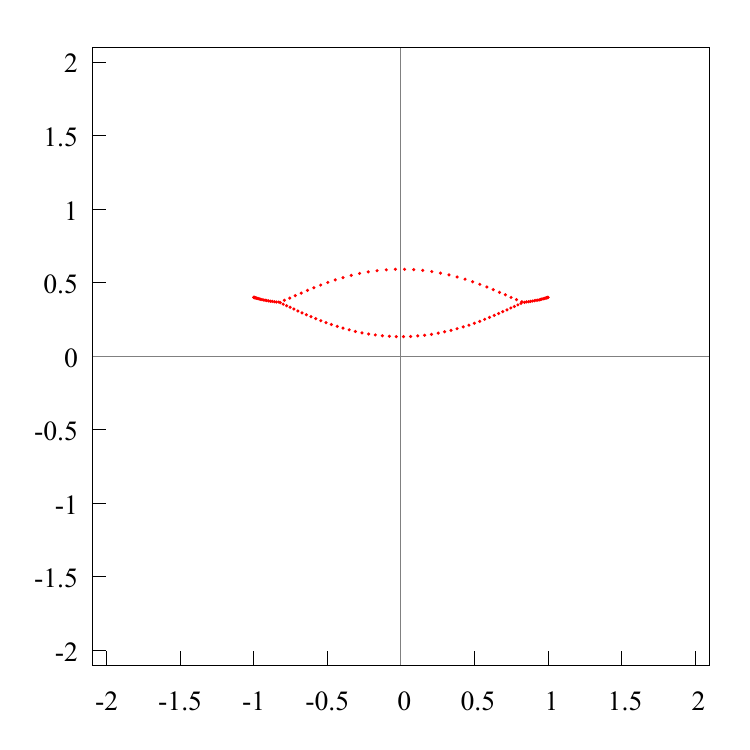}}
\vskip-6mm
\caption{The distribution of the zeros of the Hermite-Pad\'{e} polynomials $Q_{120,1}$ (red points)
for a set of three functions $[1,f_1,f_2]$, where
$f_1(z)=\sqrt{(z-(1+i\cdot 0.4))/(z-(-1+i\cdot 0.4))}$,
$f_2(z)=\sqrt[3]{(z-2)/(z+2)}$.
The set of the branch points of the functions $f_1$ and $f_2$ do not intersect each other,
thus the pair $f_1,f_2$ create an Angelesco system;
comp. \ref{Fig_ang(2_3)_4000_120_full}--\ref{Fig_ang(2_3)_4000_120_rd}.
}
\label{Fig_ang_im(_4)_2_3_4000_120_rd}
\end{figure}

\begin{remark}\label{rem2}
Instead of the pair of functions \eqref{ang1} with second order branch points
we can consider a pair of functions with arbitrary order of the branch points.
For example, the pair of functions
\begin{equation}
\label{ang2}
f_1(z)=\sqrt[3]{(z-1)/(z+1)},\qquad f_2(z)=\sqrt[3]{(z-2)/(z+2)}
\end{equation}
have the same distribution of the zeros of the corresponding Hermite--Pad\'{e} polynomials
as \eqref{ang1};
see fig. \ref{Fig_ang(3_3)_4000_120_full}--\ref{Fig_ang(3_3)_4000_120_rd}
and comp. fig. \ref{Fig_ang_im(_4)_3_3_4000_120_full}--\ref{Fig_ang_im(_4)_3_3_4000_120_rd}.
However, the pair of functions
\begin{equation}
\label{ang3}
f_1(z)=\sqrt{(z-1)/(z+1)},\qquad f_2(z)=\sqrt[3]{(z-2)/(z+2)}
\end{equation}
have another distribution of the zeros of the corresponding Hermite--Pad\'{e} polynomials,
see fig. \ref{Fig_ang(2_3)_4000_120_full}--\ref{Fig_ang(2_3)_4000_120_rd}
and comp. fig. \ref{Fig_ang_im(_4)_2_3_4000_120_full}--\ref{Fig_ang_im(_4)_2_3_4000_120_rd}.
Thus, the distribution depends non only of the type but also on
the degree of the branch points
(in \eqref{ang3} for the first function $f_1$ instead of degree $1/3$,
as in \eqref{ang2}, we took degree $1/2$).

The distribution of the zeros remains stable while moving the branch points
of the function $f_1$ along the imaginary axis, for example,
when we change the points $z=\pm1$ into the points $z=\pm1+i\cdot 0.4$:
\begin{align}
\label{ang4}
f_1(z)=\sqrt[3]{(z-(1+i\cdot 0.4))/(z-(-1+i\cdot 0.4))},\\
f_2(z)=\sqrt[3]{(z-2)/(z+2)},\nonumber
\end{align}
see fig. \ref{Fig_ang_im(_4)_3_3_4000_120_full}--\ref{Fig_ang_im(_4)_3_3_4000_120_rd}
and comp. \eqref{ang2}, and:
\begin{align}
\label{ang5}
f_1(z)=\sqrt{(z-(1+i\cdot 0.4))/(z-(-1+i\cdot 0.4))},\\
f_2(z)=\sqrt[3]{(z-2)/(z+2)},\nonumber
\end{align}
see fig. \ref{Fig_ang_im(_4)_2_3_4000_120_full}--\ref{Fig_ang_im(_4)_2_3_4000_120_rd}
and comp. \eqref{ang3}.
This confirms that the distribution of the zeros of the Hermite--Pad\'{e} polynomials
depends not only on the geometrical position of the branch points,
but also on the type of the branch points.
\end{remark}

\begin{figure}[!ht]
\centerline{
\includegraphics[width=15cm,height=15cm]{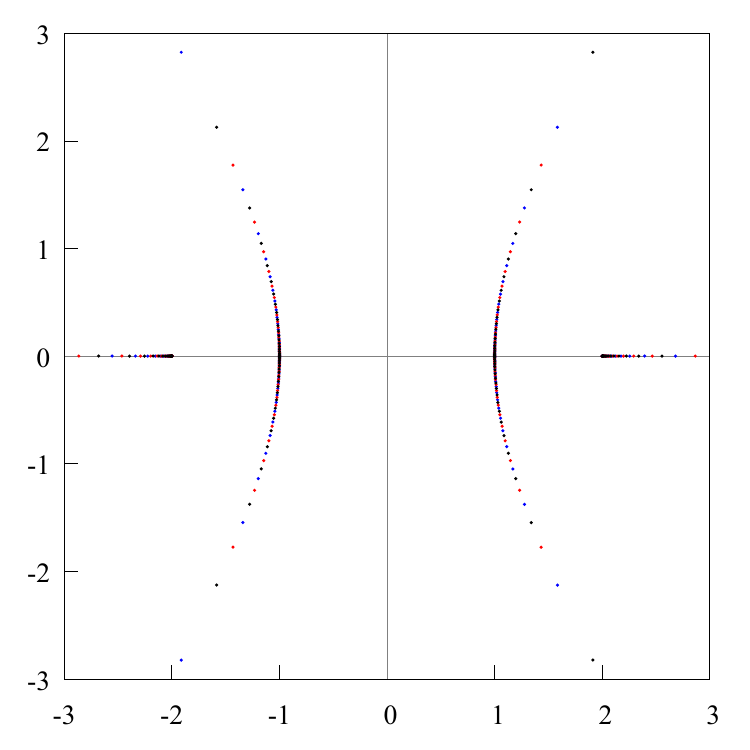}}
\vskip-6mm
\caption{The distribution of the zeros of the Hermite-Pad\'{e} polynomials $Q_{120,0}$ (blue points),
$Q_{120,1}$ (red points) and $Q_{120,2}$ (black points) for a set of three functions
$[1,f_1,f_2]$, where
$f_1(z)=\((z-1)/(z+1)\)^{1/3}\((z-2)/(z+2)\)^{1/3}$,
$f_2(z)=\((z-1)/(z+1)\)^{2/3}\((z-2)/(z+2)\)^{-1/3}$.
Here the pair $f_1,f_2$ create a Nikishin system,
because the branch points of the functions are equivalent.
}
\label{Fig_nik(2_-1)_4000_120_full}
\end{figure}

\begin{figure}[!ht]
\centerline{
\includegraphics[width=15cm,height=15cm]{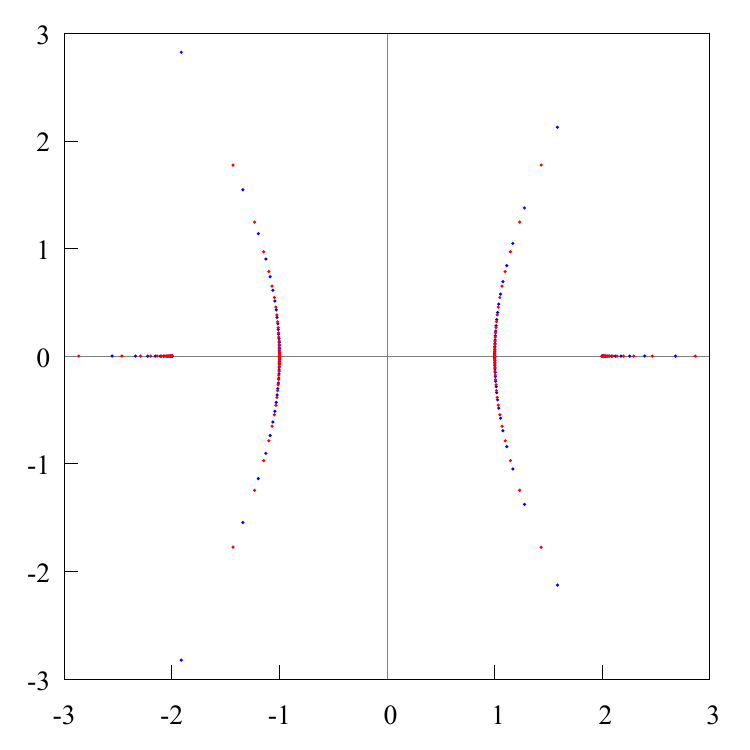}}
\vskip-6mm
\caption{The distribution of the zeros of the Hermite-Pad\'{e} polynomials $Q_{120,0}$ (blue points)
and $Q_{120,1}$ (red points) for a set of three functions
$[1,f_1,f_2]$, where
$f_1(z)=\((z-1)/(z+1)\)^{1/3}\((z-2)/(z+2)\)^{1/3}$,
$f_2(z)=\((z-1)/(z+1)\)^{2/3}\((z-2)/(z+2)\)^{-1/3}$.
The pair $f_1,f_2$ create a Nikishin system.
}
\label{Fig_nik(2_-1)_4000_120_blrd}
\end{figure}

\begin{remark}\label{rem3}
If instead of \eqref{nik1} we consider the pair of functions
\begin{equation}
\label{nik2}
f_1(z)=\(\frac{z-1}{z+1}\)^{1/3}\(\frac{z-2}{z+2}\)^{1/3},\quad
f_2(z)=\(\frac{z-1}{z+1}\)^{2/3}\(\frac{z-2}{z+2}\)^{-1/3},
\end{equation}
then the distribution of the zeros of the corresponding Hermite--Pad\'{e} polynomials
will be different than before,
see fig. \ref{Fig_nik(2_-1)_4000_120_full}--\ref{Fig_nik(2_-1)_4000_120_blrd}.
\end{remark}

\subsection{}\label{s1s2}
A numerous conjectures about the asymptotic behavior of the zeros of
the Hermite--Pad\'{e} polynomials of both first and second kind was, in
part, formulated in the fundamental work of J. Nuttall \cite{Nut84}
(see also \cite{Nut81}, \cite{Nut82}). These conjectures served as a
background for some of the numerous investigations. Among them the
prominent results of H.~Stahl~\cite{Sta85a}--\cite{Sta86b} and
A.~A.~Gonchar--E.~A.~Rakhmanov~\cite{GoRa81},~\cite{GoRa87}, first of
all should be noted, and also the results of A.~I.~Aptekarev with
co-authors \cite{ApKa86}, \cite{ApKuVa08}, \cite{Apt08}, \cite{ApLy10},
\cite{ApLyTu11}, \cite{ApYa11}. J. Nuttall's exact results
\cite{Nut81}, \cite{Nut82}, \cite{Nut84} about the asymptotics of the
Hermite--Pad\'{e} polynomials for a collection of $m$ functions
$[f_0\equiv1,f_1,\dots,f_{m-1}]$, $m\ge3$, are based mainly on the a
priori assumption of the existence of an associated $m$-sheeted Riemann
surface $\RS_m$ with a canonical decomposition into $m$ sheets
$\RS^{(j)}$, $j=1,\dots,m$. The sheets defined by an Abel integral of
third kind with purely imaginary periods and logarithmic singularities
of the form $(m-1)\log{z}$, $z=z^{(1)}\sim\infty^{(1)}$ and $-\log{z}$,
$z=z^{(j)}\sim\infty^{(j)}$, $j=2,\dots,m$, and having also property
that the first ``physical'' sheet $\RS^{(1)}$ is always connected (see
\cite{Nut81}, \cite{Nut84}, \cite{NuTr87}).\footnote{ In some cases,
considered by Nuttall, it appears that the multi-valued analytic
functions $f_1,\dots,f_{m-1}$ continue from the first sheet onto
$\RS_m\setminus\RS^{(m)}$ as single-valued meromorphic functions and
$[1,f_1,\dots,f_{m-1}]$ are in a sense independent.} In such case it is
possible to describe the asymptotical behavior of the Hermite--Pad\'{e}
polynomials in terms connected to this RS (see \cite{Nut84},
\cite{Nut86}, \cite{ApKuVa08}, \cite{Apt08}). However, the question how
to ``construct'' such a RS in the general case when $m\geq3$ remains
open.\footnote{ For $m=2$ the existence of a corresponding
hyperelliptic Riemann surface with two sheets follows directly from the
theorems by Stahl; see \cite{Nut84}, \cite{Sta97b}, \cite{ApYa11},
\cite{KovSu14}.} As for the present, results are achieved only for the
case $m=3$. There are two major methods to obtain the results of such
type. The first one mainly based on the {\it cubic equation} (see
\cite{ApKa86}, \cite{ApKuVa08}, \cite{ApTu12b}, also \cite{Nut86},
\cite{MaRaSu12}, \cite{MaRaSu13}). Using this method one can find an
explicit representation of the so-called Nuttall's
$\phi$-function\footnote{In the case of two sheeted hyperelliptic RS
such function is known as Deift's $g$-function; see~\cite{RaSu13}.}
(see \cite{ApKuVa08}, \cite{ApTu12b}, \cite{ApTu14}), i.e., an Abel
integral of third kind, which has the level curves that define the
needed canonical decomposition of RS $\RS_3$ into three sheets . The
other method consists of first solving the theoretical-potential
extremal problem, connected with the existence of the so-called {\it
Nuttall condenser} on the Riemann sphere $\myo\CC$ (see \cite{RaSu12},
\cite{RaSu13}, \cite{KovSu14}, \cite{Sue13b}). If such a condenser can
be found, then the RS $\RS_3$ with the needed decomposition into three
sheets is ``constructed'' on the base of this condenser by a specific
scheme (see \cite{RaSu13}, \cite{KovSu14}). Note that for now the
second method can be applied only when the pair $f_1,f_2$ forms a
Nikishin system \cite{RaSu13}, \cite{Sue13b}. Namely, for some general
enough functions $f_1,f_2$, which forms a {\it complex} Nikishin
system,\footnote{Under (general) Nikishin system, which contains two
functions $f_1,f_2$, we understand such a system for which the
corresponding vector theoretical-potential equilibrium problem is
defined by a Nikishin matrix; see \cite{Apt08}, \cite{ApLy10},
\cite{Lap12}.} in \cite{RaSu13} (see also \cite{RaSu12},
\cite{KovSu14}) an {\it existence} of the Nuttall condenser, consisting
of two non-intersecting plates and having certain symmetrical
properties, is proved. Such a condenser is an analogue of Stahl
compact, but for the case of Hermite--Pad\'{e} polynomials when $m=3$.
In \cite{KovSu14} a scheme for constructing a RS with three sheets and
with canonical decomposition made on the base of an already existing
Nuttall's condenser is proposed.

As mentioned before, the numerical results are obtained here for two opposite
cases:
for a pair of functions $f_1,f_2$ creating an Angelesco system (see \eqref{3})
and for a pair of functions $f_1=f,f_2=f^2$ creating a (generalized) Nikishin system (see \eqref{4}).
These new numerical results give rise to some new conjectures about
the asymptotical properties of the Hermite--Pad\'{e} polynomials of first kind.
It is well known that the Hermite--Pad\'{e} polynomials of first and second kind
are closely related, in particular, they are bi-orthogonal;
see \cite[\S 2, formula (2.1.9)]{Nut84},
also the recent works \cite{Ass01}, \cite{AsFi13}, \cite{Ass11}.
It seems that the conjectures presented here can be applied also to
the Hermite--Pad\'{e} polynomials of second kind.

\begin{figure}[!ht]
\centerline{
\includegraphics[width=15cm,height=15cm]{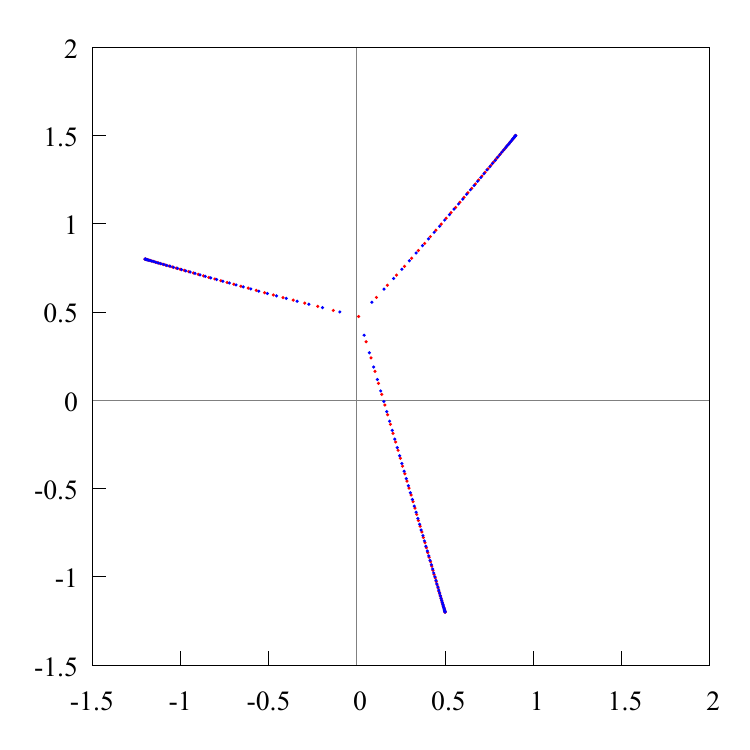}}
\vskip-6mm
\caption{Zeros and poles of the diagonal Pad\'{e} approximant $[131/131]_f$ of the function
$f(z)=1/\{(z-(-1.2+i\cdot0.8))(z-(0.9+i\cdot1.5))(z-(0.5-i\cdot1.2))\}^{1/3}$.
The zeros and poles of the diagonal Pad\'{e} approximant $[n/n]_f$ are distributed,
when $n\to\infty$, according to Stahl's theorem \cite{Sta97b}.
When $n=131$ the zeros and poles of the diagonal Pad\'{e} approximant $[131/131]_f$
of the function $f$ are distributed accordingly to the electrostatical model,
by E. A. Rakhmanov \cite{Rak12}, and they model
the Chebotarev-Stahl compact $S_{131}$, which depends on $n$.
Froissart doublets are not present when $n=131$.
}
\label{Fig_pade10_2500_131}
\end{figure}

\begin{figure}[!ht]
\centerline{
\includegraphics[width=15cm,height=15cm]{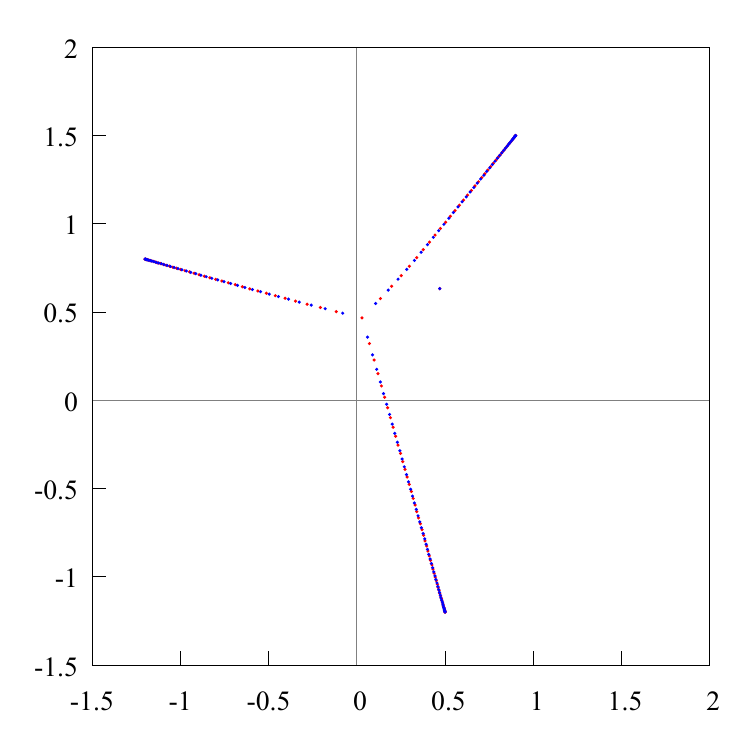}}
\vskip-6mm
\caption{Zeros and poles of the diagonal Pad\'{e} approximant $[130/130]_f$ of the function
$f(z)=1/\{(z-(-1.2+i\cdot0.8i))(z-(0.9+i\cdot1.5))
(z-(0.5-i\cdot1.2))\}^{1/3}$, distributed accordingly
to the electrostatical model by E. A. Rakhmanov \cite{Rak12}.
There is a Froissart doublet when $n=130$ (see also fig.
\ref{Fig_pade10_2500_130_red} and fig. \ref{Fig_pade10_2500_130_blu}).
Since the genus of the Riemann surface is $1$,
there might be at most one Froissart doublet.
In full compliance with the Rakhmanov model \cite{Rak12},
the Froissart doublet ``attracts'' the Stahl $S$-compact $S_{130}$,
comp. fig. \ref{Fig_pade10_2500_131}.
}
\label{Fig_pade10_2500_130}
\end{figure}

\begin{figure}[!ht]
\centerline{
\includegraphics[width=15cm,height=15cm]{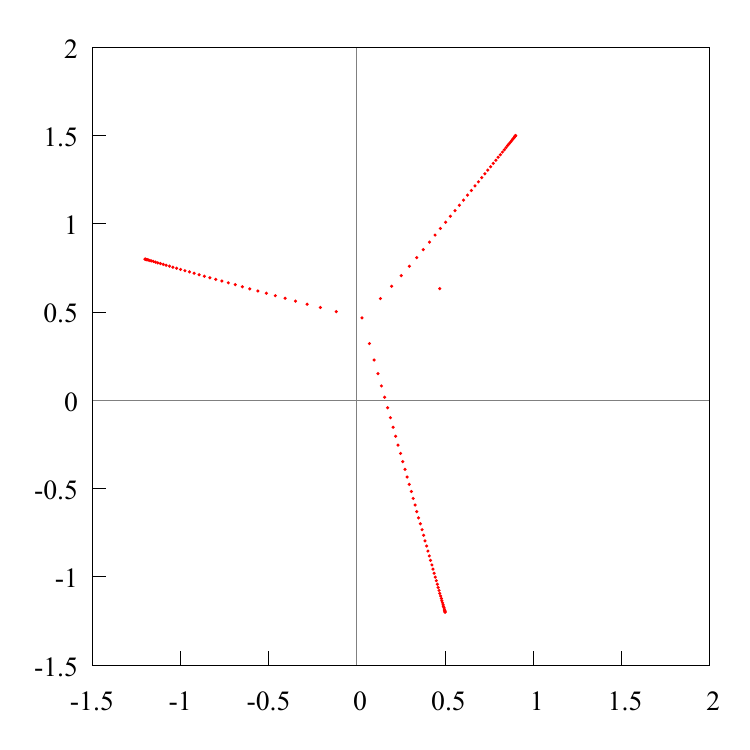}}
\vskip-6mm
\caption{The poles of the Pad\'{e} approximant $[130/130]_f$
approximate a Chebotarev point $v_{130}$
for the $S$-compact $S_{130}$. When $n\to\infty$ we have that $v_n\to v$ is a classical
Chebotarev point. There is one spurious pole of the Pad\'{e} approximant $[130/130]_f$,
it is accompanied by a spurious zero of the Pad\'{e} approximant $[130/130]_f$
(see fig. \ref{Fig_pade10_2500_130_blu}).
}
\label{Fig_pade10_2500_130_red}
\end{figure}

\begin{figure}[!ht]
\centerline{
\includegraphics[width=15cm,height=15cm]{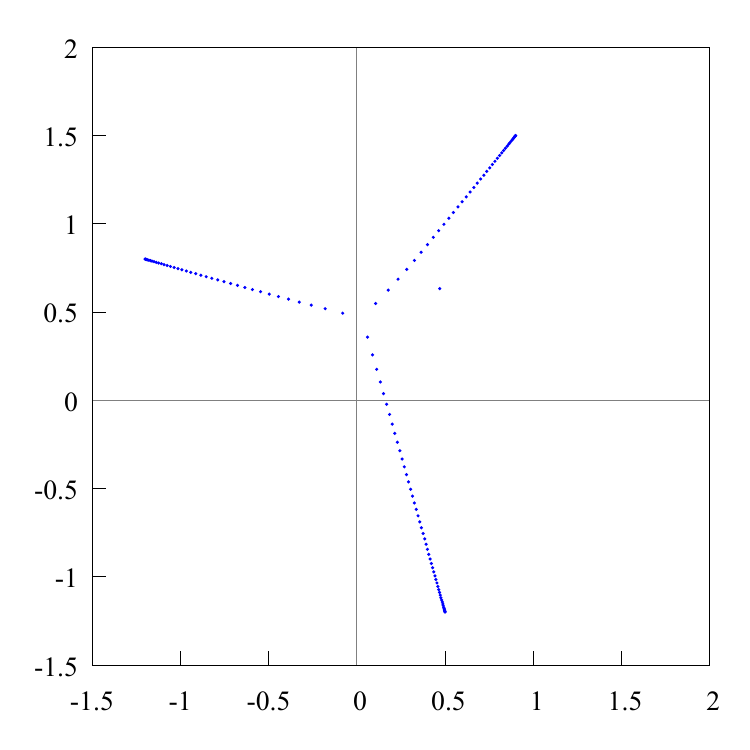}}
\vskip-6mm
\caption{The Chebotarev point should not be approximated by
zeros of the Pad\'{e} approximant $[130/130]_f$.
There is one spurious zero of the Pad\'{e} approximant $[130/130]_f$,
it is accompanied by a spurious pole of the Pad\'{e} approximant $[130/130]_f$
(see fig. \ref{Fig_pade10_2500_130_red}).
}
\label{Fig_pade10_2500_130_blu}
\end{figure}

\begin{figure}[!ht]
\centerline{
\includegraphics[width=15cm,height=15cm]{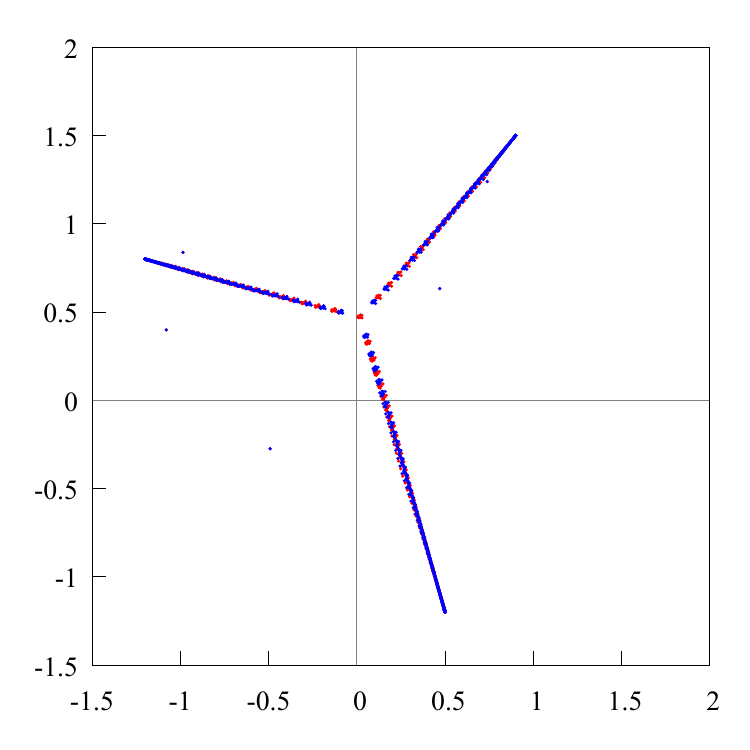}}
\vskip-6mm
\caption{Zeros and poles of the diagonal Pad\'{e} approximants $[n/n]_f$, $n=121,\dots,130$,
of the function $f(z)=1/\{(z-(-1.2+i\cdot0.8))(z-(0.9+i\cdot1.5))
(z-(0.5-i\cdot1.2))\}^{1/3}$ for each $n=121,\dots,130$
are distributed in the complex plane accordingly to the electrostatical model
by E. A. Rakhmanov \cite{Rak12}.
Since the genus of the Riemann surface is $1$,
there might be at most one Froissart doublet for each $n$.
Here are observed 5 Froissart doublets.
In full compliance with the Rakhmanov model \cite{Rak12}
the $n$-th Froissart doublet ``attracts''
the Stahl $S$-compact $S_{n}$, $n\in\{121,\dots,130\}$,
comp. fig. \ref{Fig_pade10_2500_131}.
}
\label{Fig_pade10_2500_121_130}
\end{figure}

\begin{figure}[!ht]
\centerline{
\includegraphics[width=15cm,height=15cm]{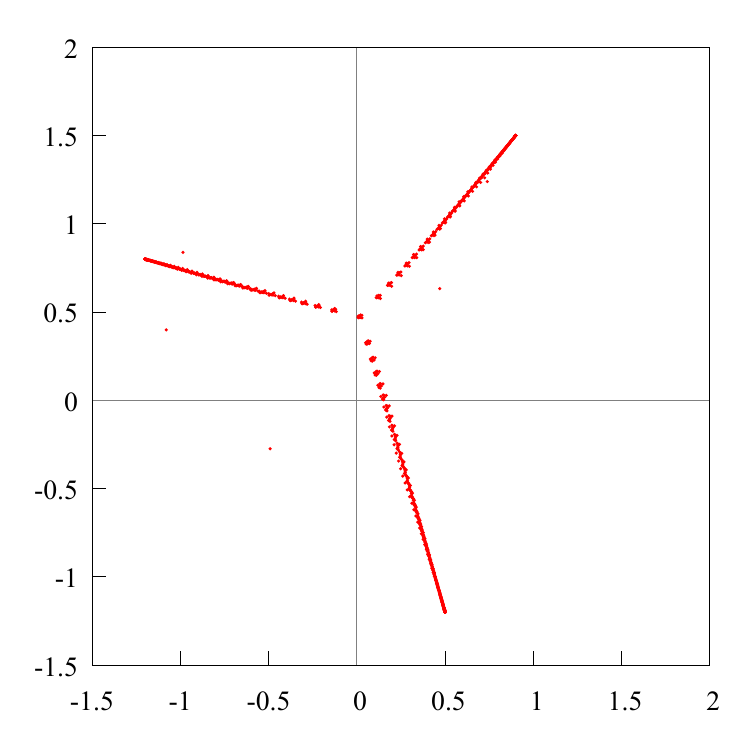}}
\vskip-6mm
\caption{For the function $f(z)=1/\{(z-(-1.2+i\cdot0.8))(z-(0.9+i\cdot1.5))
(z-(0.5-i\cdot1.2))\}^{1/3}$
poles of the Pad\'{e} approximants $[n/n]_f$, $n=121,\dots,130$,
approximate a Chebotarev point.
However, the presence of Froissart doublets
changes slightly the position of the existing points $v_n$
depending on $n\in\{121,\dots,130\}$.
This is in full accordance with the Rakhmanov model \cite{Rak12}.
}
\label{Fig_pade10_2500_121_130_red}
\end{figure}

\begin{figure}[!ht]
\centerline{
\includegraphics[width=15cm,height=15cm]{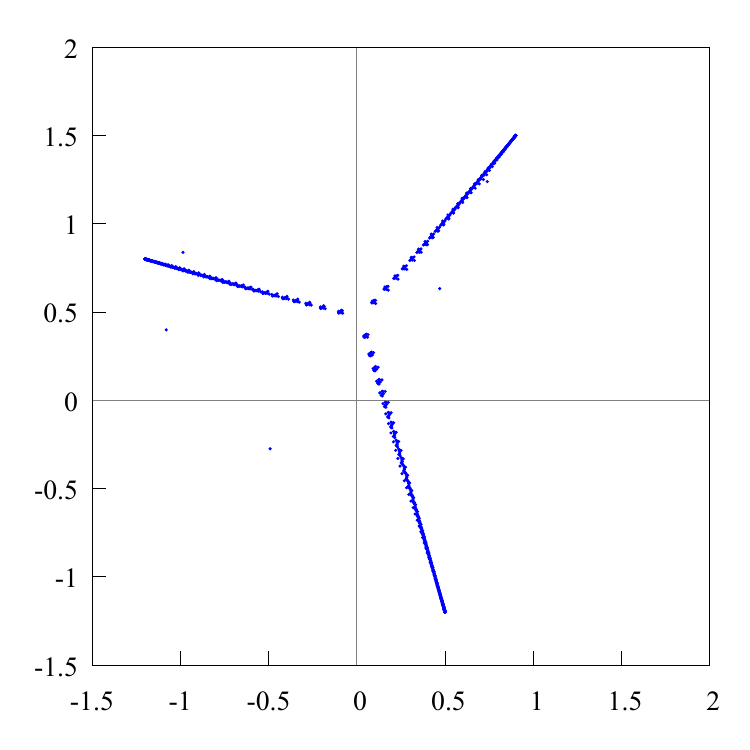}}
\vskip-6mm
\caption{Zeros of the Pad\'{e} approximants $[n/n]_f$, $n=121,\dots,130$, for the function
$f(z)=1/\{(z-(-1.2+i\cdot0.8))(z-(0.9+i\cdot1.5))(z-(0.5-i\cdot1.2))\}^{1/3}$.
The Chebotarev point should not be approximated by
zeros of the Pad\'{e} approximant $[n/n]_f$.
}
\label{Fig_pade10_2500_121_130_blu}
\end{figure}

\begin{figure}[!ht]
\centerline{
\includegraphics[width=15cm,height=15cm]{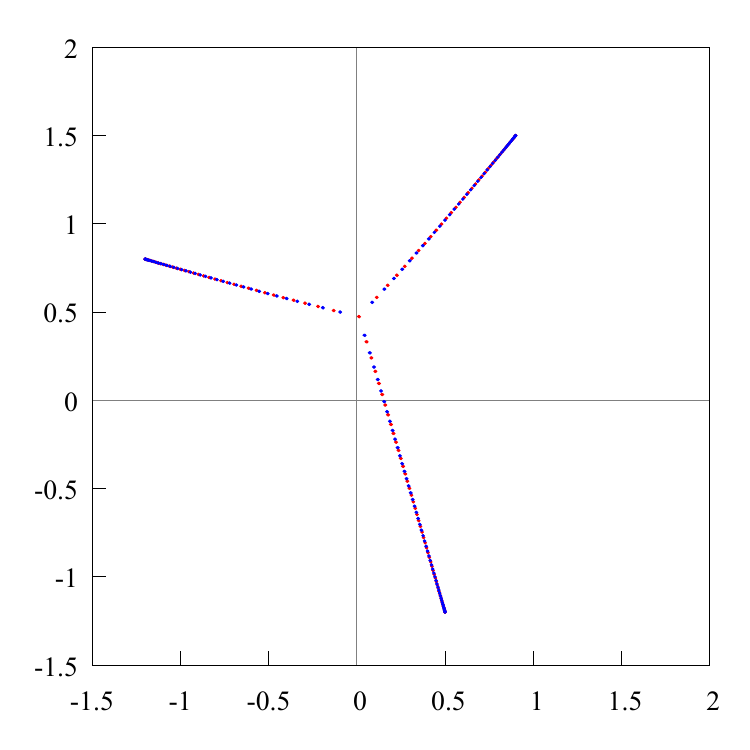}}
\vskip-6mm
\caption{Zeros and poles of the diagonal Pad\'{e} approximants $[n/n]_f$
when $n=131$ and $n=132$ of the function
$f(z)=1/\{(z-(-1.2+i\cdot0.8))(z-(0.9+i\cdot1.5))(z-(0.5-i\cdot1.2))\}^{1/3}$.
The Froissart doublets are not present when $n=131$ and $n=132$.
Therefore the $S$-compacts $S_{131}$ and $S_{132}$ are practically the same.
}
\label{Fig_pade10_2500_131_132}
\end{figure}

\begin{figure}[!ht]
\centerline{
\includegraphics[width=15cm,height=15cm]{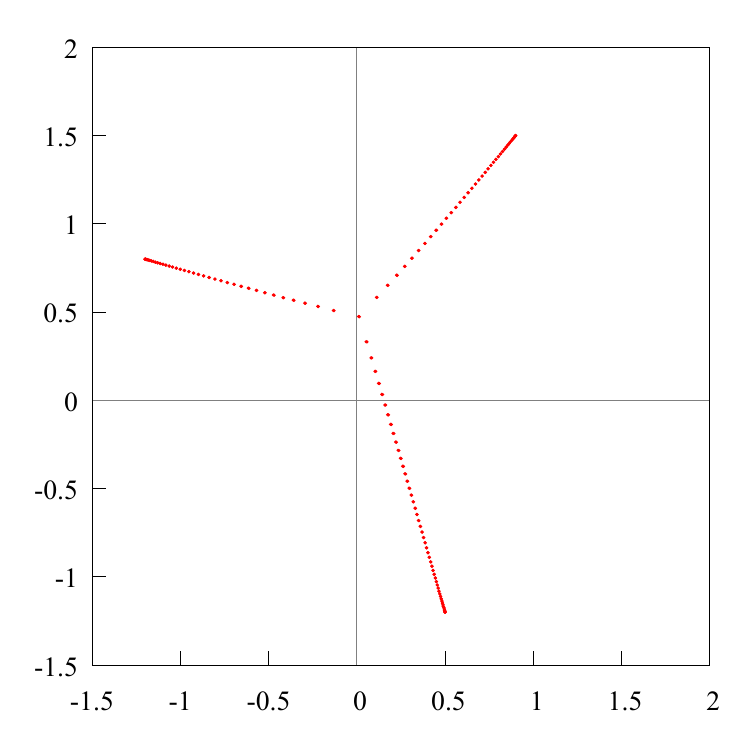}}
\vskip-6mm
\caption{Poles of the diagonal Pad\'{e} approximants $[n/n]_f$
when $n=131$ and $n=132$ of the function
$f(z)=1/\{(z-(-1.2+i\cdot0.8))(z-(0.9+i\cdot1.5))(z-(0.5-i\cdot1.2))\}^{1/3}$.
The Froissart doublets are not present when $n=131$ and $n=132$.
Therefore, the Chebotarev points $v_{131}$ and $v_{132}$, approximated by
the poles of the diagonal Pad\'{e} approximants $[n/n]_f$ when
$n=131$ and $n=132$, are very close to each other and are practically the same.
}
\label{Fig_pade10_2500_131_132_red}
\end{figure}

\begin{figure}[!ht]
\centerline{
\includegraphics[width=15cm,height=15cm]{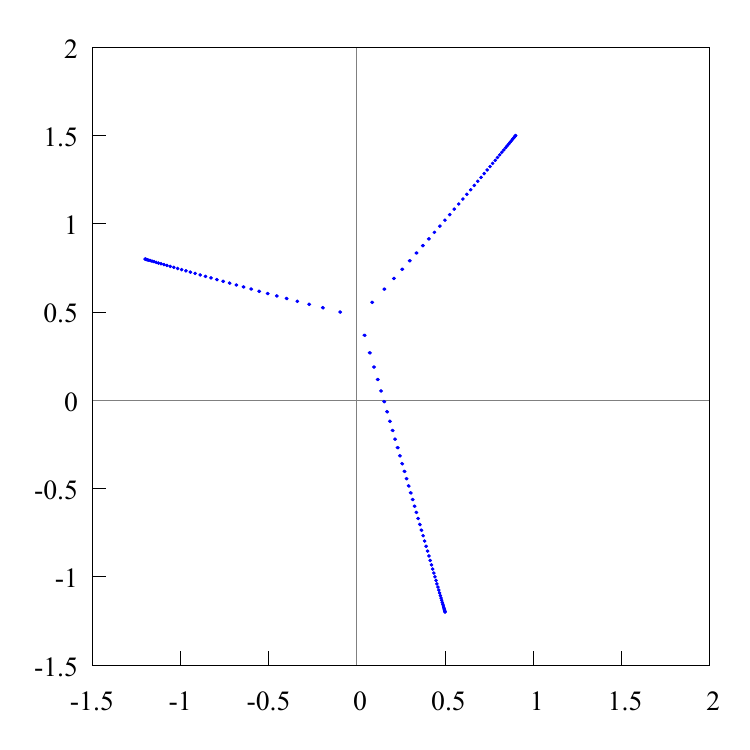}}
\vskip-6mm
\caption{Zeros of the Pad\'{e} approximants $[n/n]_f$ when $n=131$ and $n=132$ of the function
$f(z)=1/\{(z-(-1.2+i\cdot0.8))(z-(0.9+i\cdot1.5))(z-(0.5-i\cdot1.2))\}^{1/3}$
do not approximate the Chebotarev points $v_{131}$ and $v_{132}$.
The Froissart doublets are not present when $n=131$ and $n=132$.
Therefore the zeros the diagonal Pad\'{e} approximants $[n/n]_f$ when $n=131$ and $n=132$
are practically the same.
}
\label{Fig_pade10_2500_131_132_blu}
\end{figure}

\begin{figure}[!ht]
\centerline{
\includegraphics[width=15cm,height=15cm]{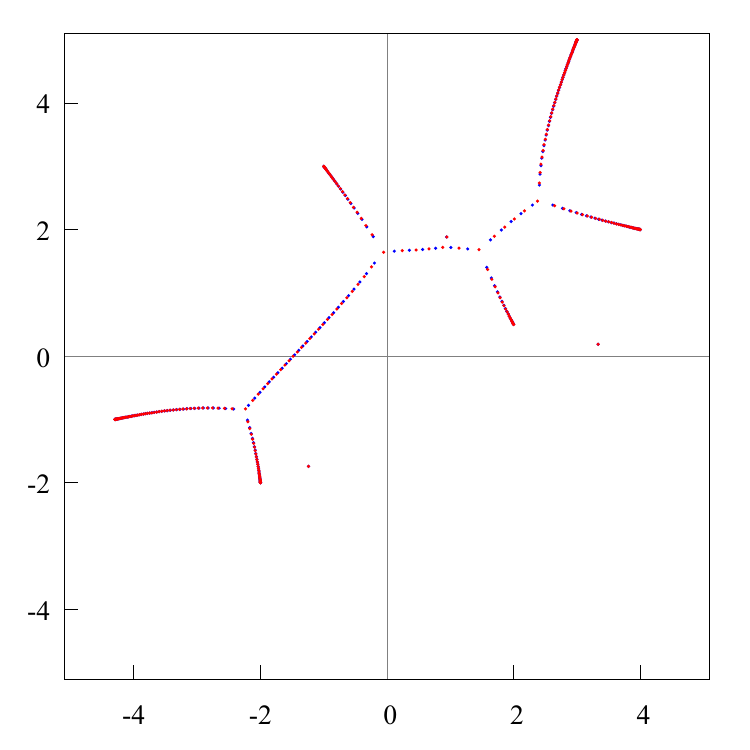}}
\vskip-6mm
\caption{Zeros and poles of the diagonal Pad\'{e} approximants $[103/103]_f$ of the function
$f(z)=1/\bigl\{(z+(4.3+i\cdot1.))
(z-(2.+i\cdot.5))\*(z+(2.+i\cdot2.))\*(z+(1.-i\cdot3.))
(z-(4.+i\cdot2.))(z-(3.+i\cdot5.))\bigr\}^{1/6}$.
In the limit when $n\to\infty$
the zeros and poles of the diagonal Pad\'{e} approximants $[n/n]_f$
are distributed accordingly with Stahl's theorem \cite{Sta97b}.
Under fixed $n=103$ these zeros and poles are distributed in a plane,
accordingly to the electrostatical model by Rakhmanov \cite{Rak12}.
Since the genus of the Riemann surface is $1$,
for each $n$ there might be no more than 4 Froissart doublets.
Here are observed 3 Froissart doublets.
When $n=103$, in full compliance with the Rakhmanov model,
the Froissart doublets ``attract'' the Stahl $S$-compact $S_{103}$.
}
\label{Fig_pade103_5000_266_full}
\end{figure}

\begin{figure}[!ht]
\centerline{
\includegraphics[width=15cm,height=15cm]{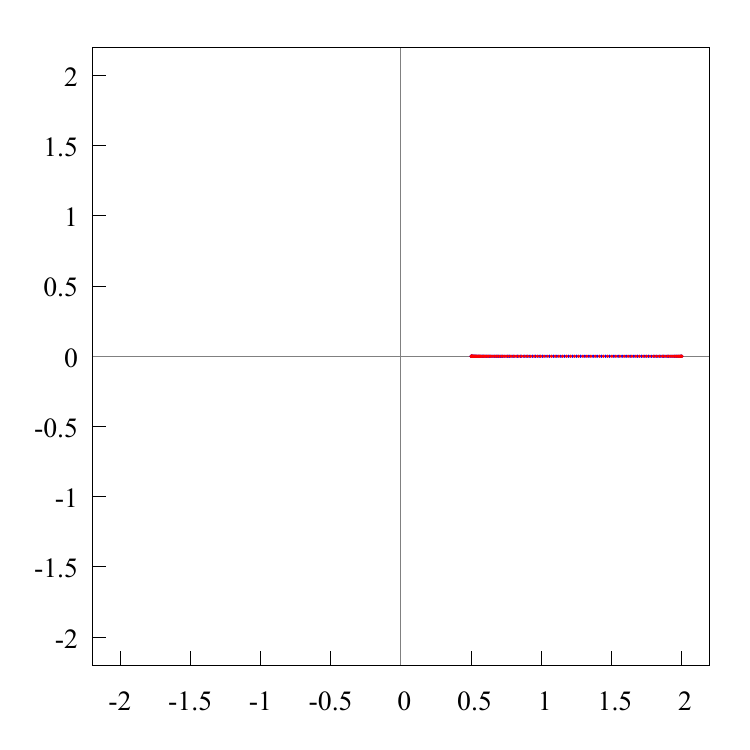}}
\vskip-6mm
\caption{The distribution of zeros and poles of the two-point Pad\'{e} approximant
(at points $z=0$ and $z=\infty$) of the function
$f(z)=\sqrt{(z-a)/(z-b)}$, $a=.5$, $b=2$.
Here are selected two ``same'' branches of the function $f$:
$f_0=\sqrt{(z-.5)/(z-2)}$, $f_\infty=\sqrt{(z-.5)/(z-2)}$.
The interval $[1/2,2]$ is the Buslaev compact.
}
\label{Fig_bus210a_3000_90_full}
\end{figure}

\begin{figure}[!ht]
\centerline{
\includegraphics[width=15cm,height=15cm]{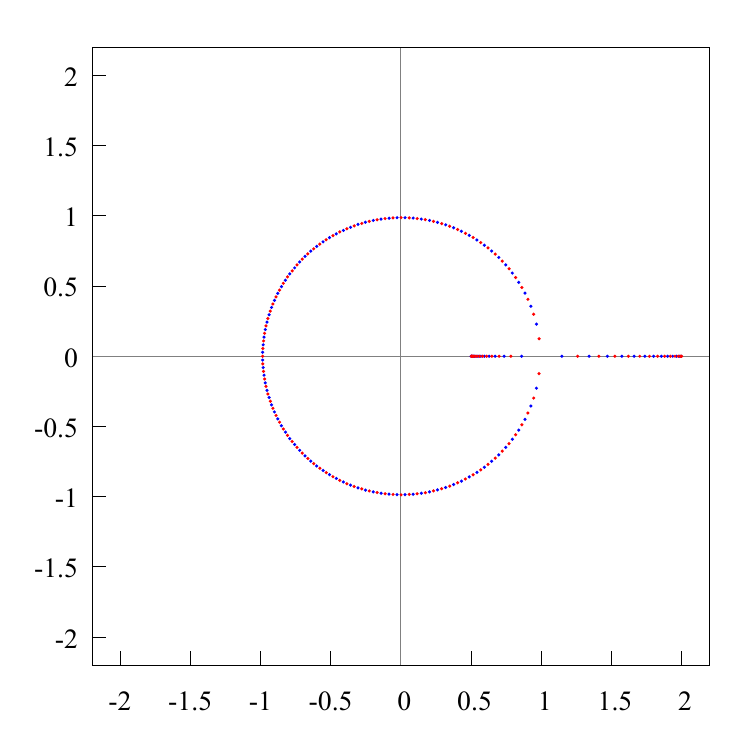}}
\vskip-6mm
\caption{The distribution of zeros and poles of the two-point Pad\'{e} approximant
(at points $z=0$ and $z=\infty$) of the function
$f(z)=\sqrt{(z-a)/(z-b)}$, $a=.5$, $b=2$.
Here are selected two ``different'' branches of the function $f$:
$f_0(z)=\sqrt{(z-.5)/(z-2)}$, $f_\infty(z)=-\sqrt{(z-.5)/(z-2)}$.
There is one Chebotarev point $v=\sqrt{ab}=1$ of zero density
on the Buslaev compact.
From this can be seen, that this point cannot be calculated
using of the two-point Pad\'{e} approximants, comp. \ref{Fig_pade10_2500_130_blu};
contrary to the case with classical Pad\'{e} approximants
(see fig. \ref{Fig_pade10_2500_130_red}, \ref{Fig_pade10_2500_130})
}
\label{Fig_bus210b_4000_120_full}
\end{figure}

For a better exposition of the new numerical results, we state
at the beginning some well known facts about the asymptotical
behavior of the zeros and poles of the classical Pad\'{e} approximants,
i.e. zeros of the Pad\'{e} polynomials (see \eqref{2}), and also for
two-point Pad\'{e} approximants \cite{Bus13}.
These results and their pictures (see
fig. \ref{Fig_pade10_2500_131}--\ref{Fig_pade103_5000_266_full} and
fig. \ref{Fig_bus210a_3000_90_full}--\ref{Fig_bus205c_4000_120_full})
are needed for the analysis of the numerical results,
connected with the behavior of the zeros of the Hermite--Pad\'{e} polynomials.

\section{Main results}\label{s2}

The main empirical results obtained in the paper are presented herein.

\subsection{Angelesco system}\label{s2s1}
\begin{figure}[!ht]
\centerline{
\includegraphics[width=15cm,height=15cm]{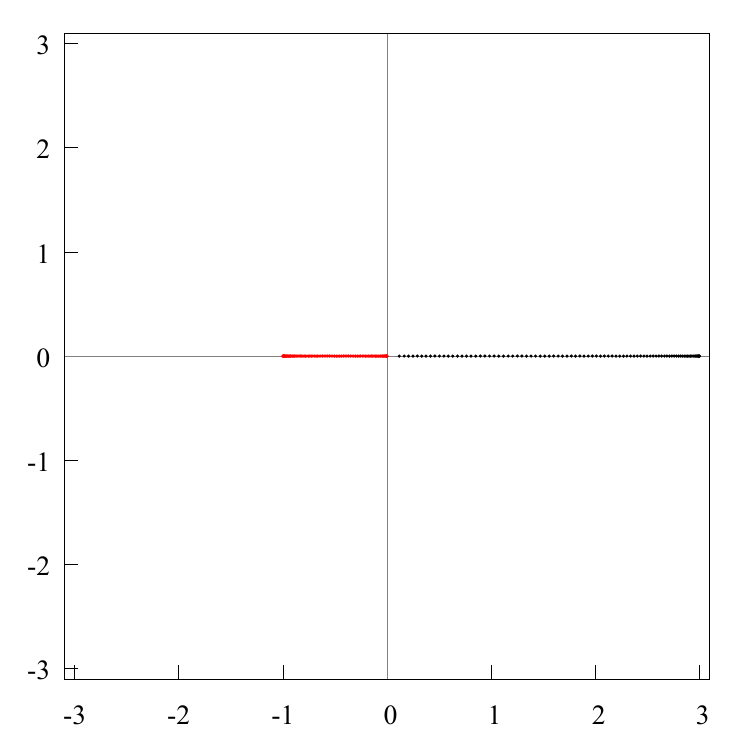}}
\vskip-6mm
\caption{The distribution of the zeros of the Hermite-Pad\'{e} polynomial $Q_{90,1}$ (red points)
and $Q_{90,2}$ (black points) for the set of three functions $[1,f_1,f_2]$, where
$f_1(z)=\sqrt{(z+1)/z}$, $f_2=\sqrt{(z-3)/z}$ (we select such a root branch,
that $\sqrt{1\,}=1$). The functions $f_1,f_2$ are of Markov type:
$f_1(z)=\myh\sigma_1(z)$, $f_2(z)=\myh\sigma_2(z)$, where
$\sigma'_1(x)=\sqrt{(1+x)/(-x)}/\pi$, $x\in(-1,0)$,
$\Delta_1=\supp{\sigma_1}=[-1,0]$,
$\sigma'_2(x)=\sqrt{(3-x)/x}/\pi$, $x\in(0,3)$,
$\Delta_2=\supp{\sigma_2}=[0,3]$. It is known \cite{GoRa81}, \cite{Apt08}, that
the zeros of the polynomials $Q_{n,1}$ and $Q_{n,2}$, when taken the limit,
are distributed on the support of the equilibrium measures
$S_1=\lambda_1$ and $S_2=\lambda_2$ and in accordance with their densities
$\lambda_1'(x)$ and $\lambda_2'(x)$. The figure clearly shows, that
$S_1=\Delta_1$, but $S_2=[a_2^*,3]\subsetneq\Delta_2$, because
$|\Delta_1|<|\Delta_2|$. The point $a_2^*\in(0,3)$ is calculated by the formula \eqref{eq.kal84},
obtained by V. A. Kalyagin \cite{Kal84} (see also \cite{Nut84}, \cite{ApKa86}).
The density $\lambda_1'(x)$ behaves, in the neighborhoods of the points $a_1$ and $b_1$,
like $(x-a_1)^{-1/2}$ and $(b_1-x)^{-1/2}$, respectively, and the density
$\lambda_2'(x)$ behaves, in the neighborhood of the point $a_2$, like $(x-a^*_1)^{1/2}$ 
and in the neighborhood of the point $b_2$ like $(b_2-x)^{-1/2}$.
}
\label{Fig_mar_1_2000_90_rdbk}
\end{figure}

The first part of the numerical results of the behavior of the
Hermite--Pad\'{e} polynomials for a collection of three functions
$[1,f_1,f_2]$ refers to the case when each of the functions $f_1$ and
$f_2$ has a pair of branch points of second order and the sets of the
branch points of the functions $f_1$ and $f_2$ do not intersect. Thus,
the pair of functions $f_1,f_2$ forms an Angelesco system (see
\eqref{3}). About the main properties of the Angelesco system, see
\cite{Nut84}, also \cite{GoRa81}. In \cite{GoRa81}, the first results
of general character about the convergence of the Hermite--Pad\'{e}
approximants of second kind for the collection of functions
$[1,f_1,\dots,f_{m-1}]$, $m\ge3$ were obtained, with the functions
$f_1,\dots,f_{m-1}$ creating an Angelesco system and
$f_j=\myh{\sigma}_j$ being Markov functions with non-intersecting
supports $\Delta_j=\supp{\sigma_j}$, $j=1,2,\dots,m-1$, along the real
axis, and containing a finite number of intervals. In \cite{GoRa81}, A.
A. Gonchar and E.~A.~Rakhmanov found a new property, which they called
{\it pushing} of the support of the equilibrium measures (see also
\cite{Nut84}, \cite{ApKa86}, \cite{Apt08}). For $m=3$ this method looks
as follows. Let the support sets $\Delta_1,\Delta_2$ of the measures
$\sigma_1,\sigma_2$ are the non-intersecting closed intervals
$\Delta_1=[a_1,b_1]$, $\Delta_2=[a_2,b_2]$ of different lengths and
contained in the real axis; for definiteness we suppose that
$|\Delta_1|<|\Delta_2|$ and that $\Delta_2$ be on the right of
$\Delta_1$. It turns out that in the case when the intervals are close
enough, the support $F_2\subset\Delta_2$ of the equilibrium measure
$\lambda_2$ is the interval $[a^*_2,b_2]$, where $a^*_2\in(a_2,b_2)$
(see \eqref{eq.kal84} and fig. \ref{Fig_mar_1_2000_90_rdbk}). Thus, the
equilibrium measures $\lambda_1,\lambda_2$ are absolutely continuous
with respect to the normalized linear Lebesgue measure, the density
$\lambda_1'(x)$ behaves, in the neighborhood of the points $a_1$ and
$b_1$, like $(x-a_1)^{-1/2}$ and $(b_1-x)^{-1/2}$, respectively, and
the density $\lambda_2'(x)$ behaves, in the neighborhood of the point
$a_2$, like $(x-a^*_1)^{1/2}$ and like $(b_2-x)^{-1/2}$ around $b_2$.
Thus, under a specific mutual positions the smaller interval pushes the
support of the equilibrium measure inside the larger interval; this
does not happen with the support of the equilibrium measure for a
smaller interval (see fig. \ref{Fig_mar_1_2000_90_rdbk}). When
$[a_1,b_1]=[-a,0]$, $[a_2,b_2]=[0,1]$, where $a\in(0,1)$, the point
$a_2^*$ is calculated by the following formula, found by V.~A.~Kalyagin
\cite{Kal84} (see also \cite[p. 5.3, formula (5.3.18)]{Nut84},
\cite{ApKa86}):
\begin{equation}
\label{eq.kal84}
a_2^*=\frac{(1-a)^3}{9(a^2-a+1)}.
\end{equation}

The authors have found, by experiments, a new property called {\it mutual
pushing} of the supports of the equilibrium measures, in the case when the functions $f_1$ and $f_2$
have a pair of branch points $a_1,b_1\in\CC$ and $a_2,b_2\in\CC$,
$\{a_1,b_1\}\cap\{a_2,b_2\}=\varnothing$, located not on one line,
but on two parallel lines, and the intervals
$\Delta_1=[a_1,b_1]$ and $\Delta_2=[a_2,b_2]$ have different lengths, i.e.
$|\Delta_1|>|\Delta_2|$ (see fig.
 \ref{Fig_ang3(14)_5000_180_rdbk}).

\begin{figure}[!ht]
\centerline{
\includegraphics[width=15cm,height=15cm]{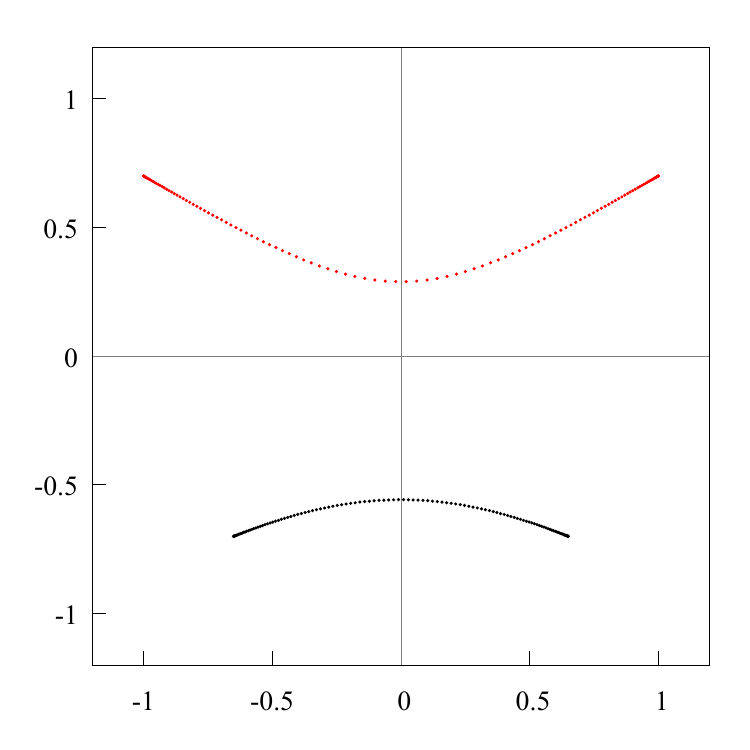}}
\caption{The distribution of the zeros of the Hermite-Pad\'{e} polynomial $Q_{120,1}$ (red points),
$Q_{120,2}$ (black points) in the plane $\CC_z$ for two functions
$f_1(z)=1/\sqrt{(z-(-1.0+i\cdot0.7))(z-(1.0+i\cdot0.7))}$,
$f_2(z)=1/\sqrt{(z-(-.65-i\cdot0.7))(z-(.65-i\cdot0.7))}$.
Since the branch points $a_1=-1.0+i\cdot0.7,b_1=1.0+i\cdot0.7$ and
$a_2=-0.65+i\cdot0.7,b_2=0.65+i\cdot0.7$ of the functions are
far enough from each other, there is no collision of the supports
of the equilibrium measures and the supports $\lambda_1$ and $\lambda_2$
are two non-intersecting arcs.
The measures $\lambda_1$ and $\lambda_2$ are absolutely continuous with respect
to the length of the arc $|dz|$, and their densities $\lambda_j'$, $j=1,2$,
behave, in the neighborhoods of the branch points $a_j,b_j$, like Chebyshev measures,
i.e. $\sim|z-a_j|^{-1/2}$ and $\sim|z-b_j|^{-1/2}$, respectively.
It is obvious that the extremal compacts $F_1$ and $F_2$ are converging to each other,
and the zeros of the polynomials $Q_{120,j}$, which are onto $F_j$,
are diverging from each other and from the branch points $a_j,b_j$, $j=1,2$.
The zeros of the polynomial $Q_{120,0}$ (blue points, see
fig. \ref{Fig_ang3(1)_5000_120_full}) create a third extremal compact
$F_0$, which separates the compacts $F_1$ and $F_2$.
}
\label{Fig_ang3(1)_5000_120_rdbk}
\end{figure}

\begin{figure}[!ht]
\centerline{
\includegraphics[width=15cm,height=15cm]{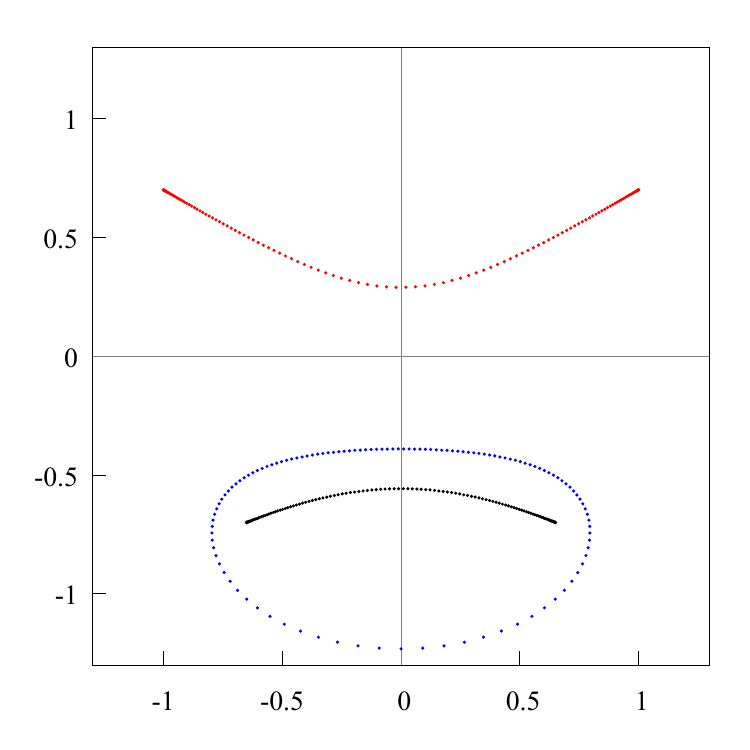}}
\caption{The distribution of the Hermite-Pad\'{e} polynomials $Q_{120,0}$ (blue points),
$Q_{120,1}$ (red points), $Q_{120,2}$ (black points) in the plane $\CC_z$
for two functions
$f_1(z)=1/\sqrt{(z-(-1.0+i\cdot0.7))(z-(1.0+i\cdot0.7))}$,
$f_2(z)=1/\sqrt{(z-(-.65-i\cdot0.7))(z-(.65-i\cdot0.7))}$.
The red points are distributed onto the support of the equilibrium measure
$\lambda_1$, when taking the limit, accordingly to its density $\lambda'_1$.
The black points are distributed onto the support of the equilibrium measure
$\lambda_2$, when taking the limit, accordingly to its density $\lambda'_2$.
Same applies for the blue points.
Both measures $\lambda_1$ and $\lambda_2$ are absolutely continuous
on the corresponding extremal curves $F_1=\supp{\lambda_1}$ and $F_2=\supp{\lambda_2}$
with respect to the length of the arc $ds$, and their densities
behave like the Chebyshev measures $|z-a_j|^{-1/2},|z-b_j|^{-1/2}$, $j=1,2$
around the endpoints of the curves.
The blue points (zeros of the polynomial $Q_{120,0}$) separate the red points
(zeros of the polynomial $Q_{120,1}$) from the black points (zeros of the polynomial $Q_{120,2}$).
}
\label{Fig_ang3(1)_5000_120_full}
\end{figure}

\subsection*{Fig. \ref{Fig_ang3(1)_5000_120_rdbk}--\ref{Fig_ang3(1)_5000_120_full}}
First, when these lines are far enough from each other, there is no
collision of the support of the equilibrium measures and the supports
of $\lambda_1$ and $\lambda_2$ are two non-intersecting arcs,
respectively\footnote{It is well known that in the case of an Angelesco
system these two arcs in a sense are ``{\it attracted}'' to each other,
and in the case of a Nikishin system are they ``{\it repelled}'' from
each other (see \cite{Apt08}, also \cite{Nut84}, \cite{NuTr87},
\cite{GoRa85}, \cite{Gon03}).}, see \cite{Nut84}, \cite{NuTr87},
\cite{Apt08}. The measures $\lambda_1$ and $\lambda_2$ are absolutely
continuous with respect to the length of the arc $|dz|$, and their
densities $\lambda_j'$, $j=1,2$, in the neighborhoods of the branch
points $a_j,b_j$, behave like Chebyshev measures, i.e.
$\sim|z-a_j|^{-1/2}$ and $\sim|z-b_j|^{-1/2}$, respectively. This can
be seen very well on figure \ref{Fig_ang3(1)_5000_120_rdbk}, where the
red points are the zeros of the Hermite--Pad\'{e} polynomial
$Q_{120,1}$ and the black points are the zeros of $Q_{120,2}$. It is
obvious that the extremal compact sets $F_1$ and $F_2$ are attracted to
each other, and the zeros of the polynomials $Q_{120,j}$, which are
onto $F_j$, are repelled from each other and from the branch points
$a_j,b_j$, $j=1,2$. The zeros of the polynomial $Q_{120,0}$ (blue
points, see fig. \ref{Fig_ang3(1)_5000_120_full}) form a third extremal
compact $F_0$, which separates the compact sets $F_1$ and $F_2$. The
distribution of the zeros of the polynomial $Q_{n,0}$, when
$n\to\infty$, is described by the third extremal measure $\lambda_0$,
$\supp\lambda_0=F_0$. Thus, the following fact is true (see
\cite{Sta88}). Let $U_n(z):=\max\{\log|Q_{n,1}(z)|,\log|Q_{n,2}(z)|\}$,
$z\in\CC$. Then $U_n$ is subharmonic in $\CC$. Therefore,
$U_n(z)=-U^{\mu_n}(z)$, where $\mu_n$ is measure, $|\mu_n|\leq n$,
$U^{\mu_n}(z)=-\int\log|z-\zeta|\,d\mu_n(\zeta)$ is the logarithmic
potential with respect to $\mu_n$. For an arbitrary polynomial
$Q\in\CC[z]$ define the measure
$$
\chi(Q)=\sum_{\zeta:Q(\zeta)=0}\delta_\zeta ,
$$
which counts the number of zeros of the polynomial $Q$. Then
\begin{equation}
\label{eq_Sta88}
\lim_{n\to\infty}\frac1n\mu_n=
\lim_{n\to\infty}\frac1n\chi(Q_{n,0})=\lambda_0
\end{equation}
(the convergence in \eqref{eq_Sta88} is understood as weak convergence of measures).

\begin{figure}[!ht]
\centerline{
\includegraphics[width=15cm,height=15cm]{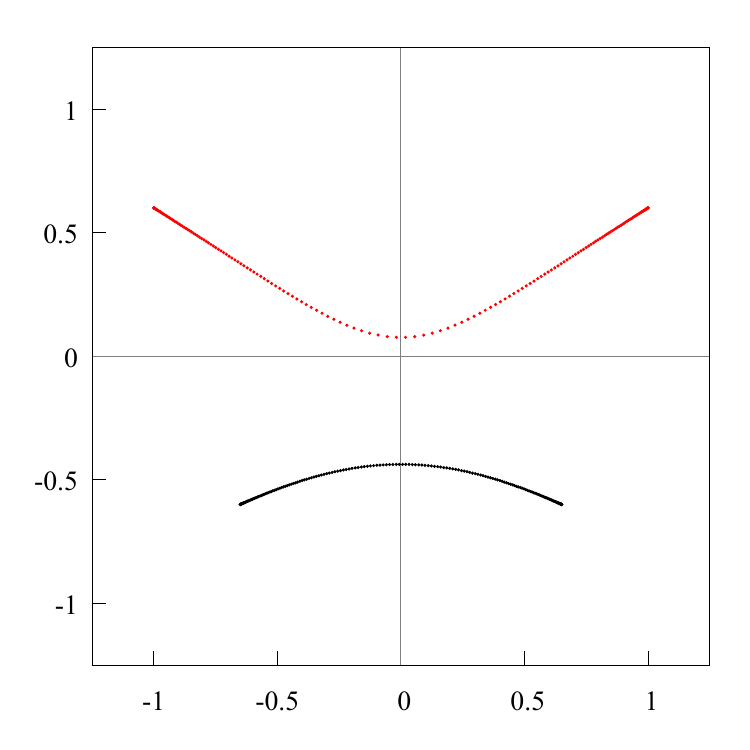}}
\caption{The distribution of the Hermite-Pad\'{e} polynomials $Q_{180,1}$ (red points),
$Q_{180,2}$ (black points) for two functions
$f_1(z)=1/\sqrt{(z-(-1.0+i\cdot0.6))(z-(1.0+i\cdot0.6))}$,
$f_2(z)=1/\sqrt{(z-(-.65-i\cdot0.6))(z-(.65-i\cdot0.6))}$.
The branch points have come closer to each other,
however they are still far enough from each other
and there is no collision of the supports of the equilibrium measures.
It is clearly seen, that the upper extremal compact $F_1$
has strongly bent towards the lower extremal compact $F_2$.
}
\label{Fig_ang3(2)_5000_180_rdbk}
\end{figure}

\begin{figure}[!ht]
\centerline{
\includegraphics[width=15cm,height=15cm]{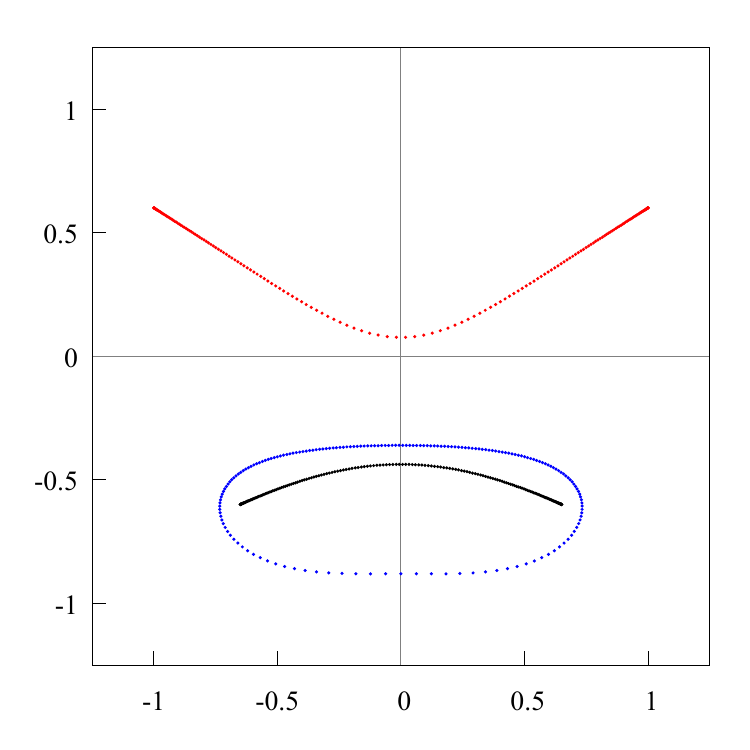}}
\caption{The distribution of the Hermite-Pad\'{e} polynomials $Q_{180,0}$ (blue points),
$Q_{180,1}$ (red points), $Q_{180,2}$ (black points) for two functions
$f_1(z)=\sqrt{(1.0-(-1.0+i\cdot0.6)\zeta)(1.0-(1.0+i\cdot0.6)\zeta)}$,
$f_2(z)=\sqrt{(1.0-(-.65-i\cdot0.6)\zeta)(1.0-(.65-i\cdot0.6)\zeta)}$.
The branch points have come closer to each other,
however they are still far enough from each other
and there is no collision of the supports of the equilibrium measures.
It is clearly seen, that the upper extremal compact $F_1$
has strongly bent towards the lower extremal compact $F_2$.
The third extremal compact $F_0$, as before, separates $F_1$ and $F_2$.
}
\label{Fig_ang3(2)_5000_180_full}
\end{figure}

\subsection*{Fig. \ref{Fig_ang3(2)_5000_180_rdbk}--\ref{Fig_ang3(2)_5000_180_full}}
Further convergence of the parallel lines along the imaginary axis leads to the following.
While the branch points are far enough from each other,
collision of the equilibrium measures does not occur.
However, on fig. \ref{Fig_ang3(2)_5000_180_rdbk} it is clearly seen,
that the upper extremal compact set $F_1$ has strongly curved towards
the lower extremal compact set $F_2$.
The third extremal compact set $F_0$, as before, separates $F_1$ and $F_2$.

\begin{figure}[!ht]
\centerline{
\includegraphics[width=15cm,height=15cm]{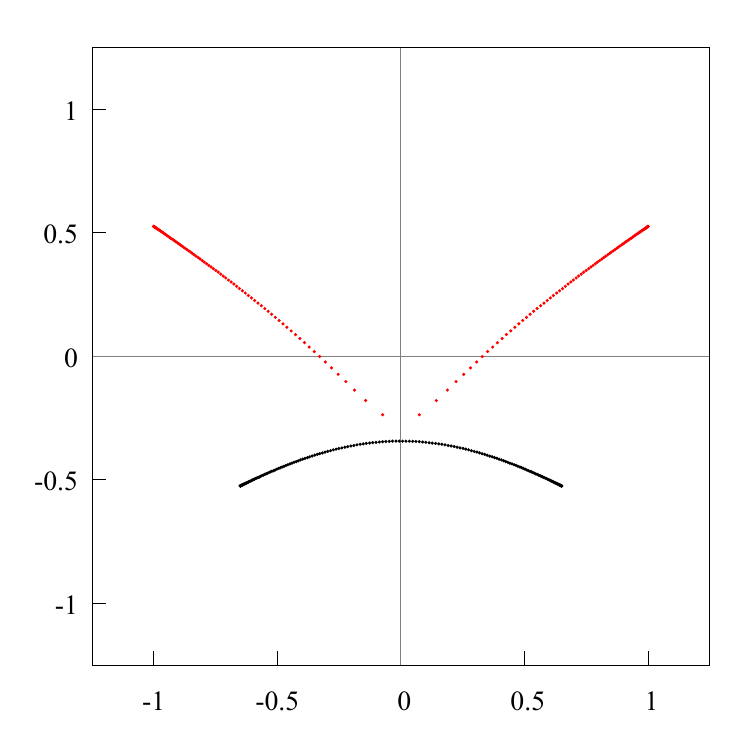}}
\caption{The distribution of the Hermite-Pad\'{e} polynomials $Q_{180,1}$ (red points),
$Q_{180,2}$ (black points) for two functions
$f_1(z)=1/\sqrt{(z-(-1.0+i\cdot0.525))(z-(1.0+i\cdot0.525))}$,
$f_2(z)=1/\sqrt{(z-(-.65-i\cdot0.525))(z-(.65-i\cdot0.525))}$.
It is clearly seen, that the upper extremal compact $F_1$
has even strongly bent towards the lower extremal compact $F_2$.
The third extremal compact $F_0$, as before, separates $F_1$ and $F_2$.
The support of the equilibrium measure of the upper extremal compact
$F_1$ starts to break down, while the second compact $F_2$ has hardly changed.
}
\label{Fig_ang3(8)_5000_180_rdbk}
\end{figure}

\begin{figure}[!ht]
\centerline{
\includegraphics[width=15cm,height=15cm]{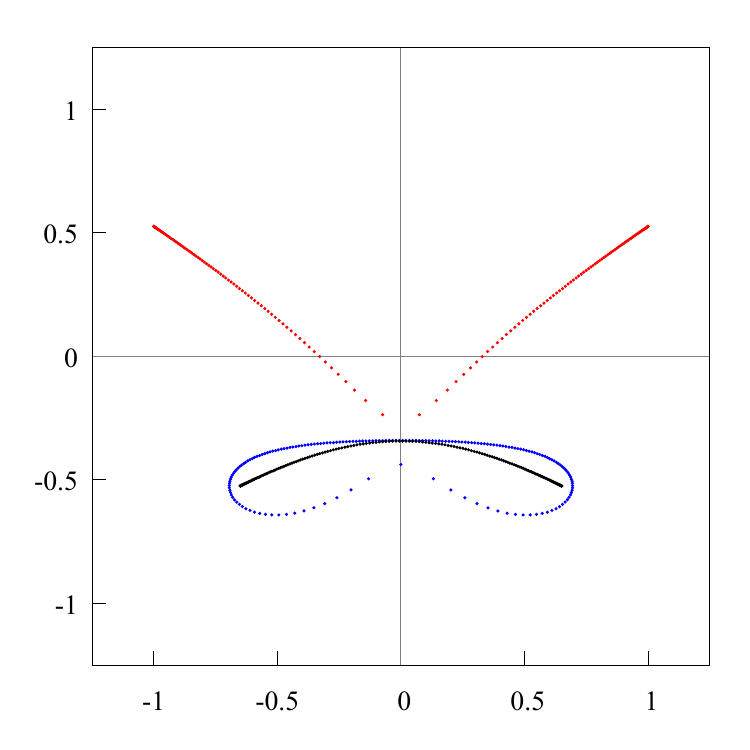}}
\caption{The distribution of the Hermite-Pad\'{e} polynomials
$Q_{180,0}$ (blue points),
$Q_{180,1}$ (red points),
$Q_{180,2}$ (black points) for two functions
$f_1(z)=1/\sqrt{(z-(-1.0+i\cdot0.525))(z-(1.0+i\cdot0.525))}$,
$f_2(z)=1/\sqrt{(z-(-.65-i\cdot0.525))(z-(.65-i\cdot0.525))}$.
It is clearly seen, that the upper extremal compact $F_1$
has even strongly bent towards the lower extremal compact $F_2$.
The third extremal compact $F_0$, as before, separates $F_1$ and $F_2$.
The support of the equilibrium measure of the upper extremal compact
$F_1$ starts to break down, while the second compact $F_2$ has hardly changed.
The third extremal compact $F_0$, as before,
separates the other two compacts from each other,
but now it touches the second compact $F_2$.
}
\label{Fig_ang3(8)_5000_180_full}
\end{figure}

\subsection*{Fig. \ref{Fig_ang3(8)_5000_180_rdbk}--\ref{Fig_ang3(8)_5000_180_full}}
Then, under further convergence of the branch points,
the upper extremal compact set $F_1$ has even strongly curved towards
the lower extremal compact set $F_2$, the support of the equilibrium measure
of the upper compact set $F_1$ starts to break down,
the second lower compact set almost does not change;
see fig. \ref{Fig_ang3(8)_5000_180_rdbk}.
The third extremal compact set $F_0$, as before,
separates the other two compact sets from each other,
but now it touches the second compact set $F_2$.

\begin{figure}[!ht]
\centerline{
\includegraphics[width=15cm,height=15cm]{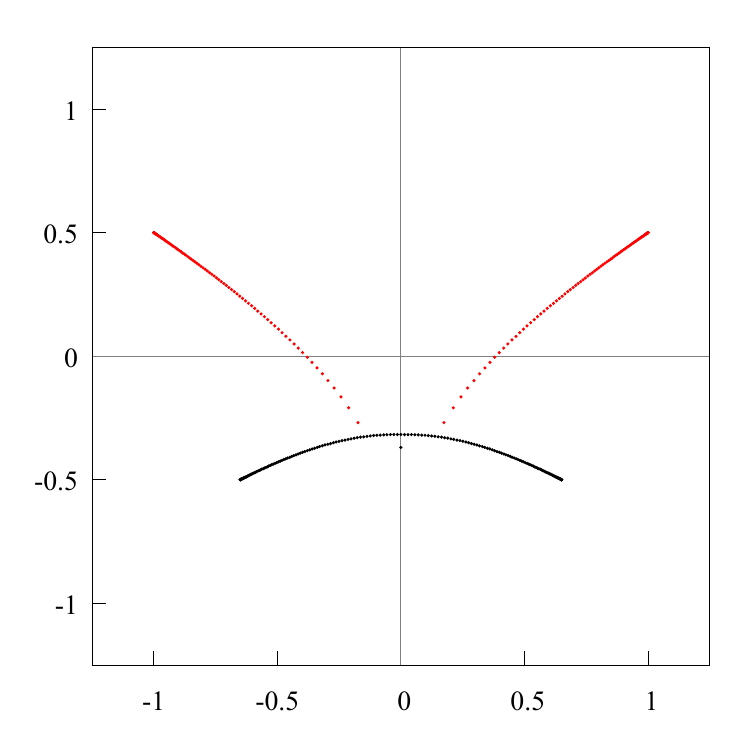}}
\caption{The distribution of the Hermite-Pad\'{e} polynomials $Q_{180,1}$ (red points),
$Q_{180,2}$ (black points) for two functions
$f_1(z)=1/\sqrt{(z-(-1.0+i\cdot0.5))(z-(1.0+i\cdot0.5))}$,
$f_2(z)=1/\sqrt{(z-(-.65-i\cdot0.5))(z-(.65-i\cdot0.5))}$.
It is clearly seen, that under this position of the pair of branch points,
the support of the equilibrium measure of the upper extremal compact
$F_1$ breaks down, while the second compact $F_2$ has hardly changed.
The red points break down the support of the equilibrium measure $\lambda_1$
of the compact $F_1$ on two arcs.
}
\label{Fig_ang3(7)_5000_180_rdbk}
\end{figure}

\begin{figure}[!ht]
\centerline{
\includegraphics[width=15cm,height=15cm]{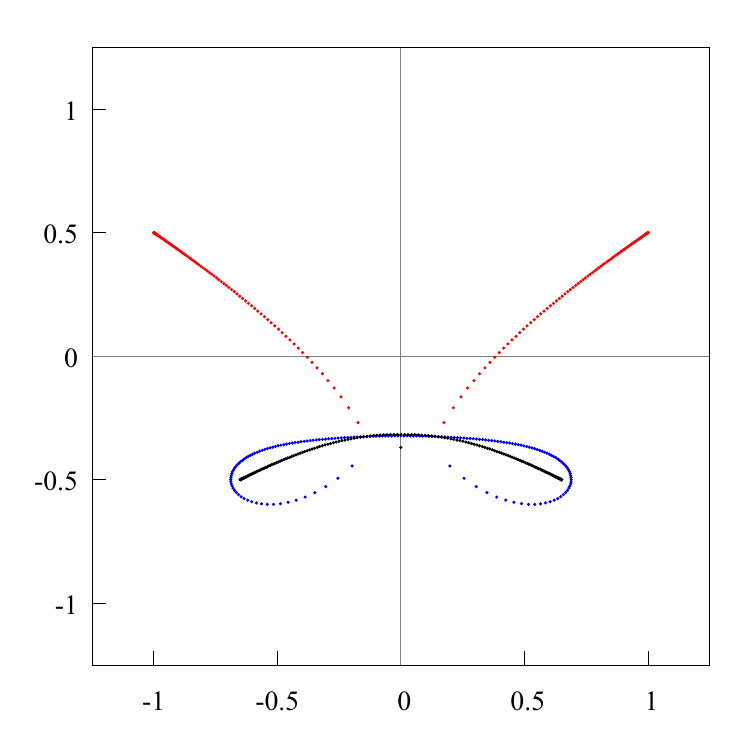}}
\caption{The distribution of the Hermite-Pad\'{e} polynomials
$Q_{180,0}$ (blue points),
$Q_{180,1}$ (red points),
$Q_{180,2}$ (black points) for two functions
$f_1(z)=1/\sqrt{(z-(-1.0+i\cdot0.5))(z-(1.0+i\cdot0.5))}$,
$f_2(z)=1/\sqrt{(z-(-.65-i\cdot0.5))(z-(.65-i\cdot0.5))}$.
It is clearly seen, that under this position of the pair of branch points,
the support of the equilibrium measure of the upper extremal compact
$F_1$ breaks down, while the second compact $F_2$ has hardly changed.
The third extremal compact $F_0$, as before,
separates the other two compacts from each other,
and touches the second compact $F_2$
}
\label{Fig_ang3(7)_5000_180_full}
\end{figure}

\subsection*{Fig. \ref{Fig_ang3(7)_5000_180_rdbk}--\ref{Fig_ang3(7)_5000_180_full}}
Finally, under certain relative positions of the pair of branch points,
the support of the equilibrium measure of the upper extremal compact
$F_1$ breaks down, while the second compact set $F_2$ has hardly changed,
see fig. \ref{Fig_ang3(7)_5000_180_rdbk}.
The third extremal compact set $F_0$, as before,
separates the other two compact sets from each other,
and touches the second compact set $F_2$,
see fig. \ref{Fig_ang3(7)_5000_180_full}.

\begin{figure}[!ht]
\centerline{
\includegraphics[width=15cm,height=15cm]{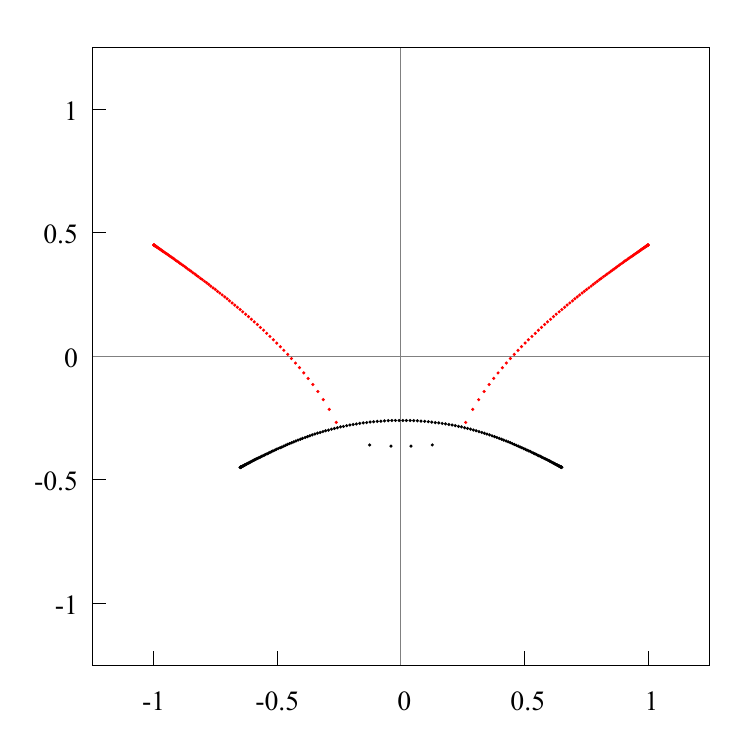}}
\caption{The distribution of the Hermite-Pad\'{e} polynomials $Q_{180,1}$ (red points),
$Q_{180,2}$ (black points) for two functions
$f_1(z)=1/\sqrt{(z-(-1.0+i\cdot0.45))(z-(1.0+i\cdot0.45))}$,
$f_2(z)=1/\sqrt{(z-(-.65-i\cdot0.45))(z-(.65-i\cdot0.45))}$.
It is clearly seen, that under this position of the pair of branch points,
the two arcs, which are the result of the breaking of the support of the measure $\lambda_1$,
have reached the second (lower) compact $F_2$. The second compact $F_2$ has started to change:
from the total set of black points (zeros of the polynomial $Q_{180,2}$) several points stand out,
which started to form another component.
Thus, a second component of the support of the equilibrium measure $\lambda_2$
started to form, i.e. the support of the equilibrium measure $\lambda_2$
started breaking down on two arcs.
}
\label{Fig_ang3(9)_5000_180_rdbk}
\end{figure}

\begin{figure}[!ht]
\centerline{
\includegraphics[width=15cm,height=15cm]{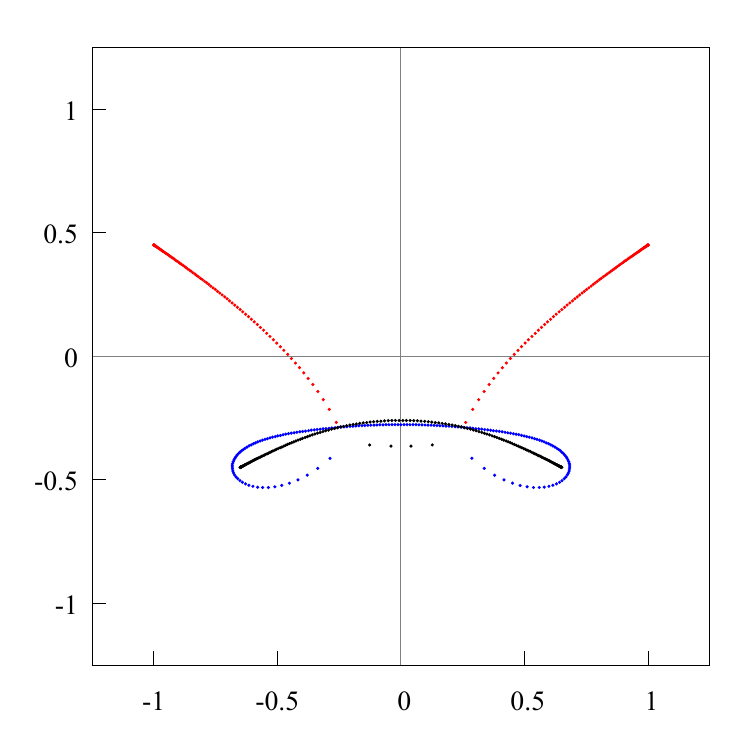}}
\caption{The distribution of the Hermite-Pad\'{e} polynomials $Q_{180,0}$ (blue points),
$Q_{180,1}$ (red points), $Q_{180,2}$ (black points) for two functions
$f_1(z)=1/\sqrt{(z-(-1.0+i\cdot0.45))(z-(1.0+i\cdot0.45))}$,
$f_2(z)=1/\sqrt{(z-(-.65-i\cdot0.45))(z-(.65-i\cdot0.45))}$.
It is clearly seen, that under this position of the pair of branch points,
the two arcs, which are the result of the breaking of the support of the measure $\lambda_1$,
have reached the second (lower) compact $F_2$. The second compact $F_2$ has started to change:
from the total set of black points (zeros of the polynomial $Q_{180,2}$) several points stand out,
which started to form another component.
Thus, a second component of the support of the equilibrium measure $\lambda_2$
started to form, i.e. the support of the equilibrium measure $\lambda_2$
started breaking down on two arcs.
The third extremal compact $F_0$, as before, ``seeks'' to separate
the other two compacts from each other, but now each of the compacts $F_1$ and $F_2$
has by two components. It is clearly seen, that the compact $F_0$
now crosses the compact $F_2$.
}
\label{Fig_ang3(9)_5000_180_full}
\end{figure}

\subsection*{Fig. \ref{Fig_ang3(9)_5000_180_rdbk}--\ref{Fig_ang3(9)_5000_180_full}}
Under further convergence of the pair of branch points,
the two arcs, which are the result of the breaking of the support of the measure $\lambda_1$,
have reached the second (lower) compact set $F_2$. The second compact set $F_2$ has started to change:
from the total set of black points (zeros of the polynomial $Q_{180,2}$)
several points stand out,
which started to form another component.
Thus, a second component of the support of the equilibrium measure $\lambda_2$
started to form, i.e. the support of the equilibrium measure $\lambda_2$
started breaking down on two arcs; see fig. \ref{Fig_ang3(9)_5000_180_rdbk}.
The third extremal compact set $F_0$, as before, ``seeks'' to separate
the other two compact sets from each other, but now each of the compact sets $F_1$ and $F_2$
has two components. It is clearly seen, that the compact set $F_0$
now crosses the compact set $F_2$; see fig. \ref{Fig_ang3(9)_5000_180_full}.

\begin{figure}[!ht]
\centerline{
\includegraphics[width=15cm,height=15cm]{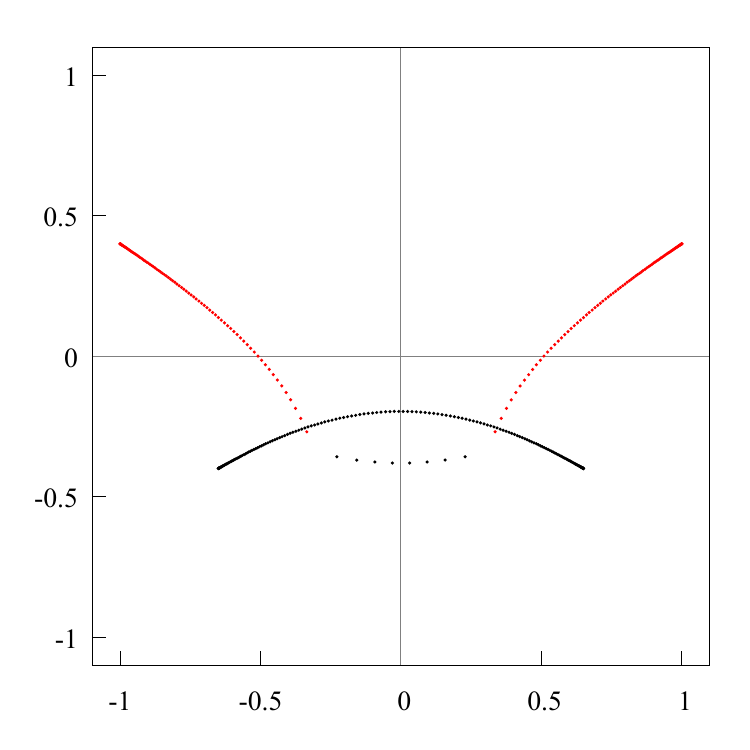}}
\caption{The distribution of the Hermite-Pad\'{e} polynomials $Q_{180,1}$ (red points),
$Q_{180,2}$ (black points) for two functions
$f_1(z)=1/\sqrt{(z-(-1.0+i\cdot0.4))(z-(1.0+i\cdot0.4))}$,
$f_2(z)=1/\sqrt{(z-(-.65-i\cdot0.4))(z-(.65-i\cdot0.4))}$.
It is clearly seen, that under this position of the pair of branch points,
the two arcs, which are the result of the breaking of the support of the measure $\lambda_1$,
cross the second compact $F_2$. The second compact $F_2$ continues to change:
from the total set of black points (zeros of the polynomial $Q_{180,2}$)
even more points stand out (than before), which form the second component of $F_2$.
Thus, the forming of the second component of the support of the equilibrium measure
$\lambda_2$ continues.
}
\label{Fig_ang3(14)_5000_180_rdbk}
\end{figure}

\begin{figure}[!ht]
\centerline{
\includegraphics[width=15cm,height=15cm]{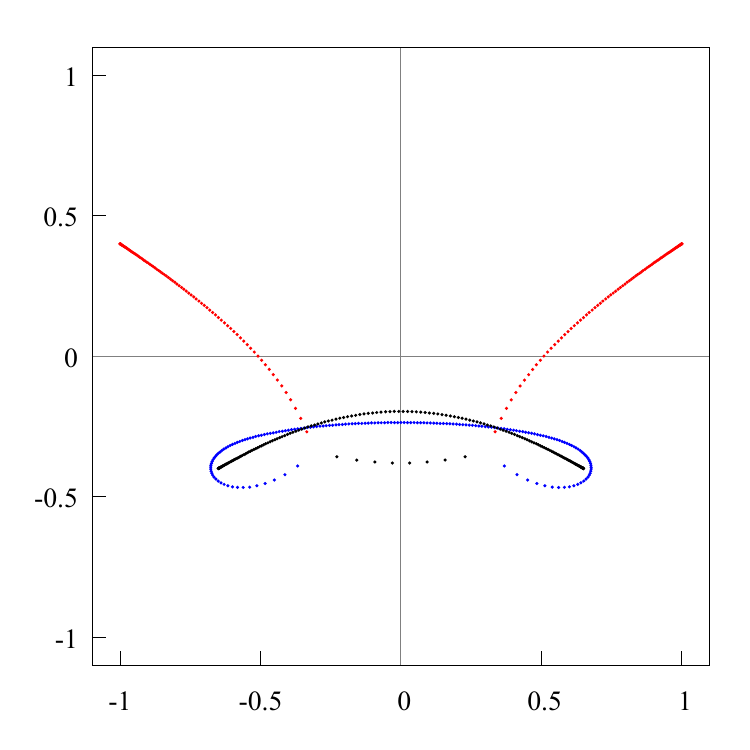}}
\caption{The distribution of the Hermite-Pad\'{e} polynomials
$Q_{180,0}$ (blue points),
$Q_{180,1}$ (red points),
$Q_{180,2}$ (black points) for two functions
$f_1(z)=1/\sqrt{(z-(-1.0+i\cdot0.4))(z-(1.0+i\cdot0.4))}$,
$f_2(z)=1/\sqrt{(z-(-.65-i\cdot0.4))(z-(.65-i\cdot0.4))}$.
It is clearly seen, that under this position of the pair of branch points,
the two arcs, which are the result of the breaking of the support of the measure $\lambda_1$,
cross the second compact $F_2$. The second compact $F_2$ continues to change:
from the total set of black points (zeros of the polynomial $Q_{180,2}$)
even more points stand out (than before), which form the second component of $F_2$.
Thus, the forming of the second component of the support of the equilibrium measure
$\lambda_2$ continues.
The third extremal compact $F_0$ crosses the compact $F_2$. As before,
it ``seeks'' to separate the other two compacts $F_1$ and $F_2$ from each other,
but now it ``fights'' with two components.
It is clearly seen, that at the junction of the red, black and blue points
appear two ``equilateral'' triangles with multicolored vertexes.
}
\label{Fig_ang3(14)_5000_180_full}
\end{figure}

\subsection*{Fig. \ref{Fig_ang3(14)_5000_180_rdbk}--\ref{Fig_ang3(14)_5000_180_full}}
Further, the two arcs, which are the result of the breaking of the support of the measure $\lambda_1$,
cross the second compact set $F_2$. The second compact set $F_2$ continues to change:
from the total set of black points (zeros of the polynomial $Q_{180,2}$)
even more points stand out (than before), which form the second component of $F_2$.
Thus, the forming of the second component of the support of the equilibrium measure
$\lambda_2$ continues; see fig. \ref{Fig_ang3(14)_5000_180_rdbk}.
The third extremal compact set $F_0$ crosses the compact set $F_2$. As before,
it ``seeks'' to separate the other two compact sets $F_1$ and $F_2$ from each other,
but now each of these compact sets has by two components.
It is clearly seen, that at the junction of the red, black and blue points
appear two equilateral triangles with multicolored vertexes,
see fig. \ref{Fig_ang3(14)_5000_180_full}.

\begin{figure}[!ht]
\centerline{
\includegraphics[width=15cm,height=15cm]{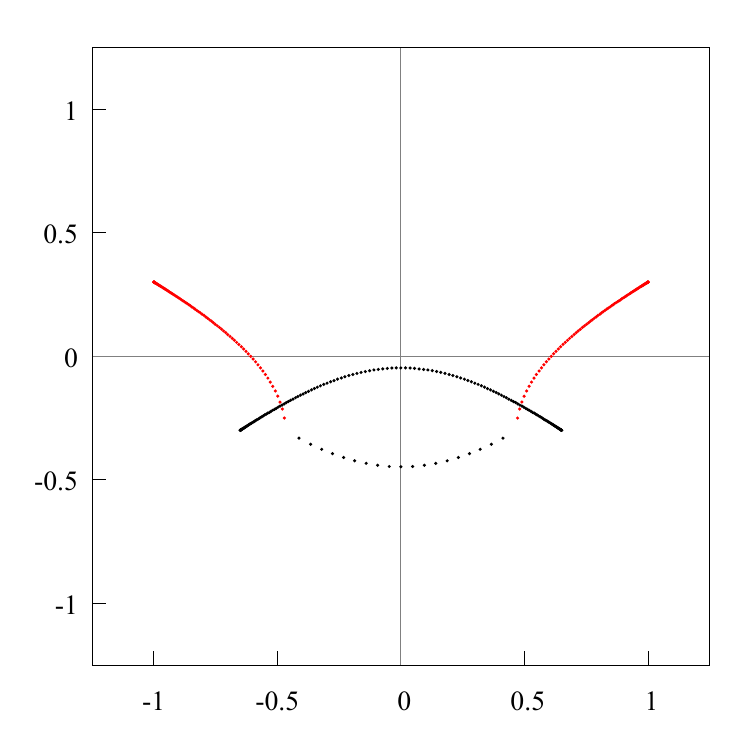}}
\caption{The distribution of the Hermite-Pad\'{e} polynomials $Q_{180,1}$ (red points),
$Q_{180,2}$ (black points) for two functions
$f_1(z)=1/\sqrt{(z-(-1.0+i\cdot0.3))(z-(1.0+i\cdot0.3))}$,
$f_2(z)=1/\sqrt{(z-(-.65-i\cdot0.3))(z-(.65-i\cdot0.3))}$.
It is clearly seen, that under this position of the pair of branch points,
the two arcs, which are the result of the breaking of the support of the measure $\lambda_1$,
even further cross the second compact $F_2$. The second compact $F_2$ continues to change:
from the total set of black points (zeros of the polynomial $Q_{180,2}$)
even more points stand out (even than before), which form the second component of $F_2$.
Thus, the forming of the second component of the support of the equilibrium measure
$\lambda_2$ continues.
}
\label{Fig_ang3(11)_5000_180_rdbk}
\end{figure}

\begin{figure}[!ht]
\centerline{
\includegraphics[width=15cm,height=15cm]{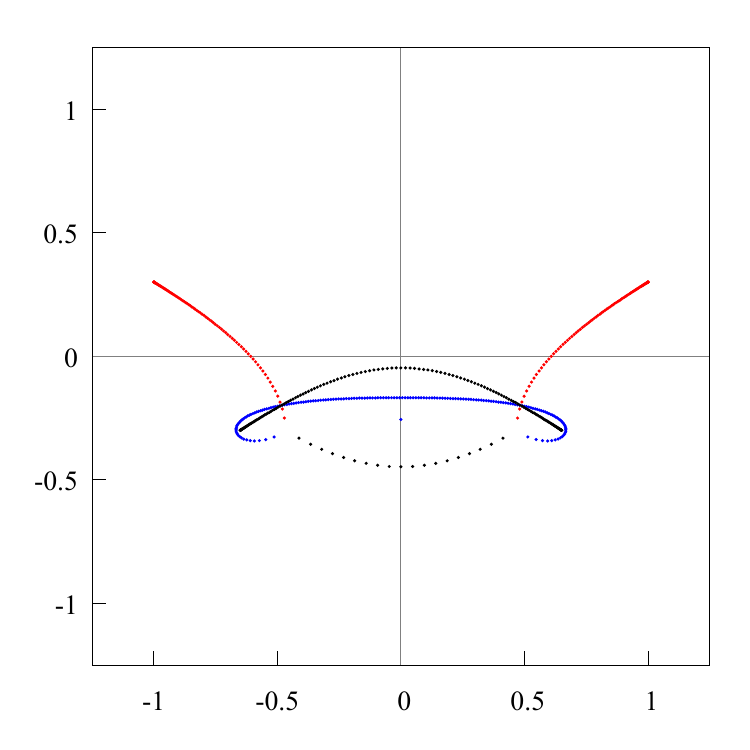}}
\vskip-10mm
\caption{The distribution of the Hermite-Pad\'{e} polynomials
$Q_{180,0}$ (blue points),
$Q_{180,1}$ (red points) and $Q_{180,2}$ (black points) for two functions
$f_1(z)=1/\sqrt{(z-(-1.0+i\cdot0.3))(z-(1.0+i\cdot0.3))}$,
$f_2(z)=1/\sqrt{(z-(-.65-i\cdot0.3))(z-(.65-i\cdot0.3))}$.
It is clearly seen, that under this position of the pair of branch points,
the two arcs, which are the result of the breaking of the support of the measure $\lambda_1$,
even further cross the second compact $F_2$. The second compact $F_2$ continues to change:
from the total set of black points (zeros of the polynomial $Q_{180,2}$)
even more points stand out (even than before), which form the second component of $F_2$.
Thus, the forming of the second component of the support of the equilibrium measure
$\lambda_2$ continues. It is clearly seen, that at the junction of the red, black and blue points
appear two ``equilateral'' triangles with multicolored vertexes.
By analogy with classical Pad\'{e} approximants and two-point Pad\'{e} approximants
(see fig. \ref{Fig_pade10_2500_130_blu} and \ref{Fig_bus210b_4000_120_full})
it is natural to assume, that the center of each triangle
has a Chebotarev point $v_1$, $v_2$ with zero density.
At the branch points $a_j,b_j$ the density of the measures $\lambda_1$
and $\lambda_2$ are proportional to $|z-a_j|^{-1/2},|z-b_j|^{-1/2}$, $j=1,2$,
respectively. There is a Froissart singlet (blue) on the imaginary axis.
}
\label{Fig_ang3(11)_5000_180_full}
\end{figure}

\subsection*{Fig. \ref{Fig_ang3(11)_5000_180_rdbk}--\ref{Fig_ang3(11)_5000_180_full}}
The two arcs, which are the result of the breaking of the support of the measure $\lambda_1$,
even further cross the second compact set $F_2$. The second compact set $F_2$ continues to change:
from the total set of black points (zeros of the polynomial $Q_{180,2}$)
even more points stand out (even than before), which form the second component of $F_2$.
Thus, the forming of the second component of the support of the equilibrium measure
$\lambda_2$ continues. It is clearly seen, that at the junction of the red, black and blue points
appear two equilateral triangles with multicolored vertexes.
By analogy with classical Pad\'{e} approximants and two-point Pad\'{e} approximants
(see fig. \ref{Fig_pade10_2500_130_blu} and \ref{Fig_bus210b_4000_120_full}),
it is natural to assume, that the center of each triangle
has a Chebotarev point $v_1$, $v_2$ with zero density.
At the branch points $a_j,b_j$ the density of the measures $\lambda_1$
and $\lambda_2$ are proportional to $|z-a_j|^{-1/2},|z-b_j|^{-1/2}$, $j=1,2$,
respectively. There is a Froissart singlet (blue) on the imaginary axis;
see fig. \ref{Fig_ang3(11)_5000_180_full}, \ref{Fig_ang3(11)_5000_180_rdbk}.

\begin{figure}[!ht]
\centerline{
\includegraphics[width=15cm,height=15cm]{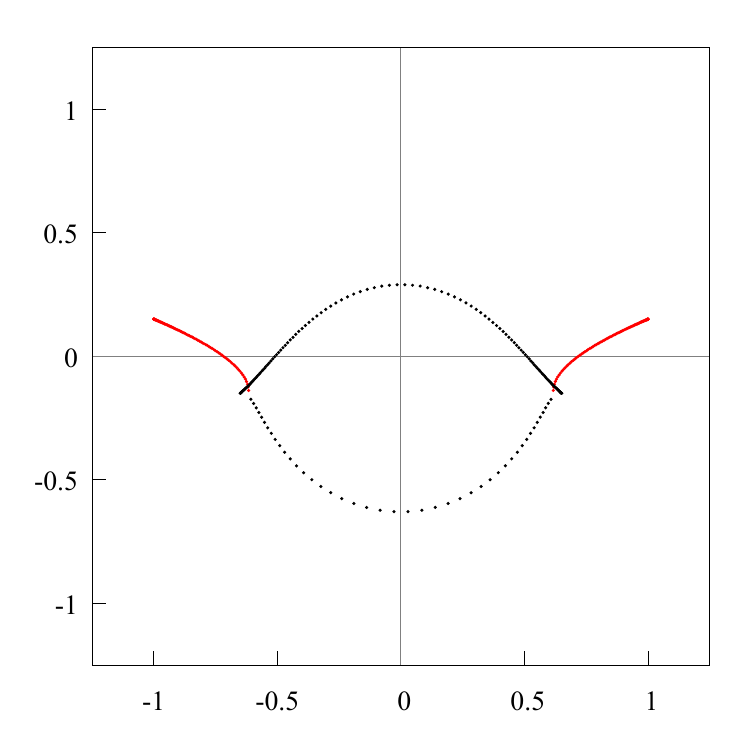}}
\caption{The distribution of the Hermite-Pad\'{e} polynomials $Q_{180,1}$ (red points),
$Q_{180,2}$ (black points) for two functions
$f_1(z)=1/\sqrt{(z-(-1.0+i\cdot0.15))(z-(1.0+i\cdot0.15))}$,
$f_2(z)=1/\sqrt{(z-(-.65-i\cdot0.15))(z-(.65-i\cdot0.15))}$.
It is clearly seen, that under this position of the pair of branch points,
the support $F_2$ of the equilibrium measure $\lambda_2$ is separated
on two practically equivalent arcs.
However, according to the distribution of the zeros of the polynomial $Q_{180,2}$,
the density $\lambda_2'$ of the equilibrium measures of each arc must be different.
On the upper arc it behaves like a Chebyshev measure, that is at the end points $a_2,b_2$
the density is proportional to $|z-a_2|^{-1/2}$ and $|z-b_2|^{-1/2}$, respectively.
The end points of the lower arc $v_1,v_2$ are the Chebotarev points
and their density is proportional to $|z-v_1|^{1/2}$ and $|z-v_2|^{1/2}$.
}
\label{Fig_ang3(13)_5000_180_rdbk}
\end{figure}

\begin{figure}[!ht]
\centerline{
\includegraphics[width=15cm,height=15cm]{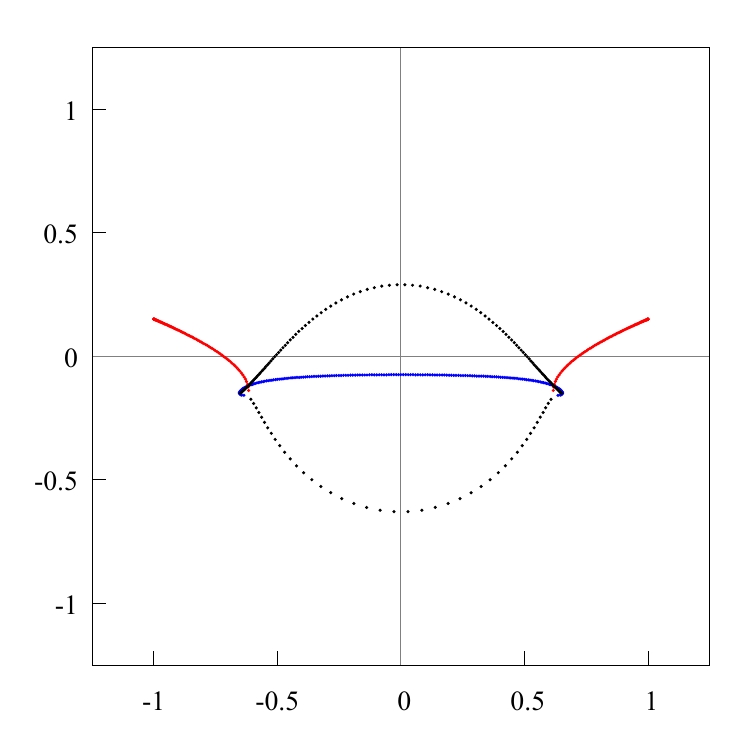}}
\caption{The distribution of the Hermite-Pad\'{e} polynomials
$Q_{180,0}$ (blue points),
$Q_{180,1}$ (red points),
$Q_{180,2}$ (black points) for two functions
$f_1(z)=1/\sqrt{(z-(-1.0+i\cdot0.15))(z-(1.0+i\cdot0.15))}$,
$f_2(z)=1/\sqrt{(z-(-.65-i\cdot0.15))(z-(.65-i\cdot0.15))}$.
It is clearly seen, that under this position of the pair of branch points,
the support $F_2$ of the equilibrium measure $\lambda_2$ is separated
on two practically equivalent arcs.
The branch points $a_2,b_2$ are on the upper arc,
the end points of the lower arc $v_1,v_2$ are the Chebotarev points
and their density is proportional to $|z-v_1|^{1/2}$ and $|z-v_2|^{1/2}$.
At the junction of the red, black and blue points
appeared two ``equilateral'' triangles with multicolored vertexes,
and the center of each has a Chebotarev point $v_1,v_2$.
}
\label{Fig_ang3(13)_5000_180_full}
\end{figure}

\subsection*{Fig. \ref{Fig_ang3(13)_5000_180_rdbk}--\ref{Fig_ang3(13)_5000_180_full}}
Finally, under certain relative positions of the pairs of branch points,
the support $F_2$ of the equilibrium measure $\lambda_2$ is separated
on two practically equivalent arcs.
However, according to the distribution of the zeros of the polynomial $Q_{180,2}$,
the density $\lambda_2'$ of the equilibrium measures of each arc must be different.
On the upper arc it behaves like a Chebyshev measure, that is at the end points $a_2,b_2$
the density is proportional to $|z-a_2|^{-1/2}$ and $|z-b_2|^{-1/2}$, respectively.
The end points of the lower arc $v_1,v_2$ are the Chebotarev points
and their density is proportional to $|z-v_1|^{1/2}$ and $|z-v_2|^{1/2}$;
see fig. \ref{Fig_ang3(13)_5000_180_rdbk}.
At the junction of the red, black and blue points
appeared two equilateral triangles with multicolored vertexes,
and the center of each has a Chebotarev point $v_1,v_2$;
see fig. \ref{Fig_ang3(13)_5000_180_full}.

\subsection{Nikishin system}\label{s2s2}
In the theory of Pad\'{e} approximants it is well known
the property called ``{\it Froissart doublets}'',
which was experimentally found (see \cite{BaGr81}, \cite{Gil78}, \cite{GiKr03})
and means that in the maximal\footnote{This notion has been introduced by
H. Stahl \cite{Sta85a}--\cite{Sta85c}. The existence of such a domain
$D=D(f)=\myo\CC\setminus{S}$ follows from the classical theorem of Stahl
\cite{Sta97b}. In this regard, the domain $D$ is usually
called {\it Stahl domain} and the symmetrical compact set $S$ -- {\it Stahl
compact} for the multi-valued analytic function $f$.}
domain of holomorphy $D=D(f)$ of the multi-valued function $f$
for some $n\in\Lambda$, $\Lambda\subset\NN$ is an infinite sequence,
are positioned pairs of zero-pole of a diagonal Pad\'{e} approximant
(that is, different from each other zeros of the Pad\'{e} polynomials $P_{n,0},P_{n,1}$,
see fig. \ref{Fig_pade10_2500_130}--\ref{Fig_pade103_5000_266_full}).
For every fixed $n\in\Lambda$, these points are different from each other,
the pole is not dependent of any singularity of the original function,
the zero and the pole are infinitely close to each other when $n\to\infty$,
that is they are asymptotically ``canceled out''.
In other words, the residue of the Pad\'{e} approximant at such a pole
converges to zero when $n\to\infty$.
Because such poles do not correspond to the singularities of the original function,
sometimes they are called ``spurious'' poles and zeros or ``defects''
of the diagonal Pad\'{e} approximant.
In a typical case these poles and zeros are dense on the Riemann surface
$\myo\CC$ when $n\to\infty$, $n\in\Lambda$
(see \cite{Sta96}, \cite{Sta98}, \cite{Sue02b}, \cite{Sue04}, \cite{Sue10}).
Thus they are sometimes called wondering poles and zeros or
floating poles and zeros. For an arbitrary
algebraic function $f$ the number of these pairs
depends mainly on the genus of the corresponding Riemann surface,
and also of the number of zeros of the functions $\Delta{f}(\zeta)=(f^{+}-f^{-})(\zeta)$,
$\zeta\in{S}$, on the Stahl compact set $S=S(f)$.
It is shown in \cite{Sue04}, that the appearance of the Froissart doublets
is due to points of an ``incorrect'' interpolation of the diagonal Pad\'{e} approximants
$[n/n]_f:=-P_{n,0}/P_{n,1}$ in the Stahl domain $D(f)=\myo\CC\setminus{S}$ with another branch
$\myt{f}$ of the original function $f$ when $n\in\Lambda$. Namely, the existence
of Froissart doublets does not allow uniform convergence of the Pad\'{e} approximants
in the Stahl domain (for details see below).

The second part of the numerical experiments are in the case,
when the sets of singularity points for the two functions
$f_1$ and $f_2$ {\it intersect} each other.
Specifically, we select a collection of three functions $[1,f,f^2]$, where
the function $f$ is of the type \eqref{4} and thus, the pair of functions $f,f^2$
forms a Nikishin system (see \cite{ApLy10}, \cite{Rak11}).
In this case another new property has been found. Namely,
the appearance of triple zeros ({\it Froissart triplets}, see below), i.e.
zeros of the Hermite--Pad\'{e} polynomials $Q_{n,0},Q_{n,1},Q_{n,2}$,
which are very close to each other, but still have different values.
These zeros are in the domain of holomorphy
of the functions $f,f^2$, do not correspond to either zeros, nor singularities
of these function, for each $n$ they are practically identical to each other
and with the transfer from $n$ to $n+1$ they shifted in the complex plane as one
unit.
It is appropriate to compare these triplets with the very well known
Froissart doublets for classical Pad\'{e} polynomials
(see \cite{Fro69}, also \cite{BaGr81}, \cite{Gil78}, \cite{GiTr87},
\cite{GiKr03}, \cite{GoGu13}, \cite{AdIb13}, \cite{Bel09}, \cite{Sta97b}, \cite{Tre13}).
Froissart doublets are sometimes called ``defects'' \cite[Chapter 2, \S\,2.2]{BaGr81},
and also spurious or wondering (floating) zeros and poles of the Pad\'{e} approximant
\cite{BaGr81}, \cite{Chu80} (see also \cite{Sue04}).

It is considered that the existence of such zeros and poles of
the Pad\'{e} approximant does not allow uniform convergence of the Pad\'{e} approximant.
Because such zeros and poles are infinitely close to each other asymptotically,
then when the limit is taken, they are practically ``canceled out''.
For some classes of hyperelliptic functions $f$,
which allow the representation $f=\myh{\sigma}$, where the support
$S_\sigma=\bigsqcup\limits_{j=1}^{2g+2}[e_{2j-1},e_{2j}]\subset\RR$ of the measure
$\sigma$ consists of finite number of non-intersecting intervals, it
was shown \cite{Sue00}, \cite{Sue06}, \cite{ApBuMaSu11} that the movement
of such poles is subject to certain regularity.
Namely, the corresponding divisor of the Nuttall $\Psi$ function
(see \cite{Sue00}, \cite{Sue06}, \cite{ApBuMaSu11})
moves along a Riemann surface $\RS_2$ with two sheets of genus $g\ge1$
and is subject to the general {\it Dubrovin system}:
\begin{equation}
\label{eqdub}
\dot z_k=-\frac{2w(\zz_k)}{\prod\limits_{j\neq k}(z_k-z_j)}
\int_{e_{2g+2}}^{\infty}\frac{\prod\limits_{j\neq k}(x-z_j)}{w(x)}\,dx,
\qquad k=1,\dots,g,
\end{equation}
where $\dot z_k=dz_k/dt$, $e_{2g+2}$ is the rightmost point of the support
$S_\sigma=\bigsqcup\limits_{j=1}^{2g+2}[e_{2j-1},e_{2j}]\subset\RR$ of the measure
$\sigma$ (it is supposed, that the support of the measure $\sigma$ has $g\geq1$ gaps,
the endpoints of the intervals of the support are numbered according to the
ascending values, and the path of integration in \eqref{eqdub} is part
of the real axis $[e_{2g+2},+\infty)$; for details about the notations
see \cite{Sue02b}). The points
$\zz_1(t),\dots,\zz_{g}(t)$, $t\in\RR_{+}$, are on the corresponding
hyperelliptic Riemann surface $\RS_2$ of genus $g$.
If $r$ is a real-valued rational function, which has poles
only outside $\myh{S}_\sigma$, then for ${f}=\myh\sigma+r$ we have: {\it
$[n/n]_{{f}}\to {f}$ in the spherical metric locally uniformly in
$\myo{\CC}\setminus\myh{S}_\sigma$} ($[n/n]_f:=-P_{n,1}/P_{n,0}$);
thus, each pole $r$ attracts exactly the number of poles of $[n/n]_{{f}}$,
as its multiplicity (see \cite{Gon03}).
In this case the movement of the poles and the points of interpolation
of the Pad\'{e} approximants $[n/n]_{{f}}$, which are in the
gaps between intervals, are subject to the system \eqref{eqdub}.

The fact, that such poles of the Pad\'{e} approximants form a pair
with the zeros that are near them and in this pair they are asymptotically close to each other
and when taking the limit are canceled out,
has led some authors to believe that their appearance
is random and is not connected with the nature of the original function $f$.
Such an approach to this property reflected on the respective terminology:
such zeros and poles of the approximant began to be called ``spurious''.
Their appearance, when using Pad\'{e} approximants, was considered especially negatively
in \cite{BaGr81}, \cite{Bel09}, \cite{Sta97}, therefore when
$n\in\Lambda$ the Pad\'{e} approximants $[n/n]_f$ themselves
sometimes were thought to be defective or entirely excluded
from the research \cite[Chapter 2, \S\,2.3]{BaGr81} or was proposed of using the
so-called ``purification'' of the spurious poles \cite{Sta96}.

In \cite{Dum08} Dumas researched the problem of the asymptotical
behavior of the sequence $\{[n/n]_d\}_{n\in\NN}$ for
elliptic\footnote{Speaking about elliptic functions,
we follow to the terminology of the monograph \cite[Chapter 10, par. 10.10]{Spr60},
where under elliptic function is understood a single-valued function
defined on an elliptical Riemann surface.}
functions of the special form
\begin{equation}
d(z)=\sqrt{(z-e_1)\dotsb(z-e_4)}-z^2+\(e_1+\dotsb+e_4\)z/2,
\label{1.1}
\end{equation}
where the points $e_1,\dots,e_4\in\mathbb C$ are pairwise different
and such root branch is selected,
that the main member of which, in a neighborhood of the point $z=\infty$,
is equal to $z^2$; thus $d\in\HH(\infty)$.
Particularly, Dumas has shown that in
{\it ``general position'' the set of poles of the Pad\'{e} approximants $[n/n]_d$
is dense in $\overline{\mathbb C}$}.

In \cite{Sue03} the result of Dumas has been extended to some
{\it classes\/} of elliptical functions.
It was found \cite[\S\,1, par .2]{Sue04} the following property:
for each $n$, from some subsequence $\Lambda$,
always exists a point $\mb=\mb_n\in{D}=\myo{\CC}\setminus{S}$,
for which
$[n/n]_d(\mb_n)=\myt{d}(\mb_n)$, where
$\myt{d}(z)=-\sqrt{(z-e_1)\dotsb(z-e_4)}-z^2+(e_1+\dotsb+e_4)z/2$ is
{\it another branch\/} of the elliptical function \eqref{1.1} in $D$.
Thus, for each $n$, from some subsequence $\Lambda$,
the diagonal Pad\'{e} approximants $[n/n]_d$ {\it interpolate\/} at
some point from the domain $D$ {\it another branch} $\myt{d}$ of the function $d$.
In the case of general position, such points of the ``incorrect'' interpolation
$\{\mb_n\}_{n\in\Lambda}$ are dense in $\myo{\CC}$.
In \cite{Sue04} were obtained much more general results in this direction.

In the work of E.~A.~Rakhmanov \cite{Rak12} was obtained an electrostatic interpretation
of the poles of the diagonal Pad\'{e} approximant $[n/n]_f$
for some algebraic functions $f$ {\it for each fixed} $n$.
In this framework, the spurious poles play a special role:
they should be considered a part of the ``external field''.

\subsection{}\label{s2s3}
\begin{figure}[!ht]
\centerline{
\includegraphics[width=15cm,height=15cm]{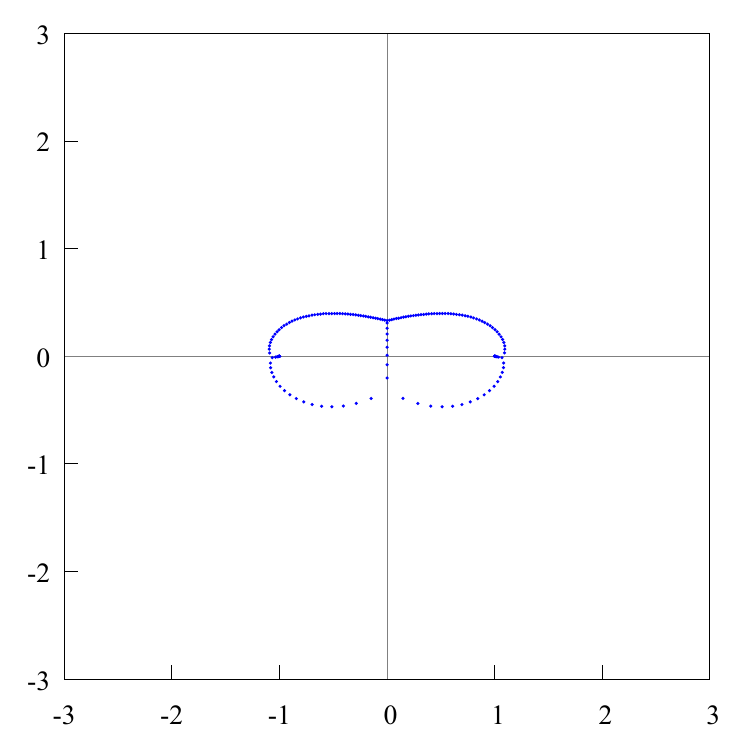}}
\vskip-6mm
\caption{The distribution of the zeros of the Hermite-Pad\'{e} polynomials $Q_{n,0}$ (blue points),
$n=121,\dots,130$, in the plane $\CC_z$
for a set of three functions $[1,f,f^2]$, where
$f(z)=(1-z^2)^{1/4}(1-i \sqrt{3}\cdot 1.6z)^{-1/2}$.
There are no Froissart zeros when $n=165$.
}
\label{Fig_nik_(1_6)_5000_165_bl}
\end{figure}

\begin{figure}[!ht]
\centerline{
\includegraphics[width=15cm,height=15cm]{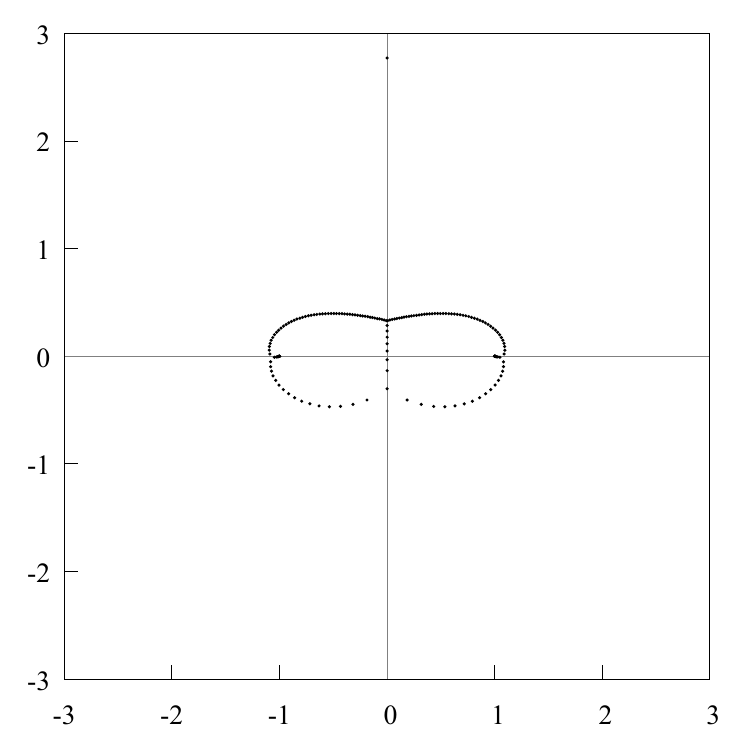}}
\vskip-6mm
\caption{The distribution of the zeros of the Hermite-Pad\'{e} polynomials $Q_{n,2}$ (black points),
$n=121,\dots,130$, in the plane $\CC_z$
for a set of three functions $[1,f,f^2]$, where
$f(z)=(1-z^2)^{1/4}(1-i \sqrt{3}\cdot 1.6z)^{-1/2}$.
There is a simple zero of the polynomial $Q_{n,2}$, $n=121,\dots,130$,
on the positive part of the imaginary axis at the point $z\approx a$, $a=i\sqrt{3}\cdot 1.6$,
corresponding to a simple pole of the function $f^2$ at the point $z=ia$.
There are no Froissart zeros when $n=165$.
}
\label{Fig_nik_(1_6)_5000_165_bk}
\end{figure}

\begin{figure}[!ht]
\centerline{
\includegraphics[width=15cm,height=15cm]{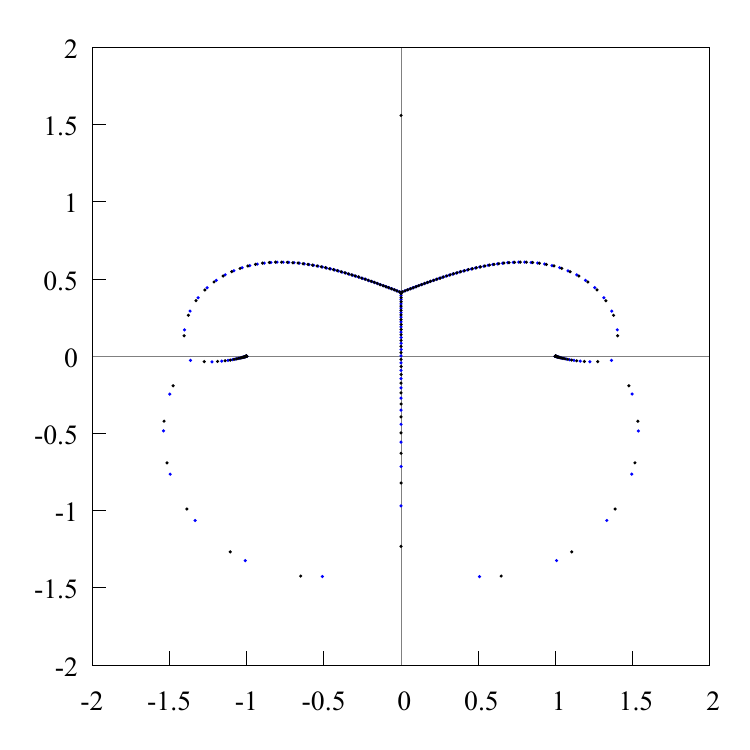}}
\vskip-6mm
\caption{The distribution of the zeros of the Hermite-Pad\'{e} polynomials $Q_{180,0},Q_{180,2}$
(blue and black points),
when $n=180$ for a set of three functions $[1,f,f^2]$, where
$f(z)=(1-z^2)^{1/4}(1-i\sqrt{3}\cdot 0.9z)^{-1/2}$.
There is a simple zero of the polynomial $Q_{180,2}$
on the positive part of the imaginary axis at the point $z\approx a$, $a=i\sqrt{3}\cdot 0.9$,
corresponding to a simple pole of the function $f^2$ at the point $z=ia$.
}
\label{Fig_nik_(_9)_5000_180_blbk}
\end{figure}

\begin{figure}[!ht]
\centerline{
\includegraphics[width=15cm,height=15cm]{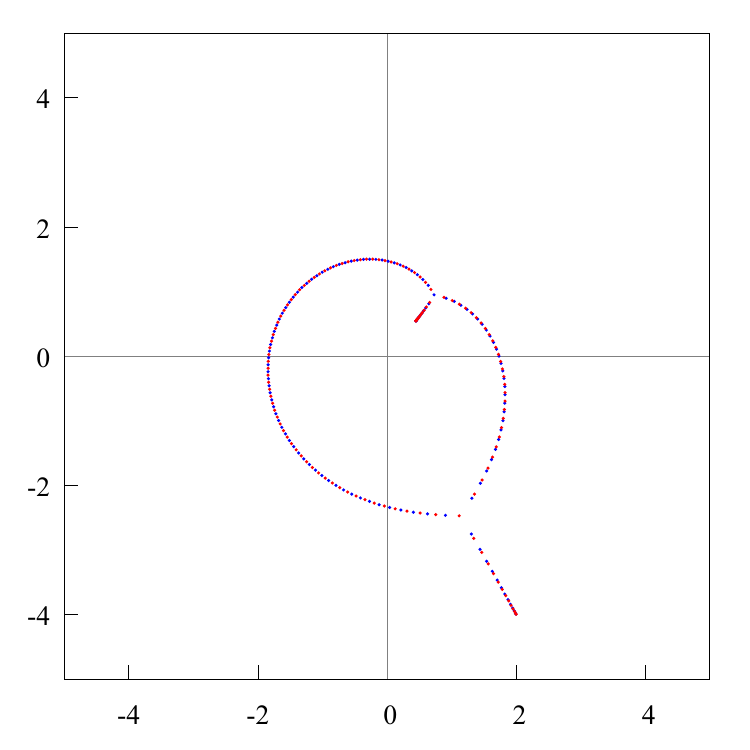}}
\vskip-6mm
\caption{The distribution of zeros and poles of the two-point Pad\'{e} approximant
of order $[120/120]$ of the function $f(z)=(z-a_1)^\alpha(z-a_2)^{-\alpha}$, $\alpha=1/4$,
$a_1=0.9-i\cdot1.1$ ,$a_2=0.1+i\cdot0.2$,
when two ``quite different'' branches are selected,
$f_0=\sqrt[4]{(z-a_1)/(z-a_2)}$, $f_\infty=-\sqrt[4]{(z-a_1)/(z-a_2)}$.
The zeros (blue points) and the poles (red points) create a Buslaev compact.
}
\label{Fig_bus205c_4000_120_full}
\end{figure}

In the present work, experimentally was obtained a new property
for the Nikishin system $f_1=f,f_2=f^2$, where
\begin{equation}
f(z)=(z^2-1)^{1/4}(z-a)^{-1/2},\quad\text{where}\quad
f(\infty)=1,\quad a\notin\RR,
\label{1.2}
\end{equation}
which is connected with the behavior of the zeros of the Hermite--Pad\'{e}
polynomials of first kind. Namely, there appear, for some $n\in\Lambda=\Lambda(f)$,
single spurious zeros of the Hermite--Pad\'{e} polynomials,
and also triple spurious zeros of the Hermite--Pad\'{e} polynomials.
Under ``spurious'' zeros of the Hermite--Pad\'{e} polynomials we understand
such zeros that, first, are different from each other (they cannot be canceled out),
second, do not correspond to either zeros, nor singularities of the original function,
and third, significantly change their location, when we transfer from $n$ to $n+1$
(sometimes they may disappear). By analogue with Pad\'{e} polynomials, it is
naturally to call such zeros {\it Froissart singlets} and {\it Froissart triplets}.
It is natural to assume the same as for Pad\'{e} approximants,
in the ``typical'' case (i.e. for the branch points $a$, which are in
``general position'', that is $a\notin i\RR$, see \eqref{1.2})
these spurious zeros (singlets and triplets)
of the Hermite--Pad\'{e} polynomials are dense only in ``their'' domain,
which is the difference from Pad\'{e} approximants.
From the numerical experiments it follows that
(see fig. \ref{Fig_nik_(1_6)_5000_165_bl}, \ref{Fig_nik_(1_6)_5000_165_bk},
\ref{Fig_nik_(_9)_5000_180_blbk}) the zeros of the polynomials $Q_{n,0}$ and $Q_{n,2}$ of
the functions of type \eqref{1.2} when $a=i\cdot 1.6$ have the same limiting distribution,
the corresponding boundary compact set $F_2$ separates the Riemann surface into
three domains, two of which have an internal boundary arc.
Thus, they have the same structure (see fig. \ref{Fig_bus205c_4000_120_full})
as in the theorem of Buslaev \cite{Bus13}, \cite{BuMaSu12}
for two-point Pad\'{e} approximants. Remark that,
the single zeros of the polynomial $Q_{n,2}$ at the point $z\approx a$
corresponds to a simple pole of the function $f^2$ at the point $z=a$.
The corresponding compact set $F_2$ has four Chebotarev's points (three of
which have zero density and one has infinite density) for the equilibrium measure,
which corresponds to the limiting distribution of the zeros of
the Hermite--Pad\'{e} polynomials $Q_{n,0},Q_{n,2}$.
It is clear, that the distribution of the zeros of the Hermite--Pad\'{e} polynomials
$Q_{n,0}$ and $Q_{n,2}$ must be equivalent, since with the mapping
$f\mapsto 1/f$ the type of the singularity of the original function
stays the same (see \eqref{1.2}).

\begin{figure}[!ht]
\centerline{
\includegraphics[width=15cm,height=15cm]{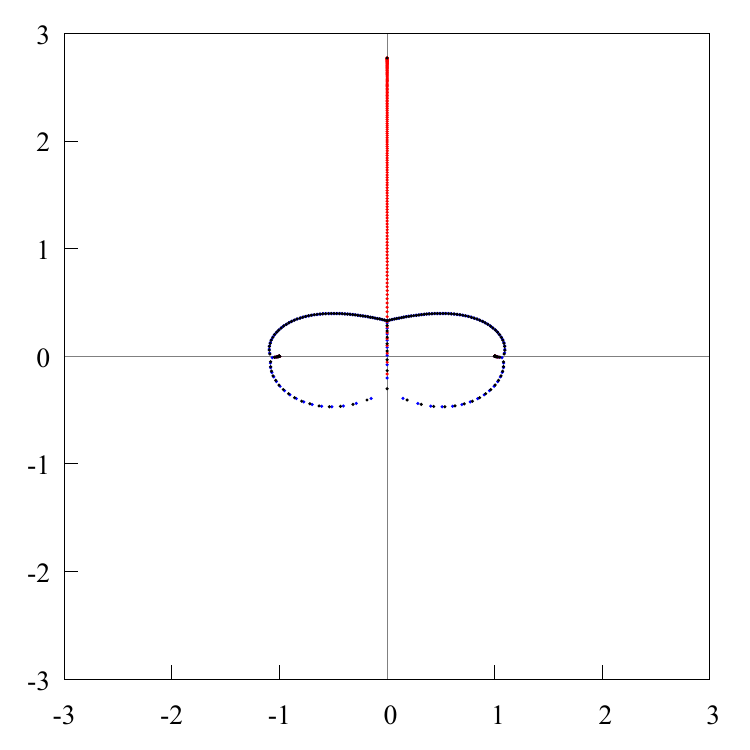}}
\vskip-6mm
\caption{The distribution of the zeros of the Hermite-Pad\'{e} polynomials $Q_{165,0}$ (blue points),
$Q_{165,1}$ (red points), $Q_{165,2}$ (black points) when $n=165$
for a set of three functions $[1,f,f^2]$, where
$f(z)=(1-z^2)^{1/4}(1-i \sqrt{3}\cdot 1.6z)^{-1/2}$.
The Riemann sphere is decomposed into 3 domains by the zeros of the Hermite-Pad\'{e} polynomials,
one of them contains the infinity point, while the other two are symmetrical with respect
to the imaginary axis. There are no Froissart zeros when $n=165$.
}
\label{Fig_nik_(1_6)_5000_165_full}
\end{figure}

\begin{figure}[!ht]
\centerline{
\includegraphics[width=15cm,height=15cm]{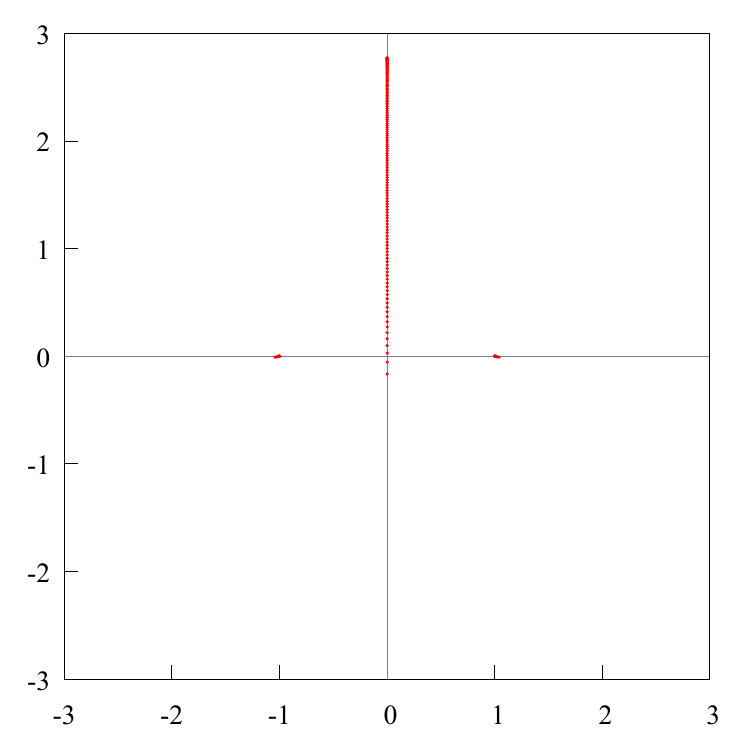}}
\vskip-6mm
\caption{The distribution of the zeros of the Hermite-Pad\'{e} polynomials $Q_{n,1}$ (red points),
$n=121,\dots,130$, in the plane $\CC_z$
for a set of three functions $[1,f,f^2]$, where
$f(z)=(1-z^2)^{1/4}(1-i \sqrt{3}\cdot 1.6z)^{-1/2}$.
There are no Froissart zeros when $n=165$.
}
\label{Fig_nik_(1_6)_5000_165_rd}
\end{figure}

\begin{figure}[!ht]
\centerline{
\includegraphics[width=15cm,height=15cm]{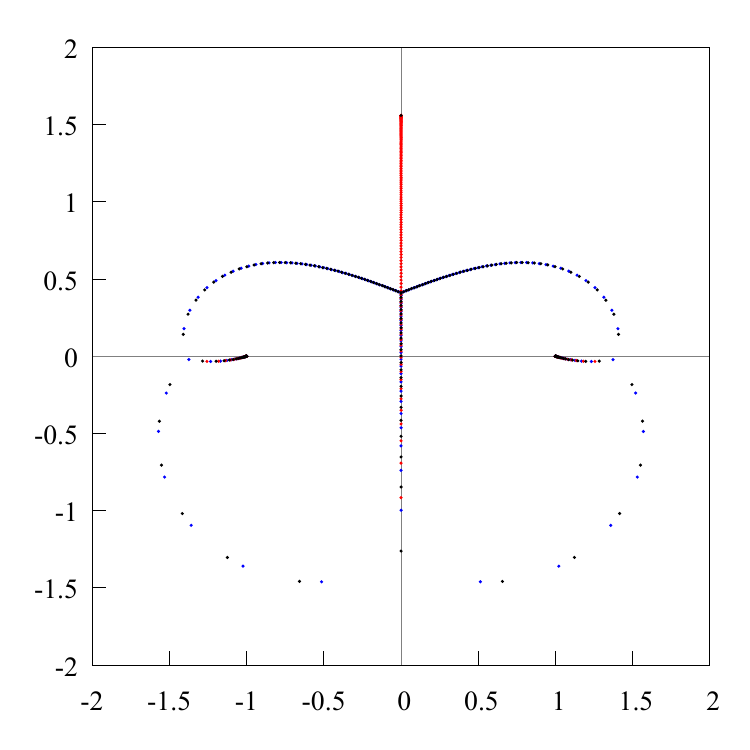}}
\vskip-6mm
\caption{The distribution of the zeros of the Hermite-Pad\'{e} polynomials $Q_{187,0}$ (blue points),
$Q_{187,1}$ (red points), $Q_{187,2}$ (black points) when $n=187$
for a set of three functions $[1,f,f^2]$, where
$f(z)=(1-z^2)^{1/4}(1-i \sqrt{3}\cdot 0.9z)^{-1/2}$.
The Riemann sphere is decomposed into 3 domains by the zeros of the Hermite-Pad\'{e} polynomials,
one of them contains the infinity point, while the other two are symmetrical with respect
to the imaginary axis.
There are no Froissart points when $n=187$.
}
\label{Fig_nik_(_9)_5000_187_full}
\end{figure}

\begin{figure}[!ht]
\centerline{
\includegraphics[width=15cm,height=15cm]{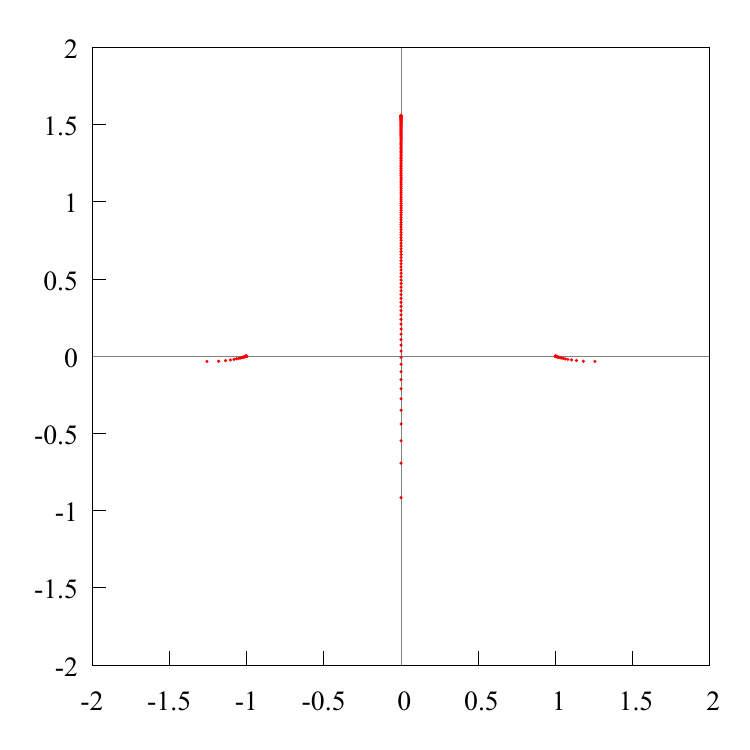}}
\vskip-6mm
\caption{The distribution of the zeros of the Hermite-Pad\'{e} polynomial $Q_{187,1}$ (red points),
when $n=187$ for a set of three functions $[1,f,f^2]$, where
$f(z)=(1-z^2)^{1/4}(1-i \sqrt{3}\cdot 0.9z)^{-1/2}$.
}
\label{Fig_nik_(_9)_5000_187_rd}
\end{figure}

\begin{figure}[!ht]
\centerline{
\includegraphics[width=15cm,height=15cm]{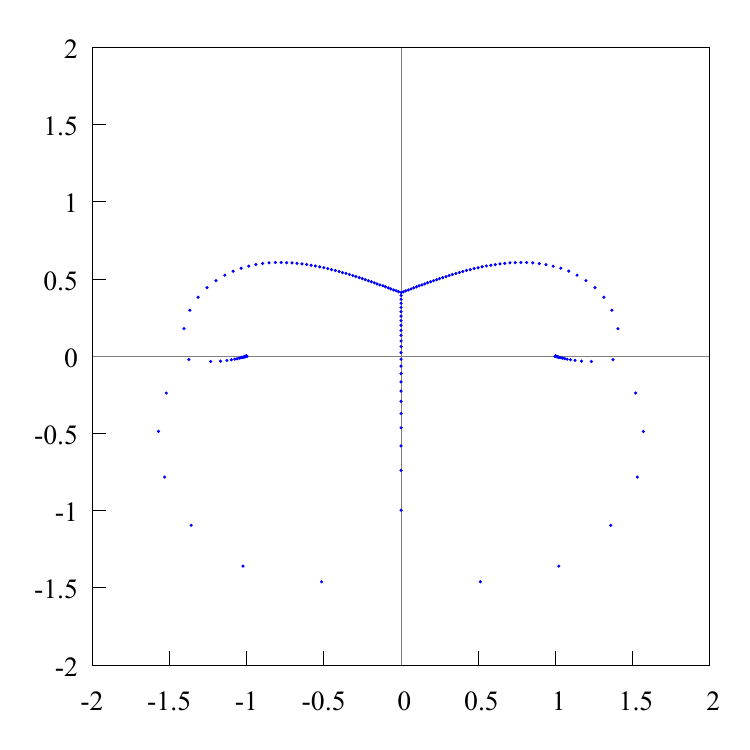}}
\vskip-6mm
\caption{The distribution of the zeros of the Hermite-Pad\'{e} polynomial $Q_{187,0}$ (blue points),
when $n=187$ for a set of three functions $[1,f,f^2]$, where
$f(z)=(1-z^2)^{1/4}(1-i \sqrt{3}\cdot 0.9z)^{-1/2}$.
}
\label{Fig_nik_(_9)_5000_187_bl}
\end{figure}

\begin{figure}[!ht]
\centerline{
\includegraphics[width=15cm,height=15cm]{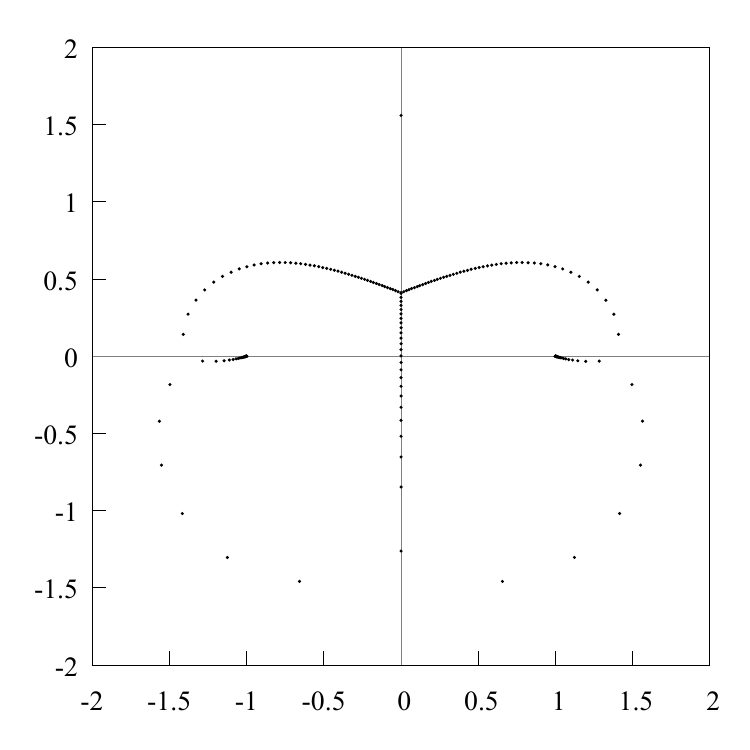}}
\vskip-6mm
\caption{The distribution of the zeros of the Hermite-Pad\'{e} polynomial $Q_{187,2}$ (black points),
when $n=187$ for a set of three functions $[1,f,f^2]$, where
$f(z)=(1-z^2)^{1/4}(1-i\sqrt{3}\cdot 0.9z)^{-1/2}$.
There is a simple zero of the polynomial $Q_{187,2}$
on the positive part of the imaginary axis at the point $z\approx a$, $a=i\sqrt{3}\cdot 0.9$,
corresponding to a simple pole of the function $f^2$ at the point $z=ia$.
}
\label{Fig_nik_(_9)_5000_187_bk}
\end{figure}

The limiting distribution of the zeros of the polynomial $Q_{n,1}$ must
be different. That is, it must correspond to the second compact set
$F_1$ (see fig. \ref{Fig_nik_(1_6)_5000_165_full},
\ref{Fig_nik_(1_6)_5000_165_rd} and fig.
\ref{Fig_nik_(_9)_5000_187_full}--\ref{Fig_nik_(_9)_5000_187_bk}). If
now we substitute $F:=F_1\cup F_2$, $F_1\cap F_2\neq\varnothing$, then
the complement $\myo\CC\setminus{F}$, as before, consists if three
domains, but now each of these domains has an internal boundary arc
(see fig. \ref{Fig_bus205c_4000_120_full}). Thus, $\myo\CC=F \sqcup
D_1\sqcup D_2\sqcup D_3$. Since the pair $f,f^2$ forms a Nikishin
system, then by the analogy with \cite{RaSu12}, \cite{RaSu13},
\cite{KovSu14}, \cite{Sue13b}, we can associate with the collection of
functions $[1,f,f^2]$ a Nuttall condenser. After then by the analogue
with \cite{RaSu13}, we can to describe in terms of the condenser the
limiting distribution of the zeros of the Hermite--Pad\'{e} polynomials
$Q_{n,0},Q_{n,1},Q_{n,2}$ and the corresponding Riemann surface with
the canonical (in Nuttall's sense) partition into three sheets. After
that, with the help of this Riemann surface, it might be possible to
find strong asymptotics for the Hermite--Pad\'{e} polynomials. The
question about how exactly to do this stays open; a heuristical method
for solving this problem for functions of the type
$f(z)=(z-a_1)^{\alpha_1}(z-a_2)^{\alpha_2}(z-a_3)^{\alpha_3}$,
$\alpha_1+\alpha_2+\alpha_3=0$, $2\alpha_j\in\CC\setminus\ZZ$, is
proposed in \cite{Sue13b}.

\begin{figure}[!ht]
\centerline{
\includegraphics[width=15cm,height=15cm]{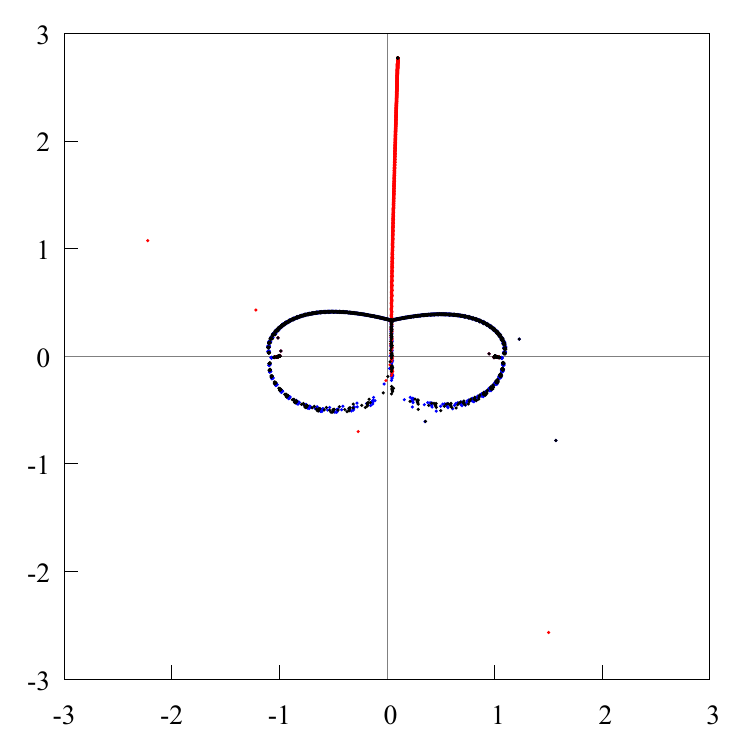}}
\vskip-6mm
\caption{The distribution of the zeros of the Hermite-Pad\'{e} polynomials $Q_{n,0}$ (blue points),
$Q_{n,1}$ (red points), $Q_{n,2}$ (black points), $n=121,\dots,130$,
for a set of three functions $[1,f,f^2]$, where $f$ is a ``perturbed'' function:
$f(z)=(1-z^2)^{1/4}(1-(0.1+i\sqrt{3}\cdot 1.6)z)^{-1/2}$.
The Riemann sphere is decomposed into 3 domains by the zeros of the Hermite-Pad\'{e} polynomials,
one of them contains the infinity point, while the other two are symmetrical with respect
to the imaginary axis. There is a pair of Froissart triplets
inside these two domans for some $n\in\{121,\dots,130\}$,
there is a pair of Froissart singlets
in the complementary domains for some $n\in\{121,\dots,130\}$,
there is one Froissart doublet on the negative part of the imaginary axis.
This follows from the analysis of the next figures \ref{Fig_nik_(1_6eps)_3000_121-130_rd},
\ref{Fig_nik_(1_6eps)_3000_121-130_bl} and \ref{Fig_nik_(1_6eps)_3000_121-130_bk}.
It is clearly seen, that the picture of the distribution of the zeros
is entirely the same compared to the unperturbed case.
}
\label{Fig_nik_(1_6eps)_3000_121-130_full}
\end{figure}

\begin{figure}[!ht]
\centerline{
\includegraphics[width=15cm,height=15cm]{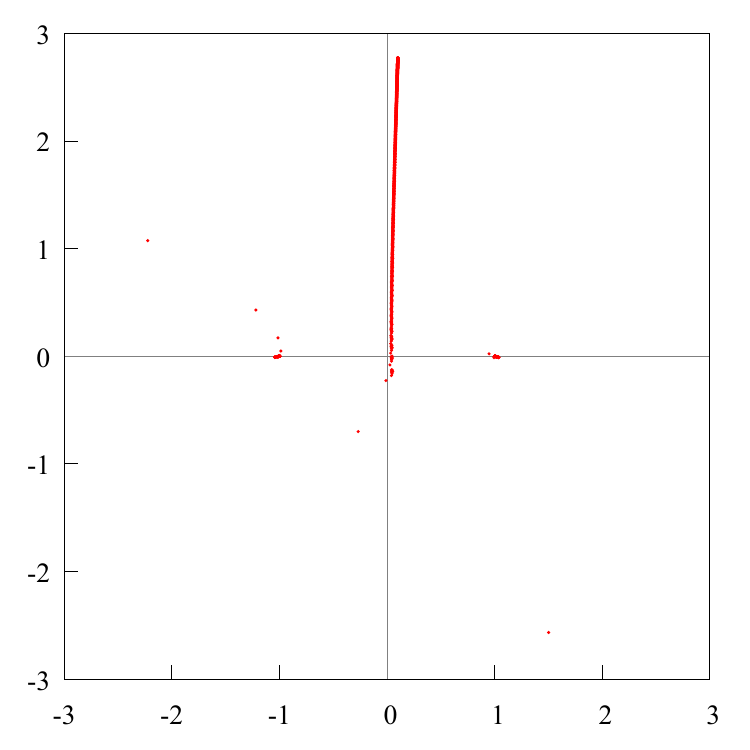}}
\vskip-6mm
\caption{The distribution of the zeros of the Hermite-Pad\'{e} polynomials $Q_{n,1}$ (red points),
$n=121,\dots,130$, for a set of three functions $[1,f,f^2]$, where
$f$ is a ``perturbed'' function:
$f(z)=(1-z^2)^{1/4}(1-(0.1+i\sqrt{3}\cdot 1.6)z)^{-1/2}$.
}
\label{Fig_nik_(1_6eps)_3000_121-130_rd}
\end{figure}

\begin{figure}[!ht]
\centerline{
\includegraphics[width=15cm,height=15cm]{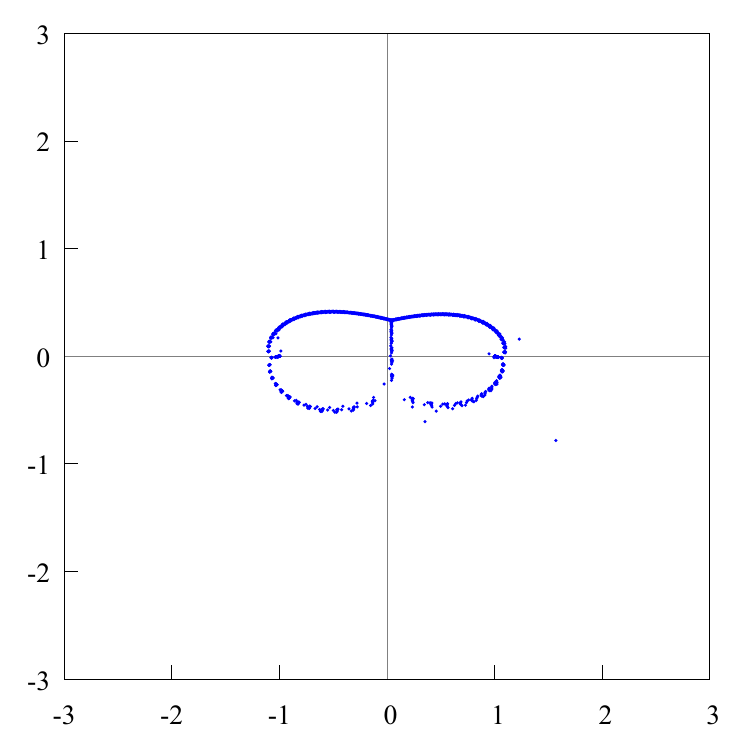}}
\vskip-6mm
\caption{The distribution of the zeros of the Hermite-Pad\'{e} polynomials $Q_{n,0}$ (blue points),
$n=121,\dots,130$, for a set of three functions $[1,f,f^2]$, where
$f$ is a ``perturbed'' function:
$f(z)=(1-z^2)^{1/4}(1-(0.1+i\sqrt{3}\cdot1.6)z)^{-1/2}$.
}
\label{Fig_nik_(1_6eps)_3000_121-130_bl}
\end{figure}

\begin{figure}[!ht]
\centerline{
\includegraphics[width=15cm,height=15cm]{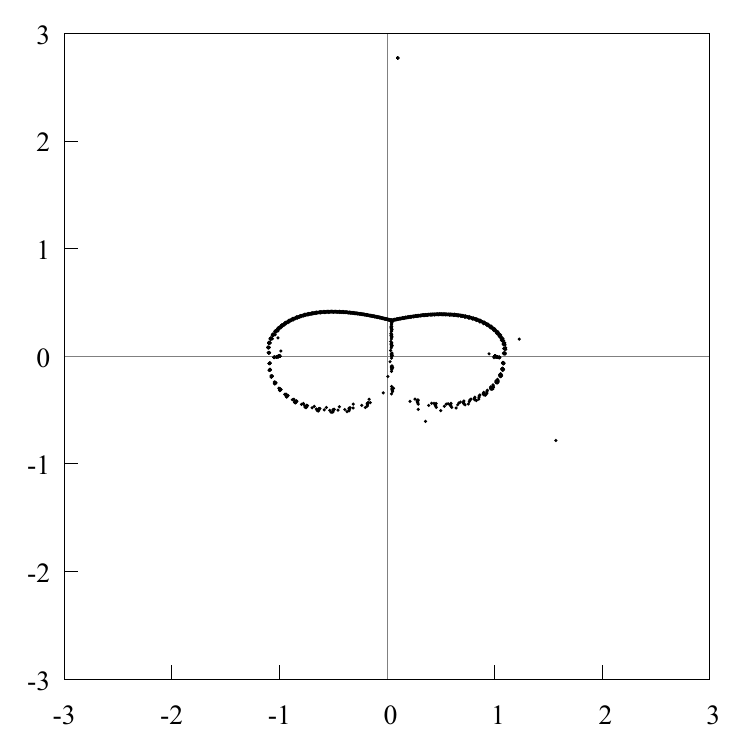}}
\vskip-6mm
\caption{The distribution of the zeros of the Hermite-Pad\'{e} polynomials $Q_{n,2}$ (black points),
$n=121,\dots,130$, for a set of three functions $[1,f,f^2]$, where
$f$ is a ``perturbed'' function:
$f(z)=(1-z^2)^{1/4}(1-(0.1+i\sqrt{3}\cdot1.6)z)^{-1/2}$.
There are a simple zeros of the polynomials $Q_{n,2}$, $n=121,\dots,130$,
on the positive part of the imaginary axis at the point $z\approx a$, $a=0.1+i\sqrt{3}\cdot1.6$,
corresponding to a simple pole of the function $f^2$ at the point $z=ia$.
}
\label{Fig_nik_(1_6eps)_3000_121-130_bk}
\end{figure}

Figures \ref{Fig_nik_(1_6eps)_3000_121-130_full}--\ref{Fig_nik_(1_6eps)_3000_121-130_bk}
show the distribution of the zeros of the Hermite--Pad\'{e} polynomials
$Q_{n,0}$ (blue points), $Q_{n,1}$ (red points), $Q_{n,2}$ (black points),
$n=121,\dots,130$, for the collection of three functions $[1,f,f^2]$, where $f$ is
``perturbed'' with respect to the function \eqref{1.2}:
\begin{equation}
\label{nik_eps}
f(z)=(z^2-1)^{1/4}(z-(0.1+i\sqrt{3}\cdot1.6))^{-1/2}.
\end{equation}
It is clear, that by comparison with the unperturbed case \eqref{1.2},
the picture of the distribution of the zeros is entirely the same.

\section{Concluding remarks}\label{s3}

Thus, the main empirical results of the paper are as follows.

\subsection{}\label{s3s1}
In the case when the pair of functions $f_1,f_2$ forms
an Angelesco system, where the functions $f_1,f_2$ have the type \eqref{3},
there has been found numerically the property of the {\it mutual pushing}
of the supports of the measures, which are equilibrium for
the extremal compact sets.

\subsection{}\label{s3s2}
In the case when the pair of functions $f_1=f,f_2=f^2$ forms a (generalized)
Nikishin
system, where the function $f$ has the type \eqref{4}, there have been
found numerically the properties of {\it Froissart singlets and triplets}
presence.

Thus, in the present paper there have been found numerically new properties
related to the behavior of the zeroes of the Hermite--Pad\'{e} polynomials
of the first kind.
These numerical phenomena should be researched and strictly justified in
future.




\clearpage
\markboth{\bf Pade approximant}{\bf Pade approximant}

\newpage

\newpage

\newpage

\newpage

\newpage

\newpage

\newpage

\newpage

\newpage

\newpage

\newpage




\clearpage
\markboth{\bf Case One, Angelesco system}{\bf Case One,  Angelesco system}

\newpage

\newpage
\begin{figure}[!ht]
\centerline{
\includegraphics[width=15cm,height=15cm]{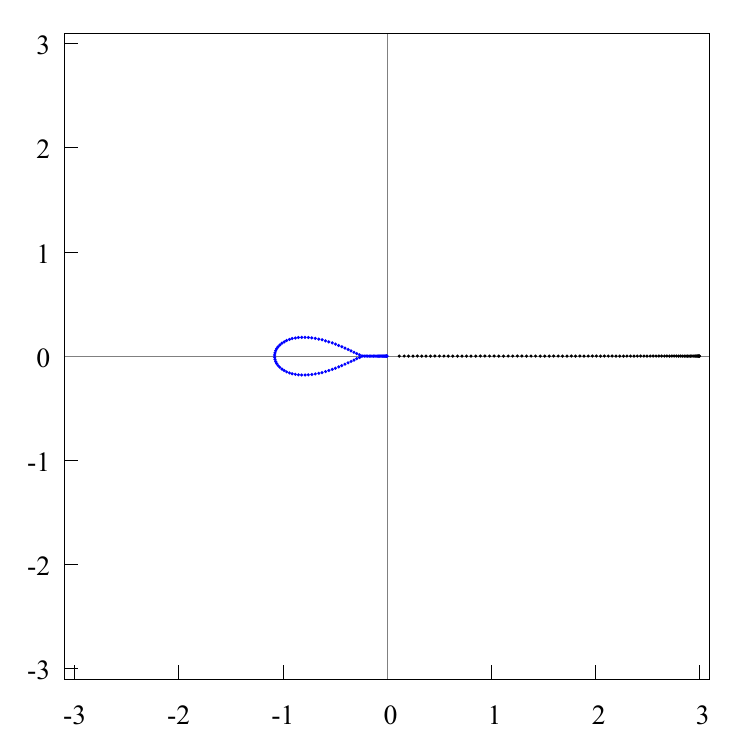}}
\vskip-6mm
\caption{The distribution of the zeros of the Hermite-Pad\'{e} polynomial $Q_{90,0}$ (blue points)
and $Q_{90,2}$ (black points) for the set of three functions $[1,f_1,f_2]$, where
$f_1(z)=\sqrt{(z+1)/z}$, $f_2=\sqrt{(z-3)/z}$.
The zeros of the polynomial $Q_{90,0}$ (blue points, comp fig.
\ref{Fig_ang3(1)_5000_120_full}) create a third extremal compact $F_0$,
which ``separates'' the compacts $F_1$ and $F_2$.
}
\label{Fig_mar_1_3000_90_blbk}
\end{figure}




\clearpage
\markboth{\bf Angelesco system and Nikishin system}{\bf Angelesco system and Nikishin system}

\newpage

\newpage

\newpage

\newpage

\newpage

\newpage

\newpage

\newpage

\newpage

\newpage

\newpage

\newpage

\newpage

\newpage

\newpage

\newpage

\newpage

\newpage




\clearpage
\markboth{\bf Case One, Angelesco system}{\bf Case One, Angelesco system}

\newpage

\newpage

\newpage

\newpage

\newpage

\newpage

\newpage

\newpage

\newpage

\newpage

\newpage

\newpage

\newpage

\newpage

\newpage

\newpage




\clearpage
\markboth{\bf Case Two, Nikishin system}{\bf Case Two, Nikishin system}

\newpage
\begin{figure}[!ht]
\centerline{
\includegraphics[width=15cm,height=15cm]{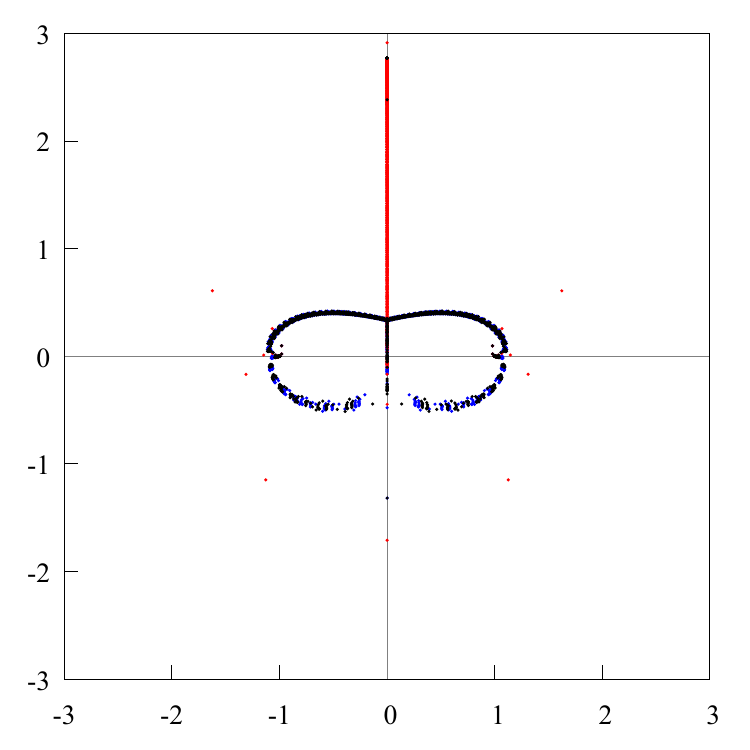}}
\vskip-6mm
\caption{The distribution of the zeros of the Hermite-Pad\'{e} polynomials $Q_{n,0}$ (blue points),
$Q_{n,1}$ (red points), $Q_{n,2}$ (black points), $n=61,\dots,80$,
for a set of three functions $[1,f,f^2]$, where
$f(z)=(1-z^2)^{1/4}(1-i \sqrt{3}\cdot 1.6z)^{-1/2}$.
The Riemann sphere is decomposed into 3 domains by the zeros of the Hermite-Pad\'{e} polynomials,
one of them contains the infinity point, while the other two are symmetrical with respect
to the imaginary axis. There is a pair of Froissart triplets
inside these two domans for some $n\in\{61,\dots,80\}$,
there is a pair of Froissart singlets
in the complementary domains for some $n\in\{61,\dots,80\}$,
there is one Froissart doublet on the negative part of the imaginary axis.
This follows from the analysis of the next figures \ref{Fig_nik_(1_6)_2000_61-80_rd},
\ref{Fig_nik_(1_6)_2000_61-80_bl}, and \ref{Fig_nik_(1_6)_2000_61-80_bk}.
}
\label{Fig_nik_(1_6)_2000_61-80_full}
\end{figure}

\newpage
\begin{figure}[!ht]
\centerline{
\includegraphics[width=15cm,height=15cm]{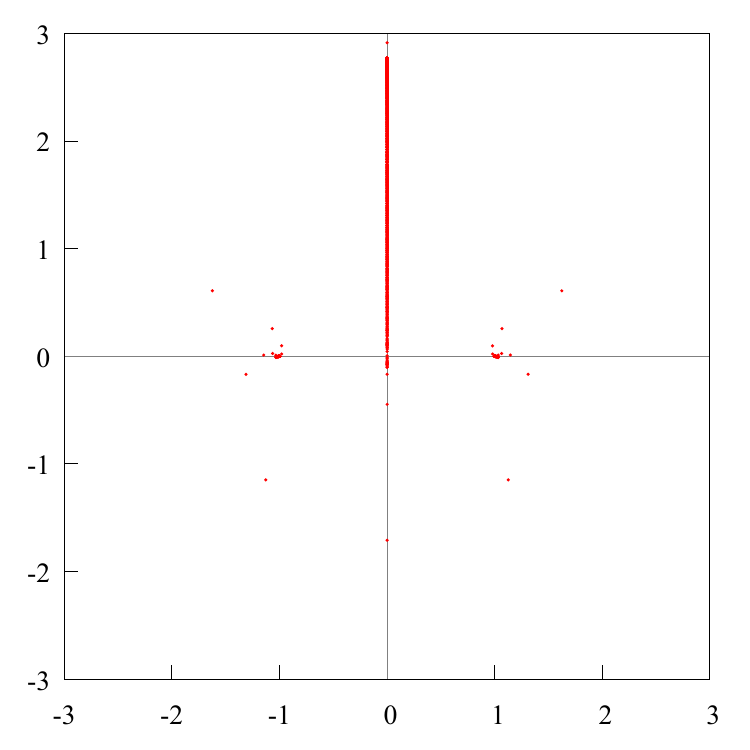}}
\vskip-6mm
\caption{The distribution of the zeros of the Hermite-Pad\'{e} polynomials $Q_{n,1}$ (red points),
$n=61,\dots,80$, for a set of three functions $[1,f,f^2]$, where
$f(z)=(1-z^2)^{1/4}(1-i \sqrt{3}\cdot 1.6z)^{-1/2}$.
}
\label{Fig_nik_(1_6)_2000_61-80_rd}
\end{figure}

\newpage
\begin{figure}[!ht]
\centerline{
\includegraphics[width=15cm,height=15cm]{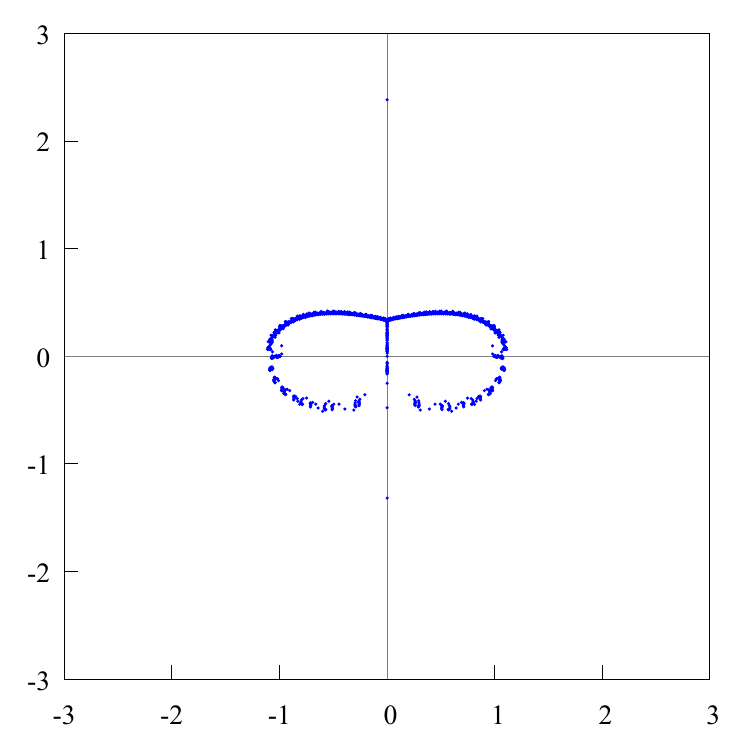}}
\vskip-6mm
\caption{The distribution of the zeros of the Hermite-Pad\'{e} polynomials $Q_{n,0}$ (blue points),
$n=61,\dots,80$, for a set of three functions $[1,f,f^2]$, where
$f(z)=(1-z^2)^{1/4}(1-i \sqrt{3}\cdot 1.6z)^{-1/2}$.
}
\label{Fig_nik_(1_6)_2000_61-80_bl}
\end{figure}

\newpage
\begin{figure}[!ht]
\centerline{
\includegraphics[width=15cm,height=15cm]{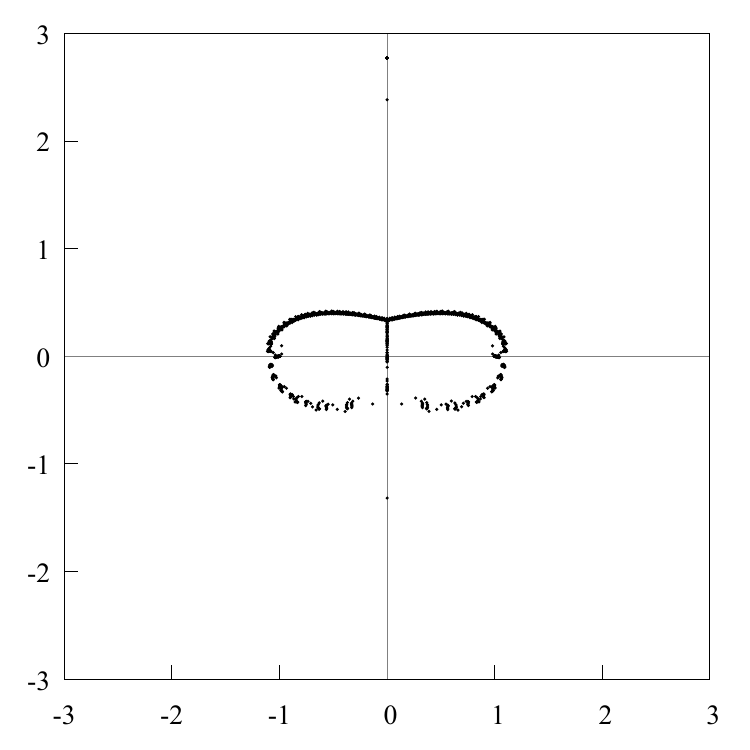}}
\vskip-6mm
\caption{The distribution of the zeros of the Hermite-Pad\'{e} polynomials $Q_{n,2}$ (black points),
$n=61,\dots,80$, for a set of three functions $[1,f,f^2]$, where
$f(z)=(1-z^2)^{1/4}(1-i \sqrt{3}\cdot 1.6z)^{-1/2}$. 
There is a simple zero of the polynomial $Q_{n,2}$, $n=61,\dots,80$,
on the positive part of the imaginary axis at the point $z\approx a$, $a=i\sqrt{3}\cdot 1.6$,
corresponding to a simple pole of the function $f^2$ at the point $z=ia$.
}
\label{Fig_nik_(1_6)_2000_61-80_bk}
\end{figure}

\newpage
\begin{figure}[!ht]
\centerline{
\includegraphics[width=15cm,height=15cm]{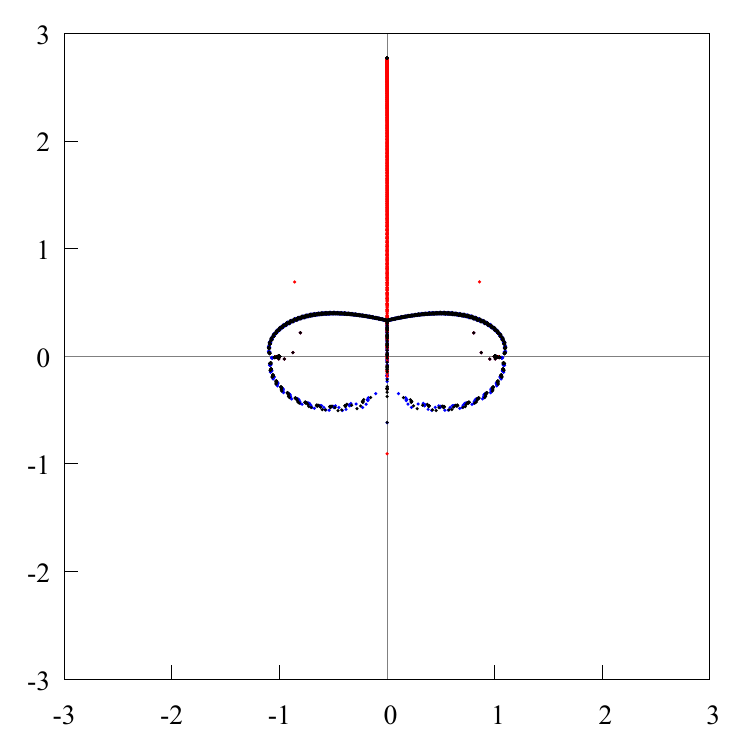}}
\vskip-6mm
\caption{The distribution of the zeros of the Hermite-Pad\'{e} polynomials $Q_{n,0}$ (blue points),
$Q_{n,1}$ (red points), $Q_{n,2}$ (black points), $n=121,\dots,130$,
for a set of three functions $[1,f,f^2]$, where
$f(z)=(1-z^2)^{1/4}(1-i \sqrt{3}\cdot 1.6z)^{-1/2}$.
The Riemann sphere is decomposed into 3 domains by the zeros of the Hermite-Pad\'{e} polynomials,
one of them contains the infinity point, while the other two are symmetrical with respect
to the imaginary axis. There is a pair of Froissart triplets
inside these two domans for some $n\in\{121,\dots,130\}$,
there is a pair of Froissart singlets
in the complementary domains for some $n\in\{121,\dots,130\}$,
there is a Froissart doublet on the negative part of the imaginary axis.
This follows from the analysis of the next figures \ref{Fig_nik_(1_6)_3000_121-130_rd},
\ref{Fig_nik_(1_6)_3000_121-130_bl} and \ref{Fig_nik_(1_6)_3000_121-130_bk}.
}
\label{Fig_nik_(1_6)_3000_121-130_full}
\end{figure}

\newpage
\begin{figure}[!ht]
\centerline{
\includegraphics[width=15cm,height=15cm]{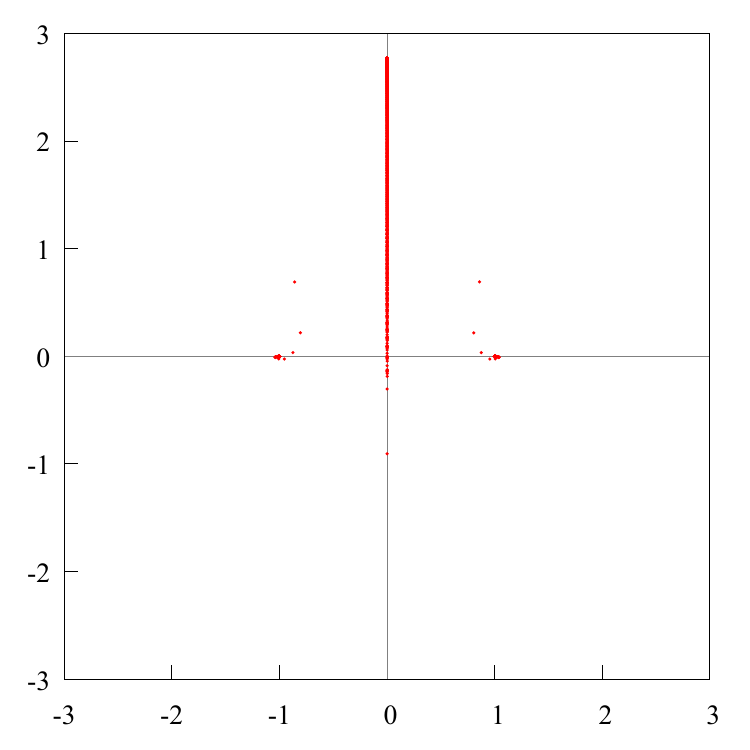}}
\vskip-6mm
\caption{The distribution of the zeros of the Hermite-Pad\'{e} polynomials $Q_{n,1}$ (red points),
$n=121,\dots,130$, for a set of three functions $[1,f,f^2]$, where
$f(z)=(1-z^2)^{1/4}(1-i \sqrt{3}\cdot 1.6z)^{-1/2}$.
}
\label{Fig_nik_(1_6)_3000_121-130_rd}
\end{figure}

\newpage
\begin{figure}[!ht]
\centerline{
\includegraphics[width=15cm,height=15cm]{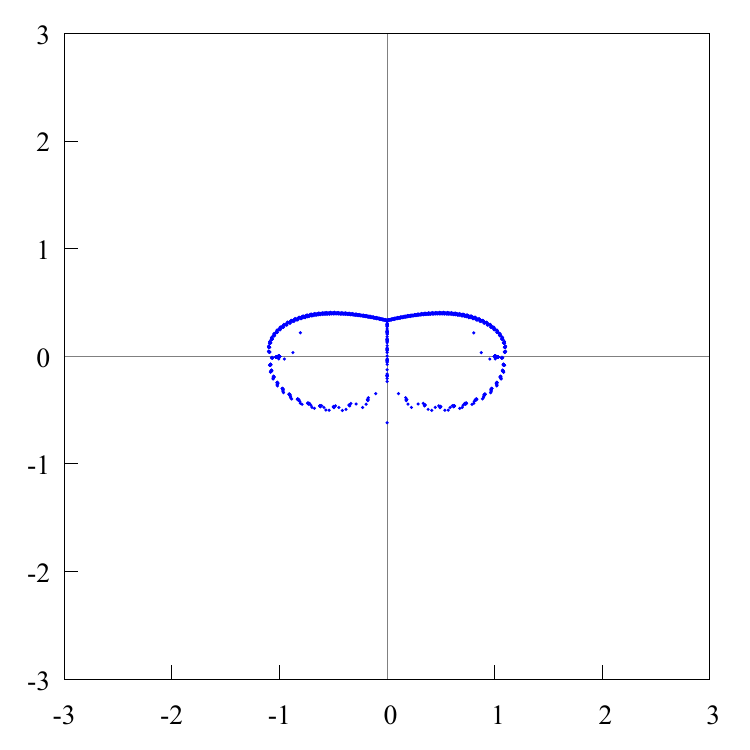}}
\vskip-6mm
\caption{The distribution of the zeros of the Hermite-Pad\'{e} polynomials $Q_{n,0}$ (blue points),
$n=121,\dots,130$, for a set of three functions $[1,f,f^2]$, where
$f(z)=(1-z^2)^{1/4}(1-i \sqrt{3}\cdot 1.6z)^{-1/2}$.
}
\label{Fig_nik_(1_6)_3000_121-130_bl}
\end{figure}

\newpage
\begin{figure}[!ht]
\centerline{
\includegraphics[width=15cm,height=15cm]{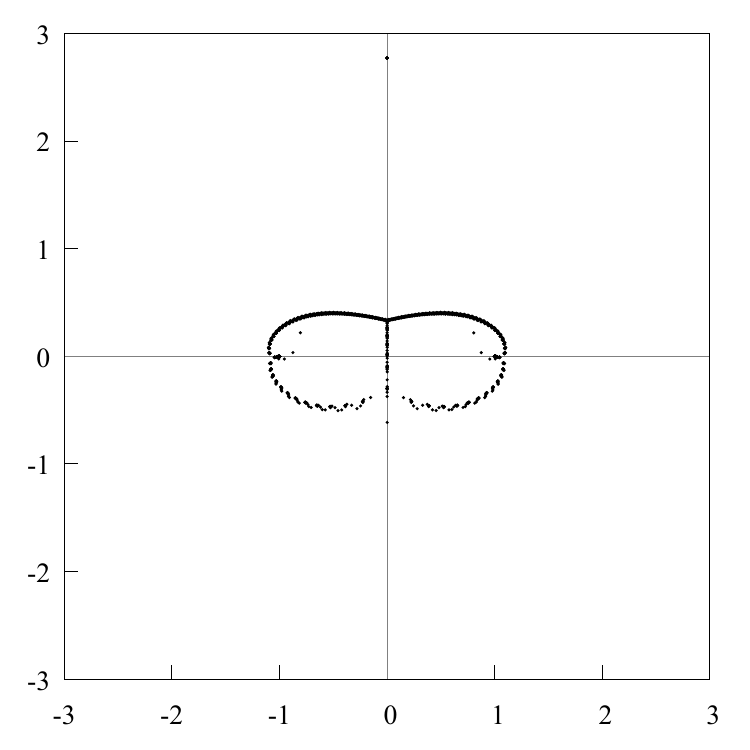}}
\vskip-6mm
\caption{The distribution of the zeros of the Hermite-Pad\'{e} polynomials $Q_{n,2}$ (black points),
$n=121,\dots,130$, for a set of three functions $[1,f,f^2]$, where
$f(z)=(1-z^2)^{1/4}(1-i \sqrt{3}\cdot 1.6z)^{-1/2}$.
There is a simple zero of the polynomial $Q_{n,2}$, $n=121,\dots,130$,
on the positive part of the imaginary axis at the point $z\approx a$, $a=i\sqrt{3}\cdot 1.6$,
corresponding to a simple pole of the function $f^2$ at the point $z=ia$.
}
\label{Fig_nik_(1_6)_3000_121-130_bk}
\end{figure}

\newpage

\newpage

\newpage

\newpage

\newpage
\begin{figure}[!ht]
\centerline{
\includegraphics[width=15cm,height=15cm]{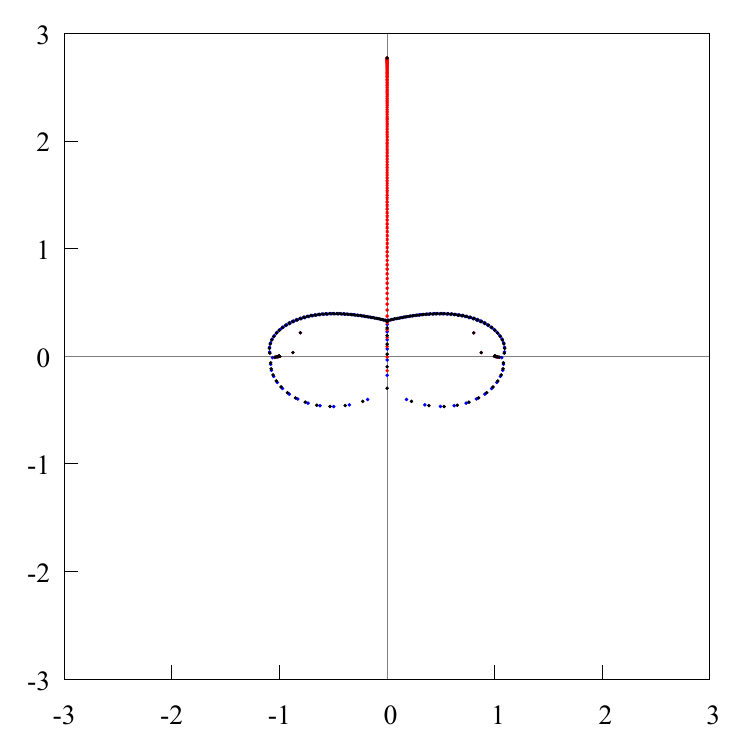}}
\vskip-6mm
\caption{The distribution of the zeros of the Hermite-Pad\'{e} polynomials $Q_{n,0}$ (blue points),
$Q_{n,1}$ (red points), $Q_{n,2}$ (black points), $n=127,128$,
for a set of three functions $[1,f,f^2]$, where
$f(z)=(1-z^2)^{1/4}(1-i \sqrt{3}\cdot 1.6z)^{-1/2}$.
The Riemann sphere is decomposed into 3 domains by the zeros of the Hermite-Pad\'{e} polynomials,
one of them contains the infinity point, while the other two are symmetrical with respect
to the imaginary axis. There is a pair of Froissart triplets
inside these two domans for some $n\in\{127,128\}$,
there is a pair of Froissart singlets
in the complementary domains for some $n\in\{127,128\}$,
there is one Froissart doublet on the negative part of the imaginary axis.
This follows from the analysis of the next figures \ref{Fig_nik_(1_6)_3000_127-128_rd},
\ref{Fig_nik_(1_6)_3000_127-128_bl} and \ref{Fig_nik_(1_6)_3000_127-128_bk}.
}
\label{Fig_nik_(1_6)_3000_127-128_full}
\end{figure}

\newpage
\begin{figure}[!ht]
\centerline{
\includegraphics[width=15cm,height=15cm]{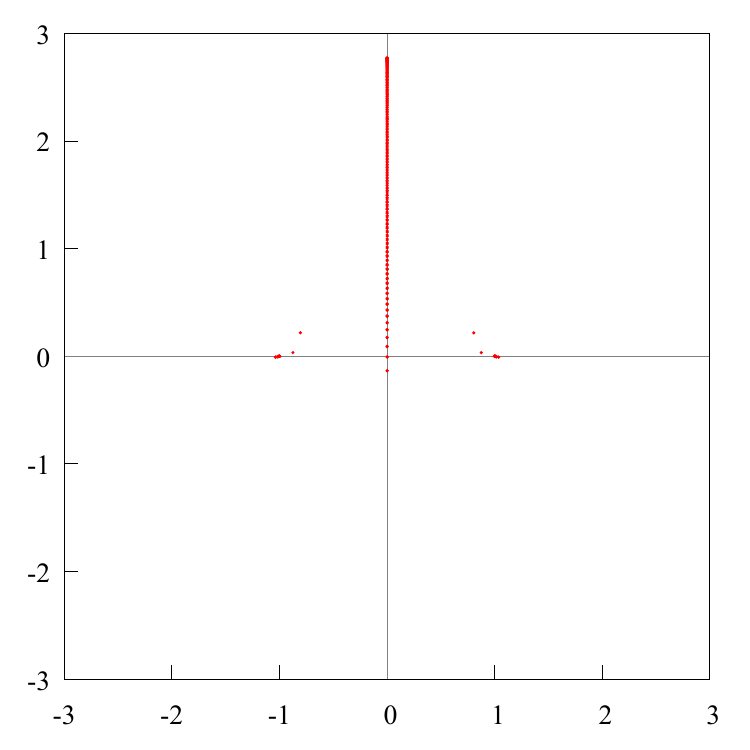}}
\vskip-6mm
\caption{The distribution of the zeros of the Hermite-Pad\'{e} polynomials $Q_{n,1}$ (red points),
$n=127,128$, for a set of three functions $[1,f,f^2]$, where
$f(z)=(1-z^2)^{1/4}(1-i \sqrt{3}\cdot 1.6z)^{-1/2}$.
}
\label{Fig_nik_(1_6)_3000_127-128_rd}
\end{figure}

\newpage
\begin{figure}[!ht]
\centerline{
\includegraphics[width=15cm,height=15cm]{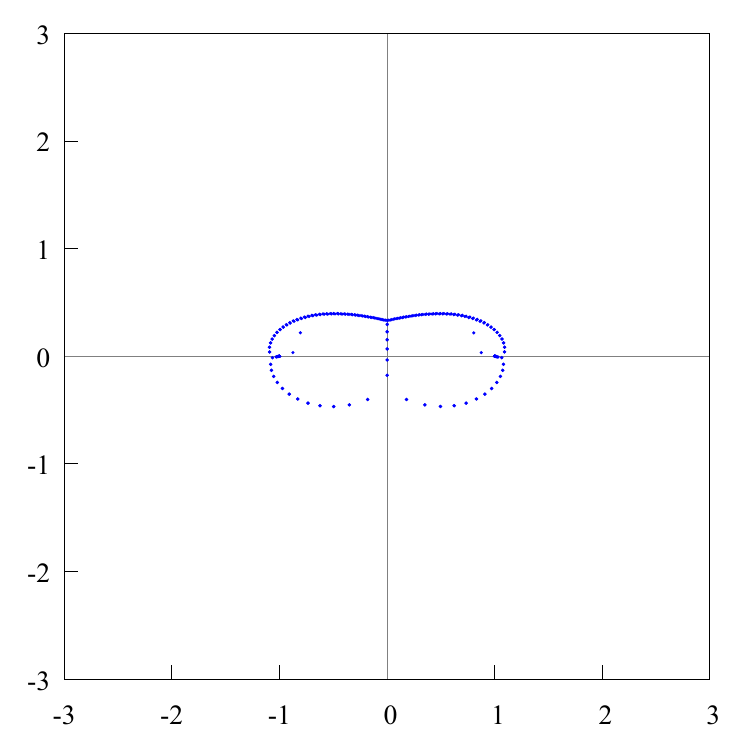}}
\vskip-6mm
\caption{The distribution of the zeros of the Hermite-Pad\'{e} polynomials $Q_{n,0}$ (blue points),
$n=127,128$, for a set of three functions $[1,f,f^2]$, where
$f(z)=(1-z^2)^{1/4}(1-i \sqrt{3}\cdot 1.6z)^{-1/2}$.
}
\label{Fig_nik_(1_6)_3000_127-128_bl}
\end{figure}

\newpage
\begin{figure}[!ht]
\centerline{
\includegraphics[width=15cm,height=15cm]{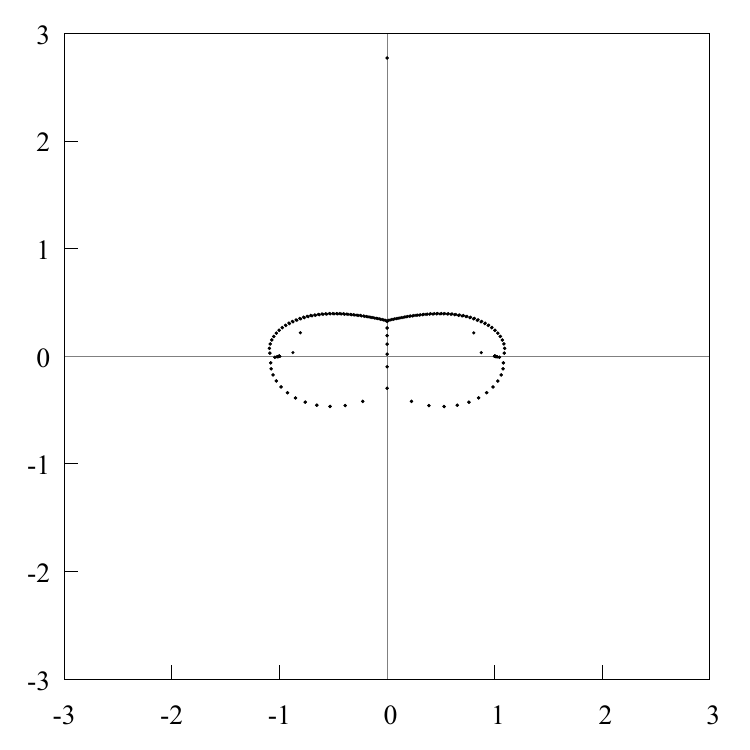}}
\vskip-6mm
\caption{The distribution of the zeros of the Hermite-Pad\'{e} polynomials $Q_{n,2}$ (black points),
$n=127,128$, for a set of three functions $[1,f,f^2]$, where
$f(z)=(1-z^2)^{1/4}(1-i \sqrt{3}\cdot 1.6z)^{-1/2}$.
There is a simple zero of the polynomial $Q_{n,2}$, $n=127,128$,
on the positive part of the imaginary axis at the point $z\approx a$, $a=i\sqrt{3}\cdot 1.6$,
corresponding to a simple pole of the function $f^2$ at the point $z=ia$.
}
\label{Fig_nik_(1_6)_3000_127-128_bk}
\end{figure}




\clearpage
\markboth{\bf $Q_{n,0}\cdot 1+Q_{n,1}f+Q_{n,2}f^2$}{\bf $Q_{n,0}\cdot 1+Q_{n,1}f+Q_{n,2}f^2$}

\newpage

\newpage

\newpage

\newpage

\newpage
\begin{figure}[!ht]
\centerline{
\includegraphics[width=15cm,height=15cm]{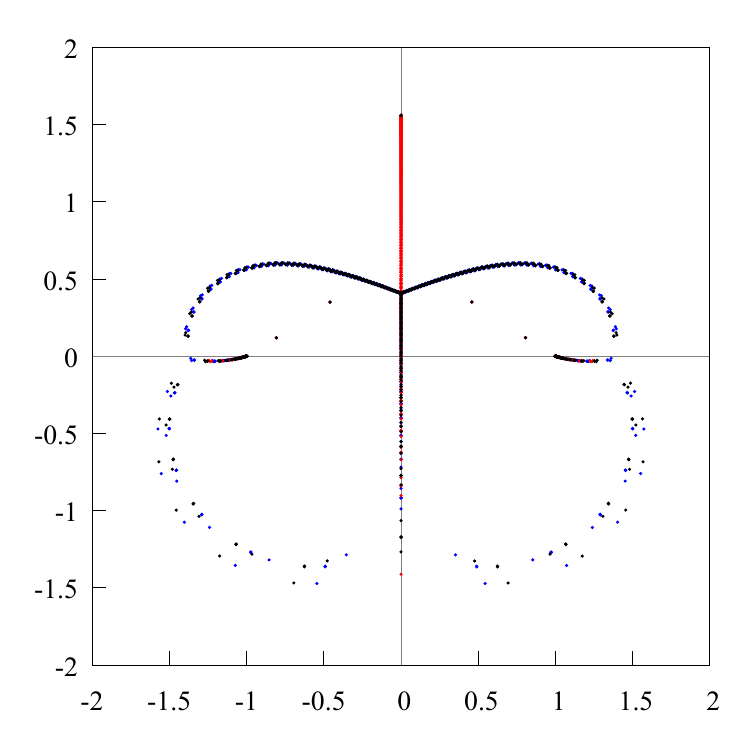}}
\vskip-6mm
\caption{
The distribution of the zeros of the Hermite-Pad\'{e} polynomials $Q_{n,0}$ (blue points),
$Q_{n,1}$ (red points), $Q_{n,2}$ (black points) when $n=166,\dots,170$
for a set of three functions $[1,f,f^2]$, where
$f(z)=(1-z^2)^{1/4}(1-i \sqrt{3}\cdot 0.9z)^{-1/2}$.
The Riemann sphere is decomposed into 3 domains by the zeros of the Hermite-Pad\'{e} polynomials,
one of them contains the infinity point, while the other two are symmetrical with respect
to the imaginary axis.
There is a pair of Froissart triplets
inside these two domans for two $n$ of $n\in\{166,\dots,170\}$.
This follows from the analysis of the next figures \ref{Fig_nik_(_9)_4000_166-170_rd},
\ref{Fig_nik_(_9)_4000_166-170_bl}, \ref{Fig_nik_(_9)_4000_166-170_bk}.
There are no Froissart points in the complementary domains for $n=166,\dots,170$.
}
\label{Fig_nik_(_9)_4000_166-170_full}
\end{figure}

\newpage
\begin{figure}[!ht]
\centerline{
\includegraphics[width=15cm,height=15cm]{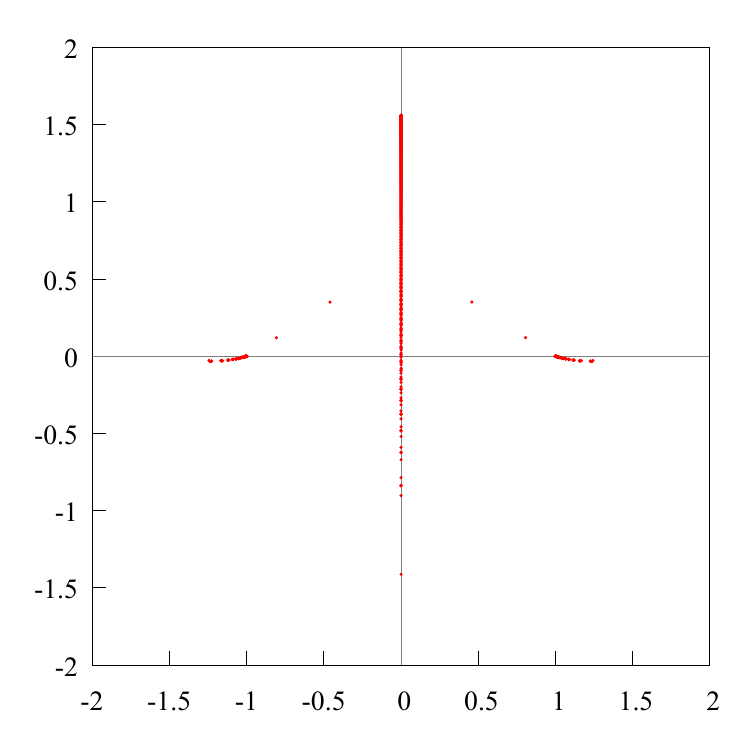}}
\vskip-6mm
\caption{The distribution of the zeros of the Hermite-Pad\'{e} polynomials $Q_{n,1}$ (red points)
$n=166,\dots,170$ for a set of three functions $[1,f,f^2]$, where
$f(z)=(1-z^2)^{1/4}(1-i \sqrt{3}\cdot 0.9z)^{-1/2}$.
}
\label{Fig_nik_(_9)_4000_166-170_rd}
\end{figure}

\newpage
\begin{figure}[!ht]
\centerline{
\includegraphics[width=15cm,height=15cm]{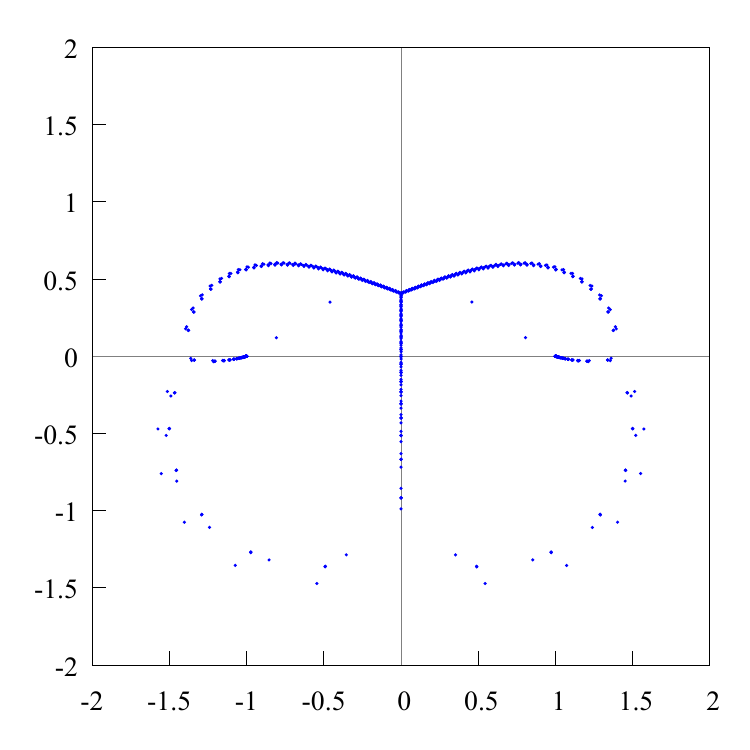}}
\vskip-6mm
\caption{
The distribution of the zeros of the Hermite-Pad\'{e} polynomials $Q_{n,0}$ (blue points),
$n=166,\dots,170$, for a set of three functions $[1,f,f^2]$, where
$f(z)=(1-z^2)^{1/4}(1-i \sqrt{3}\cdot 0.9z)^{-1/2}$.
}
\label{Fig_nik_(_9)_4000_166-170_bl}
\end{figure}

\newpage
\begin{figure}[!ht]
\centerline{
\includegraphics[width=15cm,height=15cm]{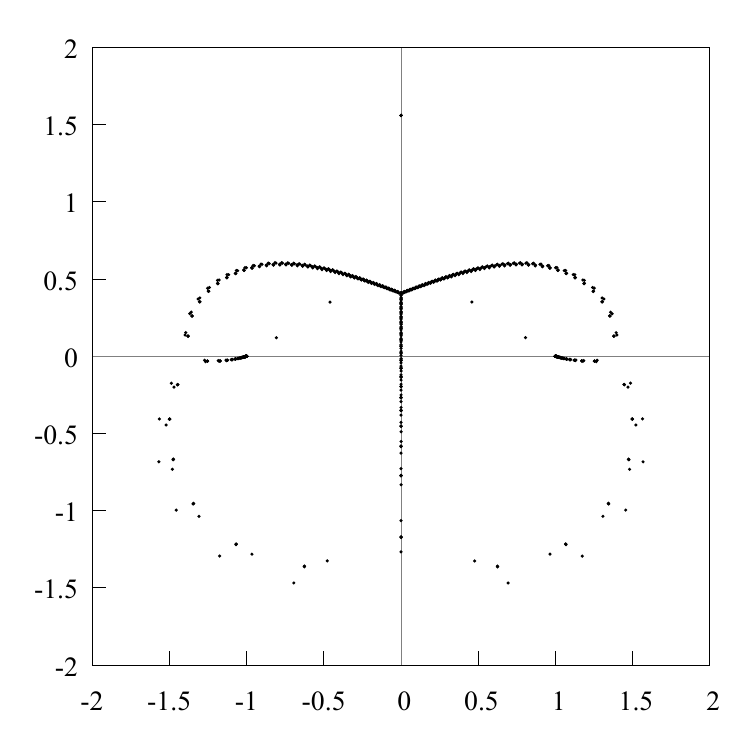}}
\vskip-6mm
\caption{The distribution of the zeros of the Hermite-Pad\'{e} polynomials $Q_{n,2}$ (black points),
$n=166,\dots,170$, for a set of three functions $[1,f,f^2]$, where
$f(z)=(1-z^2)^{1/4}(1-i\sqrt{3}\cdot 0.9z)^{-1/2}$.
There are simple zeros of the polynomials $Q_{n,2}$, $n=166,\dots,170$,
on the positive part of the imaginary axis at the point $z\approx a$, $a=i\sqrt{3}\cdot 0.9$,
corresponding to a simple pole of the function $f^2$ at the point $z=ia$.
}
\label{Fig_nik_(_9)_4000_166-170_bk}
\end{figure}

\newpage

\newpage

\newpage

\newpage

\newpage




\clearpage
\markboth{\bf Buslaev compact set for $f(z)=\sqrt{(z-a)/(z-b)}$}{\bf Buslaev
compact set for $f(z)=\sqrt{(z-a)/(z-b)}$}

\newpage

\newpage

\newpage




\end{document}